\documentclass[11pt, A4paper]{amsart}
\usepackage[latin9]{inputenc}
\usepackage{etoolbox}
\usepackage[margin=2.5cm]{geometry}
\usepackage[english]{babel}
\usepackage[T1]{fontenc}
\usepackage{mathrsfs}
\usepackage{amsthm}
\usepackage{amsmath}
\usepackage{amssymb}
\usepackage{babel}
\usepackage{graphicx}
\usepackage{amsmath}
\usepackage{subfiles}
\usepackage{tikz}
\usepackage{mathdots}
\usepackage{yhmath}
\usepackage{cancel}
\usepackage{color}
\usepackage{siunitx}
\usepackage{array}
\usepackage{multirow}
\usepackage{gensymb}
\usepackage{enumitem}
\usepackage{comment}
\usepackage{hyperref}
 \usepackage{bbm}
 \usepackage[toc]{glossaries}
\hypersetup{
    colorlinks=true,
    linkcolor=blue,
    filecolor=magenta,      
    urlcolor=blue,
    citecolor=red,
}
\usetikzlibrary{fadings}
\newcommand{\mcO}{\mathcal{O}}
\newcommand{\mcF}{\mathcal{F}}
\newcommand{\mcG}{\mathcal{G}}
\newcommand{\mcOt}{\mathcal{O}^{\times}}
\newcommand\mf[1]{\mathfrak{#1}}
\newcommand\wh[1]{\widehat{#1}}
\newcommand\un[1]{\varpi_{\mathfrak{#1}}}
\newcommand\ov[1]{\overline{#1}}
\newcommand{\Cc}{{\mathbb{C}}}

\newcommand{\Ff}{{\mathbb{F}}}

\newcommand{\Zz}{{\mathbb{Z}}}
\newcommand{\eps}{\epsilon}

\newcommand{\Rr}{{\mathbb{R}}}

\newcommand{\Aa}{{\mathbb{A}}}
\newcommand{\Af}{{\mathbb{A}_f}}

\newcommand{\Aft}{{\mathbb{A}_f^\times}}

\newcommand{\Ft}{F^\times}
\newcommand{\Qq}{{\mathbb{Q}}}
\newcommand{\Xx}{{\mathbb{X}}}

\newcommand{\Fp}{\Ff_p}

\newcommand{\Fqn}{\Ff_{q^n}}
\newcommand{\Fq}{\Ff_q}

\newcommand{\Fpt}{\Ff^\times_p}
\newcommand{\Cor}{\mathscr{C}}

\newcommand{\whO}{\widehat{\mcO}}
\newcommand{\whOt}{\widehat{\mcO}^{\times}}
\newcommand{\whK}{\widehat{K}}
\newcommand{\ind}{\mathbbm{1}}

\newcommand{\vphi}{\varphi}

\newcommand{\Pfp}{\mathrm{P}^{1}_{\Fp}}

\newcommand{\Afp}{\mathrm{A}^{1}_{\Fp}}
\newcommand{\Pfq}{\mathrm{P}^{1}_{\Fq}}

\newcommand{\Afq}{\mathrm{A}^{1}_{\Fq}}
\newcommand{\Ks}{\mathrm{K}^{sep}}
\newcommand{\Ga}{\mathrm{G}^{arith}}
\newcommand{\Gg}{\mathrm{G}^{geom}}
\newcommand\norm[1] {\left\lVert#1\right\rVert}
\DeclareMathOperator{\Gal}{Gal}
\DeclareMathOperator{\Kl}{Kl}

\DeclareMathOperator{\Aut}{Aut}

\DeclareMathOperator{\GL}{GL}
\DeclareMathOperator{\Cl}{Cl}

\DeclareMathOperator{\Nm}{Nm}
\DeclareMathOperator{\Tr}{Tr}
\DeclareMathOperator{\disc}{disc}
\DeclareMathOperator{\Spec}{Spec}
\DeclareMathOperator{\cond}{cond}
\DeclareMathOperator{\meas}{meas}

\DeclareMathOperator{\Frob}{Frob}

\DeclareMathOperator{\PGL}{PGL}

\theoremstyle{remark}
\newtheorem*{rmk}{Remark}

\theoremstyle{plain}
\newtheorem{thm}{Theorem}[section]
\numberwithin{thm}{section}
\newtheorem*{thm*}{Theorem}
\newtheorem{mainthm}{Main Theorem}
\newtheorem*{mainthm*}{Main Theorem}
\newtheorem{lemma}[thm]{Lemma}

\newtheorem{prop}[thm]{Proposition}

\newtheorem{obs}[thm]{Observation}

\newtheorem*{cor}{Corollary}

\theoremstyle{definition}
\newtheorem{defn}[thm]{Definition}

\title{Algebraic twists of GL(2) automorphic forms }
\author{Vignesh Arumugam Nadarajan}
\date{\today}
\usepackage[
backend=biber,
style=alphabetic,
sorting=nty
]{biblatex}

\begin{document}
\begin{abstract}
    In this article we consider the problem of estimating the correlation of Hecke eigenvalues of $\GL_2$ automorphic forms with a class of functions of algebraic origin defined over finite fields called trace functions. The class of trace functions is vast and includes many standard exponential sums like Gauss sums, Klostermann sums, Hyperklostermann sums etc. In particular we prove a Burgess type power saving (of $1/8$) over the trivial bound for forms with level coprime to the prime $\mf p$ over whose residue field the trace function is defined. For forms whose level is divisible by at most one power of $\mf p$ we obtain a saving of $1/12$. This generalizes the results of \cite{FKM15} to the case of number fields, with a slightly more restrictive assumption on the Fourier-M\"obius group attached to the trace function. Moreover, the implied constant for the $1/12$ power saving bound has polynomial dependence on the archimedean part of the conductor of the form. We work using the language of adeles which makes the analysis involved softer and makes the generalisation to number fields more natural. The proof proceeds by studying the amplified second moment spectral average of the correlation sum using the relative trace formula. 
\end{abstract}
\maketitle
\section{Introduction}
First, let us consider the problem of algebraic twists of cusp forms. By a cusp form we will mean a non-zero holomorphic cusp form or a non-zero Maass cusp form of weight 0. For a cusp form $f$ we have $f(z+1)=f(z)$, so $f$ admits a Fourier expansion which we will write as $$f(z)=\sum_{n=1}^{\infty}\rho_f(n)n^{\frac{k-1}{2}} e(nz)$$ for a holomorphic form of weight k and $$f(z)=\sum_{n=1}^{\infty}\rho_f(n)|n|^{-\frac{1}{2}} W_{it_f}(4\pi|n|y) e(nx).$$ for a Maass cusp form of weight 0. The normalization is chosen so that the coefficients are absolutely bounded $$|\rho_f(n)|\ll_f 1.$$

Here, $e(nz)=e^{2\pi i nz}$ and $W_{it_f}$ is a Whittaker function.

Let $K:\Zz/p\Zz\xrightarrow{}\Cc$ be a function which we will lift to $\Zz$ and let $$\|K\|_{\infty}:=max_{1\leq n\leq p}|K(n)|=M.$$ We would like to consider correlation sums of the form $$S(f,K,p):=\sum_{1\leq n\leq p}\rho_f(n) K(n).$$ We may consider the smoothed version $$S_V(f,K,p)=\sum_{1\leq n}\rho_f(n)K(n)V(n/p)$$ where $V(.)$ is a smooth compactly supported function. The trivial bound for such a sum is $$S_V(f,K,p)\ll _{V,f} Mp.$$

In their paper \cite{FKM15}, \'E. Fouvry, E. Kowalski and Ph. Michel showed that one can do much better for a special class of functions called isotypic trace functions. They proved
\begin{thm*}[FKM]
Let $f$ be a Hecke eigenform, p a prime number. Let K be an isotypic trace function of conductor cond(K).

There exists $s\geq 1$ absolute constant such that:
$$S_V(f,K;p)\ll cond(K)^s p^{1-\delta}$$ for any $\delta< 1/8$ with an implied constant depending only on V, f and $\delta$.
\end{thm*}{}
This generalizes the classical non-correlation results when $K$ is an additive character or Dirichlet character to a vast class of functions $K$ like many algebraic exponential sums, point counting functions of algebraic varieties over finite fields. The beauty of this approach lies in the uniform treatment of this vast class of functions without explicit structure like multiplicativity.

In this article, we study the generalization of this question to $\GL_2$ automorphic forms over number fields. 

Let $\phi$ be a cuspidal $\GL_2$-automorphic form over a number field $F$, of level $\mf N$ fixed. Let $V:F_{\infty}\xrightarrow[]{} \Rr$ be a compactly supported smooth function. The question we consider is to bound, $$\sum_{m\in\Ft}K(m_{\mf p})W_{\phi,f}\begin{pmatrix}
m\un p & 0 \\
0 & 1 
\end{pmatrix}V(m_{\infty})$$ where $W_{\phi,f}$ is the finite part of the global Whittaker function of $\phi$. For $m\in \frac{1}{\un p} \mcO_{\mf p}$ we denote $m_{\mf p}$ as the congruence class of $m\un p$ mod $\mf p$ with $\un p$ being the uniformizer at $\mf p$.

\subsection{Trivial bound}\label{triv}

Let $K:k(\mf p)\xrightarrow[]{}\Cc$ be a function defined over the residue field at $\mf p$. Let $\phi$ be an automorphic form that is spherical (i.e. $K_{v}$ invariant for all finite places $v$) with $\phi$ a pure tensor (i.e. $W_{\phi}=\prod_{v}W_{\phi, v}$). We have the following trivial bound:

$$|\sum_{m\in\Ft}K(m_{\mf p})W_{\phi}\begin{pmatrix}
m\pi_{\mf p } & 0 \\
0 & 1 
\end{pmatrix}V(m_{\infty})
|\ll_{\phi, F, K, V,\epsilon} \Nm(\mf p)^{\frac{1}{2}+\epsilon}$$

using Rankin-Selberg theory.

We prove the following two theorems:
\begin{mainthm}
Let $\pi$ be a cuspidal $\GL_2$-automorphic representation over $F$ of level $\mf N$ with $\mf N$ co-prime to $\mf p$. Let $\phi$ be an automorphic form in $\pi$ that is of level $\mf N$. Let $K$ be an isotypic trace function defined over the residue field $k(\mf p)$ with an associated Fourier-M\"obius group contained in the standard Borel subgroup. There exists an absolute constant $s>0$ s.t.

$$\sum_{m\in\Ft}K(m_{\mf p})W_{\phi,f}\begin{pmatrix}
m\pi_{\mf p } & 0 \\
0 & 1 
\end{pmatrix}V(m_{\infty})\ll_{\phi, F, V,\delta} \cond(K)^s\Nm(\mf p)^{\frac{1}{2}-\delta}$$

for any $\delta<\frac{1}{8}$.
\end{mainthm}

\begin{mainthm}
Let $\pi$ be a cuspidal $\GL_2$-automorphic representation over $F$ of level $\mf N$ with $0\leq v_{\mf p}(\mf N)\leq 1$. Let $\phi$ be an automorphic form in $\pi$ that is of level $\mf N$. Let $K$ be an isotypic trace function defined over the residue field $k(\mf p)$ with an associated Fourier-M\"obius group contained in the standard Borel subgroup. There exists an absolute constant $s>0$ s.t.

$$\sum_{m\in\Ft}K(m_{\mf p})W_{\phi,f}\begin{pmatrix}
m\pi_{\mf p } & 0 \\
0 & 1 
\end{pmatrix}V(m_{\infty})\ll_{\phi, F, V,\delta} \cond(K)^s\Nm(\mf p)^{\frac{1}{2}-\delta}$$

for any $\delta<\frac{1}{12}$.
\end{mainthm}

\begin{rmk}

Here are some remarks concerning the main theorems:
\begin{enumerate}
    \item By $K$ an isotypic trace function, we mean $K$ is the middle extension trace function of a $\ell$-adic sheaf $\mcF$ lisse on an open set $U\subseteq \Afp$ that is pure of weight $0$ on $U$ and is geometrically isotypic. $\cond(K)$ is an invariant associated with a trace function called its (analytic) conductor. The reader may refer to sections \S\ref{trnot} and \S\ref{trfnsize} for the definition of these terms. For an isotypic trace function $K$, we have
$$\norm{K}_{\infty} \leq \cond(K)$$

(see the remark \ref{purrmk}.)
    \item By the standard Borel subgroup, we mean the subgroup of upper triangular matrices. The Fourier-M\"obius group attached to a trace function is defined in section \S\ref{autgp}.

    \item We can show that $s=4$ would suffice in the above theorems.

    \item In main theorem 1, we can show an exponential dependency of the implied constant in the archimedean part of $\cond(\pi)$. The obstruction to the polynomial dependency in $\cond(\pi)$ is the choice of amplifier that is supported over primes. This choice due to A. Venkatesh (\cite{AVe10}) was made to simplify many details of the proof.

    \item In main theorem 2, we can show a polynomial dependency in the archimedean part of $\cond(\pi)$. 

     \item The normalizer of the diagonal torus occurs as the Fourier-M\"obius group of a Gauss sum constructed from a quadratic Dirichlet character. Looking at the proof of our theorems however, the readers can verify that our proof remains valid for this example.

    \item For an example of a trace function to which our theorem does not apply, the reader may consult \cite{FKM15}, section 11.2(3). The authors consider the symmetric square of a pullback of a Kloostermann sum and show that it has the dihedral subgroup as the Fourier-M\"obius group.    

    \item The applicability of our main theorem follows from the existence of families of trace functions with bounded conductor (and sup-norm), defined over residue fields of an unbounded sequence of primes in $F$ . See \S\ref{trfnexamples} for some examples.

    \item Our main theorems generalize the result of \cite{FKM15} to the case of number fields, with a more restrictive assumption on the trace function. Note that we obtain the same  power saving uniformly for all number fields.
\end{enumerate}    
\end{rmk}
\subsection{Acknowledgements}
The work for theorem 1 was done at EPFL during my PhD and the extension to theorem 2 was carried out at IISc TDC, Challakere when I was an IoE teaching fellow there. I thank both these institutions for providing a supportive environment for research. I would like to thank my PhD supervisor Prof. Philippe Michel, based on whose insights and guidance this work has emerged. 
\section{Notations and preliminaries on automorphic forms}

\subsection{Number fields}
Fix $F/\Qq$ a number field with ring of integers $\mcO_F$ and discriminant $\Delta_F$. Let $n$ be the degree of the number field and $(r,s)$ be its signature. Let $\Lambda_F$ be the complete $\zeta$-function of $F$, it has a simple pole at $s=1$ with residue denoted as $\Lambda_F^*(1)$.

\subsubsection{Local fields}
For $v$ a place of $F$, we call as $F_v$ the completion of $F$ at $v$. We will also write $F_{\mf p}$ for the place corresponding to the prime ideal $\mf p$. We will denote by $\mcO_v$ the ring of integers at the finite place $v$, by $\un v$ the uniformizer and $$\mf m_v:= \un v \mcO_v$$ is the maximal ideal. The field $$k_v:=\mcO_v/\mf m_v$$ is called the residue field at place $v$. We will denote the residue field at $\mf p$ by $k(\mf p)$. For $s\in \Cc$, we define the local zeta function $$\zeta_{F_v}(s)=(1-q_v^{-s})^{-1}\text{ if } v<\infty$$ (here $q_v=|\mcO_v/\un v\mcO_v|$ is the size of the residue field), $$\zeta_{F_v}(s)=\pi^{-s/2}\Gamma(s/2)\text{ if  $v$ is real}$$ and $$\zeta_{F_v}(s)=2(2\pi)^{-s}\Gamma(s) \text{ if  $v$ is complex}.$$  

\subsubsection{Adeles}
We denote the adele ring of $F$ by $\Aa$ and the group of ideles by $\Aa^{\times}$. We denote by $\Aa_f$ the finite adele ring and by $\Aa_f^{\times}$ the finite ideles. $$\whOt\simeq \prod_{v<\infty} \mcOt_v$$ will denote the maximal open compact subgroup of $\Aa_f^{\times}$.

\subsubsection{Additive characters}
We denote by $\psi$ the standard additive character on $F$. Recall that it is defined by $$\psi=\psi_{\Qq}\circ \Tr_{F/\Qq}$$ and it factors as a product of local characters $$\psi=\prod_{v}\psi_v.$$ Recall that the conductor of $\psi_v$ is defined as the largest ideal of $\mcO_v$,  on which $\psi_v$ is trivial. We will denote the conductor as $\cond(\psi_v)=\mf m_v^{d_v}$ and by $$\Delta_F=\prod_{v<\infty}\mf p_v^{d_v}$$  the discriminant of $F$.  

\subsection{Subgroups}
We fix some notations:  . 
$$G=\text{GL}_2$$
$$\overline{G}=Z/G=\text{PGL}_2$$
$$X= Z(\Aa)G(F) \setminus G(\Aa)$$
Let $H$ be a subgroup of $G$, we define $[H]$ by
$$[H]=H(F)/H(\Aa).$$

Let $R$ be a commutative ring , we define some standard subgroups of $G(R)$:

$$B(R)=\left\{\begin{pmatrix}
a & b \\
0 & d 
\end{pmatrix} : a,d\in R^{\times},  b\in R\right\}\text{ , }N(R)=\left\{\begin{pmatrix}
1 & b \\
0 & 1 
\end{pmatrix} :  b\in R\right\},$$

$$A(R)=\left\{\begin{pmatrix}
a & 0 \\
0 & d 
\end{pmatrix} : a,d\in R^{\times}\right\}\text{ , }Z(R)=\left\{\begin{pmatrix}
a & 0 \\
0 & a 
\end{pmatrix} :  a\in R^{\times}\right\}.$$

We also define
$$n(x)=\begin{pmatrix}
1 & x \\
0 & 1 
\end{pmatrix}\text{, }w= \begin{pmatrix}
0 & -1 \\
1 & 0 
\end{pmatrix}\text{ and }a(y)= \begin{pmatrix}
y & 0 \\
0 & 1 
\end{pmatrix}.$$

For any place $v$, $K_v$ is defined to be the maximal compact subgroup of $G(F_v)$ defined by

$$K_v= \begin{cases}
              G(\mcO_{F_v})               \text{ if $v$ is finite} \\
              O_2(\Rr)              \text{ if $v$ is real} \\
              U_2(\Cc)               \text{ if $v$ is complex}
           \end{cases}$$

For $v<\infty$ and $n\geq 0$, we define

$$K_{v,0}(\un v^n)=\left\{\begin{pmatrix}
a & b \\
c & d 
\end{pmatrix}\in K_v : c\in \mf m_v^n\right\}.$$

For $\mf a= \prod_{v<\infty}\mf p_v^{f_v(\mf a)}$, define

$$K_0(\mf a)=\prod_{v<\infty}K_{v,0}(\un v^{f_v(\mf a)})= \left\{\begin{pmatrix}
a & b \\
c & d 
\end{pmatrix}\in G(\whO) : c\in \mf a\whO \right\}.$$
\subsection{Measures}
The normalization of measures is as in \cite{MV10}.
For every place $v$, $dx$ be the unique self-dual Haar measure on $F_v$. This gives the measure on $N(F_v)$ by identifying it with $F_v$. For $v<\infty$, $dx_v$ gives the measure $q_v^{-d_v/2}$ to $\mcO_{F_v}$. Recall that $d_v$ is the valuation of the discriminant at the place $v$.

  For the multiplicative group $F_v^{\times}$, the Haar measure that we use is $$d^{\times}x_v=\zeta_{F_v}(1)\frac{dx_v}{|x_v|}.$$ By identification we also get the measures on $A(F_v)$ and $Z(F_v)$. We define $dx=\prod_v dx_v$ as the measure on $\Aa$ and $d^{\times}x=\prod_v d^{\times}x_v$ as the measure on $\Aa^{\times}$. 
  
  Using Iwasawa decomposition, $G(F_v)=Z(F_v)N(F_v)A(F_V)K_v$, we define the Haar measure on $G(F_v)$ by
  $$dg_v=d^{\times}zdxd^{\times}y dk.$$ $dk$ is the Haar measure on the compact group $K_v$ normalised to be a probability measure.
  The measure on the adelic points of the above subgroups is the product of the corresponding local measures. By $dg$ we denote also the quotient measure on $X$.
\subsection{Automorphic forms and representation}

Let $L^2(X)$ denote the Hilbert space of square integrable functions w.r.t. the measure defined above. Define

$$L^2_{\text{cusp}}(X)=\{\phi\in L^2(X): \int_{F\setminus \Aa}\phi(n(x)g)dx=0\text{ a.e. }g\}.$$

For $\phi\in L^2_{\text{cusp}}(X)$ we have the Fourier-Whittaker expansion,

$$\phi(g)=\sum_{m\in F^{\times}}W_{\phi}\left (\begin{pmatrix}
m & 0 \\
0 & 1 
\end{pmatrix}g\right )$$ with

$$W_{\phi}\left (g\right )=\int_{F\setminus \Aa}\phi(n(x)g)\psi(-x)dx.$$

$G(\Aa)$ acts by right translation on $L^2(X)$ defining the right regular representation which defines a unitary representation. There is a well-known decomposition into $G(\Aa)$ submodules (and these are orthogonal):

$$L^2(X)=L^2_{\text{cusp}}(X)\oplus L^2_{\text{Res}}(X) \oplus L^2_{\text{Eis}}(X).$$

In this article we will be concerned only with $L^2_{\text{cusp}}(X)$, which decomposes into irreducible unitary subrepresentations of $G(\Aa)$. It is also well known that for every irreducible subrepresentation $\pi$ there is a decomposition into a tensor product of irreducible, unitary, local representations $\pi_v$ (Flath's theorem) $$\pi=\otimes_v \pi_v.$$

\subsection{Whittaker model and factorisation of the inner product}

Let $\pi=\otimes_v \pi_v$ be a unitary representation as in the previous section. The intertwining map

$$\phi\in \pi \longmapsto W_{\phi}:=\int_{F\setminus \Aa}\phi(n(x)g)\psi(-x)dx$$

defines an equivariant embedding to the space of smooth functions. The image is called the Whittaker model of $\pi$ and is denoted $W(\pi,\psi)$. Note that functions in the image satisfy $$W(n(x)g)=\psi(x)W(g).$$ 

It can be shown that this space is isomorphic as $G(\Aa)$ modules to $\otimes_v W(\pi_v,\psi_v)$. $W(\pi_v,\psi_v)$ is the image of the intertwining map 

$$\phi_v\in \pi_v \longmapsto W_{\phi_v}:=\int_{F_v}\phi(w n(x)g)\psi(-x)dx$$

to the space of smooth functions on $G(F_v)$.

Further a pure tensor $\otimes_v \phi_v\in \otimes_v \pi_v$ maps to the pure tensor $\otimes_v W_{\phi_v}\in \otimes_v W(\pi_v,\psi_v)$.

Let us define an inner product on $W(\pi_v,\psi_v)$ by

$$\theta_v(W_v, W_v'):=\zeta_{F_v}(2)\frac{\int_{F_v^{\times}}W_v(a(y))\ov{W_v'}(a(y))d^{\times}y}{\zeta_{F_v}(1)L(\pi_v,Ad,1)}.$$

The normalization is chosen so that $\theta_v(W_v, W_v)=1$ if $\pi_v$, $\psi_v$ are unramified and $W_v(1)=1$ [see \cite{JS81}, prop. 2.3]. We define an invariant inner product $\langle,\rangle_v$ on $\pi_v$ so that the equivariant isomorphism $\pi_v\xrightarrow[]{} W(\pi_v, \psi_v)$ becomes an isometry. 

For $\pi=\bigotimes\pi_v$ a cuspidal representation, we can compare the global norm to the local norms of a pure tensor: Let $\phi=\bigotimes \phi_v\in \pi$ be a pure tensor,

$$||\phi||^2=2\Delta_F^{1/2}\Lambda^*(\pi, \text{Ad}, 1)
\prod_v\langle \phi_v,\phi_v\rangle_v$$

where $\Lambda(\pi, \text{Ad}, s)$ is the complete adjoint L-function and $\Lambda^*(\pi, \text{Ad}, 1)$ is the first non-vanishing coefficient in its Laurent expansion at $s=1$.

\subsection{Hecke operators}\label{Hecke}

Let $\mf l$ be a prime and $(\pi_{\mf l},V)$ be a representation of $\ov{G}(F_{\mf l})$. Further let us assume that $\pi_{\mf l}$ is unramified i.e. the set of $K_{\mf l}:=\ov{G}(\mcO_{\mf l})$ fixed vectors $$V^{K_{\mf l}}\neq \emptyset.$$ The elements of $V^{K_{\mf l}}$ are also called spherical vectors. 

It can be shown that the subspace of spherical vectors is one dimensional. Consider the function $h_{\mf l}: \ov{G}(F_{\mf l})\xrightarrow[]{} \Cc$ defined as
$$h_{\mf l}=\frac{1}{\sqrt{\Nm(\mf l)}}\ind_{K_{\mf l}\begin{pmatrix}
\un l & 0 \\
0 & 1 
\end{pmatrix}K_{\mf l}}=\frac{1}{\sqrt{\Nm(\mf l)}}\left(\sum_{x\in k(\mf l)}\ind_{\begin{pmatrix}
\un l & x \\
0 & 1 
\end{pmatrix}K_{\mf l}}+\ind_{\begin{pmatrix}
1 & 0 \\
0 & \un l 
\end{pmatrix}K_{\mf l}}\right )$$

This function is bi-$K_{\mf l}$ invariant and therefore the operator $\pi_{\mf l}(h_{\mf l})$ leaves the subspace of spherical vectors invariant. This operator is called the \textit{Hecke operator}. Since the subspace of spherical vectors is one dimensional, the Hecke operator is represented by a scalar called the \textit{Hecke eigenvalue}.

It is known that any infinite dimensional irreducible subrepresentation of the right regular representation on $L^2(\ov{G}(F_{\mf l})$ that is unramified is isomorphic to the so-called unramified principal series representation $\pi(\chi,\ov{\chi})$ where $\chi:F_{\mf l}^{\times}\xrightarrow[]{} \Cc^{\times}$ is an unramified quasicharacter. The parameters $$\alpha:=\chi(\un l) \text{ and } \beta=\ov{\chi(\un l)}$$ are called the \textit{Satake parameters}. It can be easily checked that the Hecke eigenvalue in this case is given by $$\lambda=\alpha+\beta.$$
\subsection{Adelic Poisson summation formula}
\begin{defn}[Schwartz-Bruhat function]
    A function $f:\Aa\xrightarrow[]{} \Cc$ is said to be a Schwartz-Bruhat function if it can be written as a finite linear combination of factorizable functions $f=\prod_v f_v$ where:

\begin{enumerate}
    \item For $v|\infty$, $f_v$ is a Schwartz function.
    \item For $v< \infty$ $f_v$ is a locally constant, compactly supported function.
\end{enumerate}

\end{defn}
\subsubsection{Adelic Fourier transform}
For $f:\Aa\xrightarrow[]{} \Cc$ a Schwartz-Bruhat function, we define the adelic Fourier transform by

$$\wh{f}(x)=\int_{\Aa} f(y)\psi(-xy) dy$$

where $\psi$ is the standard additive character defined above and the measure is the additive Haar measure normalized as above.

We have the Fourier inversion formula:

$$\wh{\wh{f}}(x)=f(-x).$$

and the adelic Poisson summation formula
 For $f:\Aa\xrightarrow[]{} \Cc$ a Schwartz-Bruhat function, we have 

$$\sum_{a\in F}f(a)=\sum_{a\in F}\wh{f}(a).$$

Note that $F$ is a discrete subset of $\Aa$.

\subsection{Some results on counting lattice points}

We will need the following bounds for future estimates.

\begin{lemma}[Lattice point counting - dense lattice limit] \label{lpc}
Let $I\lhd F$ be a non-zero fractional ideal and $f:F_{\infty}\xrightarrow{}\Rr$ be a Schwartz function. For $\Nm(I)\leq 1$ and $R\in\Rr_{>0}$ a scaling parameter, we have the following estimate:

$$\left|\sum_{m\in I}f\left(\frac{m}{R}\right)\right|\ll_{F}||f||_{1,n+1}\frac{R^n}{\Nm(I)}$$
where by $f(m)$ we mean $f(m_{\infty})$ and $||f||_{1,n+1}$ is the Sobolev norm controlling the $L^1$ norm of the first $n+1$ derivatives. Recall that $n$ is the degree of the number field $F$.
\end{lemma}
\begin{rmk}
\smallskip
    \begin{enumerate}
        \item One can give a more precise asymptotic formula which involves the length of the shortest vector in the lattice. In the case of ideal lattices, arithmetic-geometric mean inequality in \ref{agm} controls the length of the shortest vector in terms of the covolume of the lattice.
        
        \item We are interested in the limit $|\Nm(I)|\xrightarrow[]{}0$ i.e. the ideal lattice becoming dense. This motivates applying Poisson summation to obtain a sum over the dual lattice which we will see becomes sparse in this limit, giving the term $m=0$ as the main contribution.
    \end{enumerate}
\end{rmk}
\proof

In this case the Poisson summation formula reads, 

$$\sum_{m\in I}f\left(\frac{m}{R}\right)=\frac{R^n}{|\Nm(I)|}\sum_{m\in I^*}\hat{f}(mR)$$

where $$\hat f(y):=\int_{F_{\infty}}f(z)e(-\Tr_{F/\Qq}(yz))d\mu(z)$$ with $d\mu$ being the self dual Haar measure for the trace pairing form. 

It can be checked that $I^*$ the dual lattice under the trace pairing form is also a fractional ideal. For any fractional ideal lattice the minimum is controlled in terms of the covolume. Let $v\in I^*\setminus{0}$, we have $v=(\sigma_i(r))_{i=1}^n$ with $\sigma_i:F\xrightarrow{}F_{\infty}$ the geometric embeddings. The $L^1$ norm of $v$ satisfies
$$|v|=\sum_{i=1}^n|\sigma_i(r)|\geq n\left(\prod_{i=1}^n|\sigma_i(r)|\right)^{\frac{1}{n}}=n\left(\Nm(r)\right)^{\frac{1}{n}}.$$

Since $r\in I$, $\Nm(I^*)|\Nm(r)$, so in conclusion 
\begin{equation*}\label{agm}
|v|=\sum_{i=1}^n|\sigma_i(r)|\geq n\left|\Nm(I^*)\right|^{\frac{1}{n}}.
\end{equation*}
We deduce the following bounds:
$$\sum_{m\in I}f(m)=\frac{R^n}{|\Nm(I)|}\hat{f}(0)+\frac{R^n}{|\Nm(I)|}\sum_{m\in I^*\setminus{0}}\hat{f}(mR).$$

Since $\hat{f}$ is a Schwartz function, for any $N>n$ we can show that
$$\sum_{m\in I^*\setminus{0}}|\hat{f}(mR)|\ll_{N} \frac{||f^{(N)}||_1}{(R\min(L))^N} \leq \frac{||f^{(N)}||_1}{R^N n\left|\Nm(I^*)\right|^{\frac{N}{n}}}.$$

By choosing $N=n+1$ and noting that $\Nm(I^*)=1/(\disc(F).\Nm(I))$ we get

$$\sum_{m\in I}f(m)=\frac{R^n}{|\Nm(I)|}\big(\hat{f}(0)+o_{F}(||f^{(n+1)}||_1)\big)$$

completing the proof. \qed

We will also need to sum Schwartz functions over ideal lattices whose covolume is tending to infinity. As we saw above the minimum of the lattice in this case also tends to infinity.

\begin{lemma}[Lattice point counting - sparse lattice limit] \label{lpcs}
Let $I\lhd F$ be a non-zero fractional ideal and $f:F_{\infty}\xrightarrow{}\Rr$ be a Schwartz function. For $\Nm(I)\geq 1$ and $R\in\Rr_{>0}$ a scaling parameter, we have the following estimate:

$$\left|\sum_{m\in I}f\left(\frac{m}{R}\right)\right|\ll_{F}||f||_{1,n+1}\frac{R^n}{\Nm(I)}$$
where by $f(m)$ we mean $f(m_{\infty})$ and $||f||_{1,n+1}$ is the Sobolev norm controlling the $L^1$ norm of the first $n+1$ derivatives. Recall that $n$ is the degree of the number field $F$.
\end{lemma}

\begin{lemma}[Counting units]
\label{unic}
Let $m_0\in F$ and $I=(m_0)$ be the corresponding principal fractional ideal. Let $R>0$ be a positive real. We have the following bound:

$$\sum_{\substack{(m)=I\\|m|_{\infty}<R}}1\ll_F \log(R)^{n-1}.$$

\end{lemma}

\proof

We have
$$\sum_{\substack{(m)=I\\|m|_{\infty}<R}}1=\sum_{\substack{u\in \mcOt\\|m_0u|_{\infty}<R}}1.$$

The archimedean embeddings of $F$ are denoted by $$\sigma_1,...,\sigma_{r},\sigma_{r+1},...,\sigma_{r+s},...,\sigma_{r+2s}.$$ Here $\sigma_i$ $1\leq i\leq r$ are real and the remaining are complex non-real embeddings which satisfy $\sigma_{r+j}=\ov{\sigma_{r+s+j}}$. Also $r+2s=n$ is the degree of the number field.

With this notation the $\text{Log}$ map is defined as
$\text{Log}:F^{\times}\xrightarrow[]{}\Rr^{r+s}$ 
$$\text{Log}(m)=(\log(|\sigma_1(m)|),...,\log(|\sigma_{r}(m)|),\log(|\sigma_{r+1}(m)|),...,\log(|\sigma_{r+s}(m)|)).$$

By the Dirichlet structure theorem on units,
$\text{Log}$ restricted to $\mcOt$ has a finite kernel -- the roots of unity in $F$ and the image $\text{Log}(\mcOt):=\Lambda$ is a lattice contained in the hyperplane $\{(x_1,...,x_{r+s})\in \Rr^{r+s}:\sum x_i=0\}$.

Therefore applying the Log map we have the following equality

$$\sum_{\substack{u\in \mcOt\\|m_0u|_{\infty}<R}}1=|\mu_F|\sum_{\substack{v\in v_0+\Lambda\\|v|<\log(R)}}1$$

where $\mu_F$ is the roots of unity contained in $F$ and $v_0=\text{Log}(m_0)\in \Rr^{r+s}$. By the lattice counting lemma -- lemma \ref{lpc}, approximating the indicator function of the translated ball by a smooth cutoff function:

$$\sum_{\substack{v\in \Lambda}}\ind_{B(\log(R))}(v_0+v)\ll_F \log(R)^{n-1}$$ completing the proof.

\qed

\begin{lemma}[Counting divisors]
\label{divc}
Let $I$ be a fractional ideal and $k\in I^2$. Let $R>0$ be a positive real. We wish to count the number of divisors of $k$ which are in $I$ of size bounded by $R$:

$$\sum_{\substack{m,n\in I\\|m|_{\infty}<R\\|n|_{\infty}<R}}\ind_{mn=k}\ll_F (\Nm(k))^{\text{o}(1)}.R^{\text{o}(1)}.$$

\end{lemma}

\proof

We have 
$$m\in I \iff (m)\subseteq I$$
and
$$n\in I \iff (k)I^{-1}\subseteq (m).$$

Therefore,

$$    \sum_{\substack{m,n\in I\\|m|_{\infty}<R\\|n|_{\infty}<R}}\ind_{mn=k}=\sum_{\substack{\frac{|k|_{\infty}}{R}\leq |m|_{\infty}<R}}\ind_{(k)I^{-1}\subseteq(m)\subseteq I} $$

\begin{multline*}
    \sum_{\substack{\frac{|k|_{\infty}}{R}\leq |m|_{\infty}<R}}\ind_{(k)I^{-1}\subseteq(m)\subseteq I} \leq  \sum_{\substack{ |m|_{\infty}<R}}\ind_{(k)I^{-1}\subseteq(m)\subseteq I}\\ \leq  \sum_{\substack{\mf J\lhd \mcO\\ \Nm(\mf J)< \frac{R^n}{\Nm(I)}}}\ind_{(k)I^{-2}\subseteq \mf J} \ind_{[\mf J. I]=1}   \sum_{\substack{ |m|_{\infty}<R}}\ind_{(m)=\mf J}.
\end{multline*}

By the unit counting lemma (lemma \ref{unic}) we have 

$$\sum_{\substack{ |m|_{\infty}<R}}\ind_{(m)=\mf J}\ll_F R^{\text{o}(1)}.$$

We need to estimate the sum

$$\sum_{\substack{\mf J\lhd \mcO\\ \Nm(\mf J)< \frac{R^n}{\Nm(I)}}}\ind_{(k)I^{-2}\subseteq \mf J} \ind_{[\mf J. I]=1}.$$

To do this let us define a divisor counting function. Let $\mf a\lhd \mcO$ be an ideal,

$$
d(\mf a,X):= \sum_{\substack{\mf J\lhd \mcO\\ \Nm(\mf J)< X}} \ind_{\mf a \subseteq \mf J} =\sum_{\substack{d|\Nm(\mf a)\\ d<X}}  \sum_{\substack{\mf J\lhd \mcO\\ \Nm(\mf J)=d}} \ind_{\mf a \subseteq \mf J}.
$$

It follows that 
$$
d(\mf a,X)= \sum_{\substack{d|\Nm(\mf a)\\ d<X}}  \sum_{\substack{\mf J\lhd \mcO\\ \Nm(\mf J)=d}} \ind_{\mf a \subseteq \mf J} \leq \sum_{\substack{d|\Nm(\mf a)\\ d<X}} d^{\text{o}(1)} \leq \tau(\Nm(\mf a)) X^{\text{o}(1)}.
$$
The first inequality follows since

$$ \sum_{\substack{\mf J\lhd \mcO\\ \Nm(\mf J)=d}} \ind_{\mf a \subseteq \mf J} \leq  \sum_{\substack{\mf J\lhd \mcO\\ \Nm(\mf J)=d}} 1 \leq \tau(d)^{n} \leq d^{\text{o}(1)}. $$

where the last inequality is a standard fact in analytic number theory.

So we deduce

$$\sum_{\substack{\mf J\lhd \mcO\\ \Nm(\mf J)< \frac{R^n}{\Nm(I)}}}\ind_{(k)I^{-2}\subseteq \mf J} \ind_{[\mf J. I]=1}\leq d((k)I^{-2}, \frac{R^n}{\Nm(I)})\leq (\Nm(k) R)^{\text{o}(1)}.$$

Putting this together with the count of units, we conclude the proof.
\qed
\section{Preliminaries on trace functions}

In the next 3 sections we will collect some basic material about trace functions which will be useful to us. We follow the treatment in \cite{FKM15} and \cite{FKMS19}). Fix $p\neq \ell$ two prime numbers, $q=p^n$ and an isomorphism (of fields) $i:\overline{\Qq_{\ell}}\xrightarrow{}\Cc$.

\begin{defn}(\textbf{Trace sheaf})
A constructible $\ell$-adic sheaf $\mathcal{F}$ on $\mathbb{A}^1_{\mathbb{F}_q}$ is called a trace sheaf if it is a middle-extension sheaf whose restriction to any non-empty open subset $U$ where $\mathcal{F}$ is lisse, $\mathcal{F}|_U$ is pointwise $\iota$-pure of weight 0.

\end{defn}

The category of constructible $\ell$-adic sheaves on $\Pfq$ can be described in terms of the category of $\ell$-adic representations of the group $$\Gal(\overline{\Fq(X)}^{sep}|\Fq(X)).$$ Since this point of view may be easier to understand we will unpack our definitions a bit using the terminology of Galois representations. Let us state some results from \cite{FKMS19} that we will need.

Let $\Pfq$ denote the projective line on $\Fq$ and $\Afq$ the affine line. For this section $\text{K}=\Fq(X)$ is the field of functions. We will denote by $\Ks$ the separable closure of $\text{K}$ and $\Ga=\Gal(\Ks|\text{K})$ the arithmetic Galois group. Moreover $\Gg=\Gal(\Ks|\ov{\Fq}\text{K})\subseteq \Ga$ is the geometric Galois group, where $\ov{\Fq}\subseteq \Ks$ is an algebraic closure of $\Fq$. We have the exact sequence:

$$1\xrightarrow[]{} \Gg \xrightarrow[]{} \Ga \xrightarrow[]{} \Gal(\ov{\Fq
}|\Fq)\xrightarrow[]{} 1$$

\begin{defn}
Let $U\subseteq \Afq$ be an open set. A sheaf $\mcF$ is lisse on $U$ means it is associated to a continuous finite dimensional representation

$$\rho_{\mcF}: \Ga \xrightarrow[]{} \GL(V_{\mcF})$$

which is unramified at every closed point $x\in U$. The dimension of $V_{\mcF}$ is denoted $\text{rk}(\mcF)$.
\end{defn}

Here unramified means the action of any inertia subgroup at $x$ on $V_{\mcF}$ is trivial. For details the reader may consult \cite{FKMS19}.
\begin{defn}(\textbf{Trace function})
Given a sheaf $\mcF$ lisse on an open set $U\subseteq \Afq$ the associated trace function is defined by

$$x\in U(\Fq) \longmapsto K_{\mcF}(x)=\Tr(\Frob_x|\mcF_x)$$

$\Frob_x$ is the Frobenius element of $x$. It is only defined up to conjugation but the trace is conjugacy invariant. The trace is an element of $\ov{\Qq_l}$ but we take its image under $i$.
\end{defn}{}

For our purpose we will use the so-called middle extension, this means we will view $K_{\mcF}$ as a function on all of $\Fq$ by taking the trace of the Frobenius at every point on the subspace invariant under the inertia group at the point.

\subsection{Some notions and constructions}\label{trnot}
We borrow terminology from representation theory to describe the corresponding properties of $\ell$-adic sheafs:

A $\ell$-adic sheaf is irreducible (resp. isotypic) if the corresponding representation is.

A $\ell$-adic sheaf is geometrically irreducible (resp. geometrically isotypic) if the restriction to $\Gg$ of the corresponding representation is. 

Next we want to discuss Fourier transform of trace functions. For this discussion it will be important to restrict ourselves a little:

\begin{defn}[Fourier Sheaf]
An isotypic lisse sheaf is said to be Fourier iff none of the geometric irreducible components are geometrically isomorphic to the Artin-Schreier sheaf.
\end{defn}
The Artin-Schreier sheaf is an irreducible lisse sheaf on $\Afq$ whose trace function is the additive character on $\Fq$.
For a function $f:\Fq\xrightarrow[]{} \Cc$ we recall its (unitarily normalized) Fourier transform

$$\wh{f}(x)=\frac{1}{\sqrt{q}}\sum_{y\in\Fq}f(y)e(\frac{-xy}{q}).$$

Deligne proved that the Fourier transform of trace functions arise as the trace function of a 'geometric Fourier transform' of the corresponding sheaf. Here is the precise statement due to Laumon:

\begin{thm}
Let $\mcF$ be a Fourier sheaf lisse on an open set $U$ and pure of weight $0$. Let us denote by $K_{\mcF}$ the middle-extension trace function of $\mcF$. There exists a Fourier sheaf $\wh{\mcF}$ lisse on an open set $\wh{U}$ pure, of weight $0$ s.t. its middle extension trace function $K_{\wh{\mcF},n}$ is the Fourier transform $\wh{K_{\mcF,n}}$

$$\wh{K_{\mcF,n}}=\frac{1}{\sqrt{q^n}}\sum_{y\in\mathbb{F}_{q^n}}f(y)e(\frac{-\Tr_{\mathbb{F}_{q^n}/\Fq}xy}{q}).$$

The map $\mcF\mapsto \wh{\mcF}$ is called the geometric Fourier transform. 
\end{thm}

See below for the definition of pure of weight $0$. Moreover the conductor of the geometric Fourier transform can be shown to be bounded in terms of the conductor of the original sheaf. We also have the proposition:

\begin{prop}
The geometric Fourier transform maps irreducible (resp. isotypic) Fourier sheaves to irreducible (resp. isotypic) Fourier sheaves.
\end{prop}

See proposition 7.8 in \cite{FKMS19}.

Let us look at some examples of trace functions.

 \subsection{Examples of trace functions}\label{trfnexamples}
 For more details on these examples consult section 4.2 of \cite{FKMS19}.
\begin{enumerate}
    \item \textbf{Trivial sheaf}
    The $\ell$-adic sheaf corresponding to the trivial representaion $\ov{\Qq_{\ell}}$ is lisse everywhere, pure of weight 0, rank 1 and has the trace function
    
    $$K_{\mcF,n}\equiv 1.$$
    
    \item \textbf{Kummer sheaf}
    Let $\chi:\Fpt\xrightarrow[]{}\Cc$ be a Dirichlet character. There exists a $\ell$-adic sheaf called the Kummer sheaf that is lisse on $\Pfp\setminus \{0,\infty\}$, pure of weight 0, rank 1 and has the trace function
    
    $$K_{\mcF}(x)=\chi(x)\text{  }x\in\Fpt$$
    $$K_{\mcF,n}(x)=\chi(\Nm_{\mathbb{F}_{p^n}|\Fp} (x))\text{  }x\in\mathbb{F}_{p^n}^{\times}.$$
    
    \item \textbf{Artin-Schreier sheaf}
    Let $\psi:\Fp\xrightarrow[]{}\Cc$ be a non-trivial additive character. There exists a $\ell$-adic sheaf called the Artin-Schreier sheaf that is lisse on $\Pfp\setminus \{\infty\}$, pure of weight 0, rank 1 and has the trace function
    
    $$K_{\mcF}(x)=\psi(x)\text{  }x\in\Fp$$
    $$K_{\mcF,n}(x)=\psi(\Tr_{\mathbb{F}_{p^n}|\Fp} (x))\text{  }x\in\mathbb{F}_{p^n}.$$
    
    \item \textbf{Hyperkloostermann sheaf}
    
      Let $\psi:\Fp\xrightarrow[]{}\Cc$ be a non-trivial additive character. Define for every $m\geq 2$ the hyperkloostermann sum as the m-fold multiplicative convolution of the $\psi$ 
      
      $$Kl_{m,\psi}(a;p)=\frac{1}{p^{m-1/2}}\sum_{x_1x_2...x_m=a}\psi(\frac{x_1+...+x_m}{p})$$ where $x_1,...,x_m\in\Fpt$.

      $$Kl_{m,\psi}(a;p^n)=\frac{1}{p^{n.(m-1)/2}}\sum_{x_1x_2...x_m=a}\psi(\frac{x_1+...+x_m}{p^n})$$ where $x_1,...,x_m\in\mathbb{F}_{p^n}^{\times}
      $. 
     
     Deligne showed that $\forall$ $m\geq 2$ there exists a $\ell$-adic sheaf $\mathcal{K}\ell_{m,\psi}$ lisse on $\Pfp\setminus \{0,\infty\}$, pure of weight 0, rank n, geometrically irreducible and has the trace function
     
     $$K_{\mathcal{K}\ell_{m,\psi},n}(x)=Kl_{m,\psi}(x;p^n)$$
     for $n\geq 1$, $x \in\mathbb{F}_{p^n}^{\times}
      $.
 \item The class of trace functions is closed under taking sums, products, under precomposition with rational functions. Also with a mild restriction this class is closed under taking Fourier transforms. 
     \end{enumerate}{}

\subsection{Size of trace functions}\label{trfnsize}

The two quantities attached to a trace function that is of importance are its sup-norm and a measure of its complexity called its conductor. We do not discuss the conductor in detail but refer the reader to \cite{FKMS19}. We will also need the notion of the conductor of a trace function $K$ (denoted $\cond(K)$) -- this is the smallest conductor among all $\ell$-adic sheaves whose trace function is $K$.
\begin{defn}[Conductor]
    Let $\mcF$ be a a $\ell$-adic sheaf lisse on an open set $U$, the conductor of $\mcF$ is defined as

    $$C(\mcF)=\text{rk}(\mcF)+|\text{P}^1\setminus U(\ov{\Fq})|+\sum_{x\in D_{ram}(\ov{\Fq})}\text{Swan}_x(\mcF).$$

    Here $D_{ram}$ is the ramification locus of $\mcF$ and $\text{Swan}_x(\mcF)$ is the Swan conductor of $\mcF$ at $x$.
\end{defn}
We refer the reader to section 3.7 of \cite{FKMS19} for the definition of the latter two quantities.
\begin{defn}[Purity]
A $\ell$-adic sheaf $\mcF$ lisse on an open set $U\subseteq \Afq$ is said to be $i-$punctually pure of weight $w$, if for every $x\in U$, the eigenvalues of $\Frob_x$ on $V_{\mcF}$ (applying the map $i$) are complex numbers of absolute value less than or equal to $q^{w/2}$. In particular for every $x\in U$,

$$|K_{\mcF}(x)|\leq \text{rk}(\mcF)q^{w/2}.$$
\end{defn}

\begin{rmk}\label{purrmk}
\begin{enumerate}
    \item If we start a  sheaf $\mcF$ lisse on an open set $U\subseteq \Afq$, $i-$punctually pure of weight $w$ . We can always reduce to a sheaf $i-$punctually pure of weight 0 lisse on $U$ by making an appropriate twist. See remark 3.11 in \cite{FKMS19}.
    
     \item If we start a sheaf $\mcF$ lisse on an open set $U\subseteq \Afq$ $i-$punctually pure of weight $w$. The bound 
     
     $$|K_{\mcF}(x)|\leq \text{rk}(\mcF)q^{w/2}.$$
     
     was shown to hold for all $x\in\Afq$ for the middle extension by Deligne. See remark 3.12 in \cite{FKMS19}.
\end{enumerate}
\end{rmk}

Let $\mcF$ be a $\ell$-adic sheaf lisse on an open set $U$ pure, of weight $0$. In this article we will encounter sums of the shape

$$\sum_{x\in\Fq}K_{\mcF}(x)$$

The trivial bound for such a sum is $(\text{rk}(\mcF))q$. Thanks to Deligne's proof of the generalized Riemann hypothesis we can do better:

\begin{thm}[Square root cancellation for trace functions]
Let $\mcF$ be a $\ell$-adic sheaf lisse on an open set $U$ pure, of weight $0$ that is geometrically irreducible or isotypic with no trivial components. We have

$$\sum_{x\in\Fq}K_{\mcF}(x)\ll C(\mcF)^2q^{1/2}$$ 

where $C(\mcF)$ is the conductor of $\mcF$ and the implied constants in the above bound are absolute.
\end{thm}

Restricting ourselves to geometrically isotypic sheaves is not a serious restriction, see proposition 5.1 in \cite{FKMS19} for more details.

 \subsection{Correlation sums and quasiorthogonality relation} \label{corrsum}
 
 Let $\mcF$ and $\mcG$ be two sheaves that are lisse on an open set $U$, pure of weight $0$. We can apply the results of the previous section to the sheaf $\mcF\otimes D(\mcG)$ which is also lisse on $U$, pure of weight $0$. Note that $$K_{\mcF\otimes D(\mcG)}=K_{\mcF}\ov{K_{\mcG}}.$$
 
 We will be interested in calculating correlation sums of the form, 
 
 $$\mathcal{C}(\mcF,\mcG)=\frac{1}{q}\sum_{a\in\mathbb{F}_{q}}K_{\mcF}(a)\ov{K_{\mcG}(a)}.$$
 
 Applying the theorem from the previous section we get
 
 \begin{thm}[Quasi-orthogonality relation]
 Let $\mcF$ and $\mcG$ be as above and let us further assume that they are geometrically irreducible. There exists a complex number $\alpha_{\mcF,\mcG}\in S^1$ s.t.
 
  $$\mathcal{C}(\mcF,\mcG)=\alpha_{\mcF,\mcG} \delta_{\mcF \simeq_{geom}\mcG}+O(\frac{C(\mcF)^2C(\mcG)^2}{q^{1/2}}).$$
 
 \end{thm}
 
 Here we have used the fact that $C(\mcF\otimes D(\mcG))\leq C(\mcF)C(\mcG)$ and the results of the previous section.
 \medskip
 
 Let $\mcF$ be a sheaf that is lisse on an open set $U$, pure of weight $0$ and geometrically irreducible. Assume further that $\mcF$ is Fourier. Let $\gamma\in \PGL_2(\mathbb{F}_{q})$ acting as a fractional linear transformation on $\Pfq$,
 
 $$\gamma=\begin{pmatrix}
a & b \\
c & d 
\end{pmatrix}:x\mapsto \frac{ax+b}{cx+d}$$

In this article we will encounter the following kind of correlation sums:

$$\mathcal{C}(\mcF,\gamma):=\mathcal{C}(\wh{\mcF},\gamma*\wh{\mcF})=\frac{1}{q}\sum_{a\in\mathbb{F}_{q}}\wh{K_{\mcF}}(a)\ov{\wh{K_{\mcF}}(\gamma(a))}.$$

On applying the above theorem, we see that the size of the above sum is controlled by the following invariant:

\begin{thm}
    [Automorphism group of the sheaf]
Let $\mcF$ be as above. There exists an algebraic group $\Aut_{\mcF}\subseteq \PGL_2$ defined over $\Fq$, the  group of automorphisms of $\mcF$. It satisfies $\Aut_{\mcF}(\Fqn)\subseteq \PGL_2(\Fqn)$ consist of all $\gamma$ s.t.

$$\gamma^*\mcF\simeq_{\text{geom}}\mcF$$

 Let $B$ denote the subgroup of upper triangular matrices, we define
$$B_{\mcF}:=B\cap \Aut_{\mcF}.$$
\end{thm}

\begin{defn} [Fourier-M\"obius group]\label{autgp}
    Let $\mcF$ be as above. The Fourier-M\"obius group is defined to be the group of automorphisms of $\wh{\mcF}$, the geometric Fourier transform of $\mcF$.
\end{defn}

 \begin{cor} \label{corrsumsize}
 If $\gamma$ is not an element of the Fourier-M\"obius group, we have
 
 $$\mathcal{C}(\mcF,\gamma)=O_{C(\mcF)}(q^{-1/2}).$$
 \end{cor}
 
 \section{The relative trace formula}
\subsection{Introduction}
The main tool that we will use to compute and bound the spectral average is the relative trace formula. This approach relies on the fact that the correlation sum we are interested in bounding is a trace function(algebraic) twist of a unipotent period. In this chapter we will outline the generalities related to relative trace formula, we will apply this formula to our situation in the next chapter.

We can extend the right regular action of $G$ on $L^2([G])$ to an action of smooth compactly supported functions on $L^2([G])$. So for a smooth compactly supported function $f$ we associate an operator $R(f)$ on $L^2([G])$. The relative trace formula computes the relative trace of the operator $R(f)$ i.e. integral of the kernel $K(x,y)$ of $R(f)$ over $[H]$ for a subgroup $H<G\times G$. 

The formula relates the spectral side i.e. the period computed using the spectral expansion of the kernel to the so called geometric side -- a sum over certain integrals over orbits of $H$. This allows one to express spectral average of period integrals over automorphic forms as a sum of relative orbital integrals. We refer the reader to chapter 18 of \cite{GH23} for some history and generalities of trace formulae. 

Our starting point are the works of A. Knightly and C. Li \cite{KL06}, \cite{KL13} where the authors derive the Kuznetsov trace formula (of weight 0) by taking $G=\PGL_2$ and $H=N\times N$. We replace the Hecke operator by a Hecke operator twisted by a trace function in order to express the second moment of the correlation sum as a sum of relative orbital integrals and in turn as a sum of Klostermann sums. Given the Bruhat decomposition, the structure of the orbits is also easy to determine. We will spell out the setup of the relative trace formula in our case in the beginning of the next chapter.
\subsection{The kernel function}
\begin{defn}
We say $f:G(\Aa)\xrightarrow[]{}\Cc$ is \textit{m times smooth and factorizable} if it is a product of local functions

$$f=f_{\infty}\prod_pf_p$$ satisfying

    \begin{enumerate}
        \item  $f$ is invariant under $Z(\Aa)$.
        \item  $f$ has compact support modulo the center.
        \item  $f_{\infty}$ is m-times continuously differentiable.
        \item  $f_{p}$ is locally constant.
        \item  For a.e. $p$, $f_p=\mathbbm{1}_{Z(F_p)K_p}$.
    \end{enumerate}

\end{defn}

Note that such an $f$ is in $L^1(\overline{G}(\Aa))$. We will denote the set of functions satisfying the above property by $C_c^m(G(\Aa))$ following \cite{KL13} and by  $C_c^{\infty}(G(\Aa))$ the set of smooth factorizable functions.

Recall, that for $\phi\in L^2([\overline{G}])$ and $f\in C^{m}_c(G(\Aa))$, we have the induced right action
    
    $$R(f)\phi(x)=\int_{\overline{G}(\Aa)}f(y)\phi(xy)dy.$$
    
    By unfolding, we get
    $$R(f)\phi(x)=\int_{[\overline{G}]}\sum_{\gamma\in \overline{G}(F)} f(x^{-1}\gamma y)\phi(y)dy.$$
    The operator $R(f)$ therefore has the kernel 
    $$K(x,y)=\sum_{\gamma\in \overline{G}(F)} f(x^{-1}\gamma y).$$
    
 Note that if  $f\in C^{m}_c(G(\Aa))$ the sum is locally finite and is therefore continuous.
 \begin{thm*}[Thm 6.1 in \cite{KL13}]
 For $m\geq 3$ and $f\in C^{m}_c(G(\Aa))$ the operator $R(f)$ on $L^2([\overline{G}])$ is Hilbert-Schmidt, which in this case means that $K\in L^2([\overline{G}]\times[\overline{G}])$.
 \end{thm*}
\subsection{The spectral side}
Assume $f\in C^{\infty}_c(G(\Aa))$, we have (see \cite{KL13} equation 6.11)
$$K(x,y)=K_{cusp}(x,y)+K_{res}(x,y)+K_{Eis}(x,y)$$
for \textbf{a.e. $x, y\in [\overline{G}]$}
where

$$K_{cusp}(x,y)=\sum_{\varphi}R(f)\varphi(x)\overline{\varphi(y)}$$
the sum running over an orthonormal basis of the cuspidal part of the spectrum.

$$K_{res}(x,y)=\frac{3}{\pi}\sum_{\chi^2=1}R(f)\phi_{\chi}(x)\overline{\phi_{\chi}(y)}$$
the sum running over an orthonormal basis of the residual part of the spectrum.
$$K_{Eis}(x,y)=\frac{1}{4\pi}\sum_{\psi\in H(0)}\int_{-\infty}^{\infty}E(R(f)\psi_{it},x)\overline{E(\psi_{it},y)}dt$$
the sum running over an orthonormal basis of $H(0):=\bigoplus_{\chi}H(\chi,\chi^{-1},0).$

Like in \cite{KL13} we wish to integrate this expression over $(N(\Aa)/N(F))^2$ which has measure $0$ in $\Xx^2$. So we use the following result:

\begin{thm*}\label{spectexp}
 Assume $f=h_1*h_2^{*}$ for $h_1,h_2\in C^{\infty}_c(G(\Aa))$ where $h_2^*(y)=\overline{h_2(y^{-1})}$, we have
 \begin{equation}
     K(x,y)=K_{cusp}(x,y)+K_{res}(x,y)+K_{Eis}(x,y).\label{1}
 \end{equation}

for \textbf{every} $x, y\in [\overline{G}]$
where

$$K_{cusp}(x,y)=\sum_{\varphi}R(h_1)\varphi(x)\overline{R(h_2)\varphi(y)}$$
the sum running over an orthonormal basis of the cuspidal part of the spectrum.

$$K_{res}(x,y)=\frac{3}{\pi}\sum_{\chi^2=1}R(h_1)\phi_{\chi}(x)\overline{R(h_2)\phi_{\chi}(y)}$$
the sum running over an orthonormal basis of the residual part of the spectrum.
$$K_{Eis}(x,y)=\frac{1}{4\pi}\sum_{\psi\in H(0)}\int_{-\infty}^{\infty}E(R(h_1)\psi_{it},x)\overline{E(R(h_2)\psi_{it},y)}dt$$
the sum running over an orthonormal basis of $H(0):=\bigoplus_{\chi}H(\chi,\chi^{-1},0).$

All the sums and integrals above converge absolutely.
\end{thm*}

In \cite{KL13}, the authors prove an everywhere, absolute convergence statement assuming that the archimedean part is  binvariant under the maximal compact subgroup. They do this by using the fact that such a function $f$ can be written as $h_1*h_2+k_1*k_2$ for some $h_1,h_2,k_1,k_2$ satisfying the same properties as $f$. So, their proof can be adapted to prove the above theorem. 

\begin{rmk}
It turns out that $$K_{res}(x,y)\equiv\frac{3}{\pi}\int_{\overline{G}(\Aa)}f(g)dg.$$
\end{rmk}

 We are interested in computing 
    $$I=\int\limits_{[N]\times [N]}K(x,y)\psi_m(x)\ov{\psi_n(y)}dxdy$$
    
    where $\psi_m(n(x))=\psi(mx)$, $\psi$ is the standard additive character on $N(F)/N(\Aa)\simeq F/\Aa$ and $m,n\in F^{\times}$. 
    
    The spectral side is obtained by computing $I$ using the spectral expansion \eqref{1}.
 \subsection{The geometric side}   
    We will now compute the integral using the geometric expansion of the kernel
     \begin{align*}     
    I&=\int\limits_{[N]\times [N]}K(x,y)\psi_m(x)\ov{\psi_n(y)}dxdy\\
     &=\int\limits_{[N]\times [N]}\sum_{\gamma\in \overline{G}(F)}f(x^{-1}\gamma y)\psi_m(x)\ov{\psi_n(y)}dxdy
     \end{align*}
     
      Breaking the sum into orbits and unfolding as in \cite{KL06}, we get
     $$I=\sum_{[\delta]}I_{\delta}$$
     where the sum runs over orbits of $\ov{G}(F)$ under conjugation action of $N(F)$ i.e. 
     $$[\delta]=\{x^{-1}\delta y|x,y\in N(F)\}\text{ for }\delta\in \overline{G}(F).$$

       $$I_{\delta}=\int\limits_{H_{\delta}(F)/H(\Aa)}f(x^{-1}\delta y)\psi_m(x)\ov{\psi_n(y)}dxdy,$$
    
      $$H=N\times N\text{ and     }H_{\delta}=\{(x,y)\in H|x^{-1}\delta y=\delta\}\text{ for }\delta\in\overline{G}.$$

     Observe that $I_{\delta}=0$ unless $\psi_m(x)\ov{\psi_n(y)}=1$ for all $(x,y)\in H_{\delta}(\Aa)$. 
     \begin{defn}
     We will say $[\delta]$ is \textbf{relevant}, if $\psi_m(x)\ov{\psi_n(y)}=1$ for $(x,y)\in H_{\delta}(\Aa).$ 
     \end{defn}
     We will now calculate $I_{\delta}$ for relevant orbits.
\subsection{A parametrization of the orbits}
     
Bruhat decomposition for $G$ is as follows:

$$G=NA\coprod NAwN$$ 
where $N$ is the unipotent radical, $A$ the diagonal torus and $w=\begin{pmatrix}
0 & 1 \\
1 & 0 
\end{pmatrix}$ is the Weyl element.

We get therefore the following description of the orbits:

$$N(F)\backslash \overline{G}(F) / N(F)=\bigg\{\bigg [\begin{pmatrix}
\lambda & 0 \\
0 & 1 
\end{pmatrix}\bigg ]\bigg |\lambda\in \Ft \bigg \}\coprod \bigg\{\bigg [\begin{pmatrix}
0 & \mu \\
1 & 0 
\end{pmatrix}\bigg ]\bigg |\mu\in \Ft \bigg \}.$$

     For $$\delta= \begin{pmatrix}
\lambda & 0 \\
0 & 1 
\end{pmatrix} $$

we calculate the stabilizer explicitly:

recall $$H_{\delta}=\{(x,y)\in H|x^{-1}\delta y=\delta\}$$ we find

$$H_{\delta}(\Aa)=\Big\{\Big(\begin{pmatrix}
1 & \lambda t \\
0 & 1 
\end{pmatrix},\begin{pmatrix}
1 & t \\
0 & 1 
\end{pmatrix}\Big)|t\in\Aa\Big\}.$$

So $\delta$ is relevant only if $\psi((n-m\lambda)t)\equiv 1$ for all $t\in\Aa$, i.e.$$n=m\lambda.$$ We will call the relative orbital integral corresponding to this orbit as the \textit{diagonal contribution} to the geometric side.

 The next case is $$\delta= \begin{pmatrix}
0 & \mu \\
1 & 0 
\end{pmatrix} .$$

By explicit calculation, the stabilizer is found to be

$$H_{\delta}(\Aa)=\{(e,e)\}.$$

So all $\delta$'s are relevant.  We will call the sum of the relative orbital integrals corresponding to these orbits as the \textit{non-diagonal contribution} to the geometric side.
\subsection{Summary of the relative trace formula for the unipotent subgroup}\label{summ}

Putting together the information in the previous two sections, the geometric side reads

 \begin{align*}
     I &=\int\limits_{[N]\times [N]}K(x,y)\psi_m(x)\ov{\psi_n(y)}dxdy \\ &= \sum_{\lambda\in \Ft} \ind_{n=m\lambda}\int\limits_{\{t(\lambda ,1)|t\in F\}\backslash \Aa \times\Aa}f(\begin{pmatrix}
1 & -t_1 \\
0 & 1 
\end{pmatrix} \begin{pmatrix}
\lambda & 0 \\
0 & 1 
\end{pmatrix} \begin{pmatrix}
1 & t_2 \\
0 & 1 
\end{pmatrix})\psi(mt_1-nt_2)dt_1dt_2
\\
&+ \sum_{\mu\in \Ft} \int\limits_{ \Aa \times\Aa}f(\begin{pmatrix}
1 & -t_1 \\
0 & 1 
\end{pmatrix} \begin{pmatrix}
0 & \mu \\
1 & 0 
\end{pmatrix} \begin{pmatrix}
1 & t_2 \\
0 & 1 
\end{pmatrix})\psi(mt_1-nt_2)dt_1dt_2.
\end{align*}

Multiplying out the matrices and making a change of variables we get,

 \begin{multline*}
     I = \sum_{\lambda\in \Ft} \ind_{n=m\lambda}\int\limits_{\{0\}\times F\backslash \Aa \times\Aa}f( \begin{pmatrix}
\lambda & x \\
0 & 1 
\end{pmatrix} )\psi(-mx)dxdy
\\
+ \sum_{\mu\in \Ft} \int\limits_{ \Aa\times \Aa}f(\begin{pmatrix}
-t_1 & \mu-t_1t_2 \\
1 & t_2 
\end{pmatrix})\psi(mt_1)\psi(nt_2)dt_1dt_2.
\end{multline*}

\begin{multline*}
     I = \sum_{\lambda\in \Ft} \ind_{n=m\lambda} \meas(F\backslash \Aa)\int\limits_{ \Aa }f( \begin{pmatrix}
\lambda & x \\
0 & 1 
\end{pmatrix} )\psi(-mx)dx
\\ + \sum_{\mu\in \Ft} \int\limits_{ \Aa\times \Aa}f(\begin{pmatrix}
-t_1 & \mu-t_1t_2 \\
1 & t_2 
\end{pmatrix})\psi(mt_1)\psi(nt_2)dt_1dt_2.
\end{multline*}

This is the form we will use to compute the geometric side explicitly for our choice of $f$. We will also explain the change of variables we used to obtain this form when we make these computations later. 

The spectral side reads for $m,n\in \Ft$

 $$I =\int\limits_{[N]\times [N]}K(x,y)\psi_m(x)\ov{\psi_n(y)}dxdy = I_{cusp}+I_{Eis}$$

 where
 $$I_{cusp} =\int\limits_{[N]\times [N]}K_{cusp}(x,y)\psi_m(x)\ov{\psi_n(y)}dxdy $$
 and 
  $$I_{Eis} =\int\limits_{[N]\times [N]}K_{Eis}(x,y)\psi_m(x)\ov{\psi_n(y)}dxdy .$$

The residual contribution is zero since the residual part of the kernel is a constant as remarked earlier. More explicitly we have

 $$I_{cusp} =\sum_{\pi}\sum_{\varphi\in \mathscr{B}(\pi)} W_{R(h)\varphi}\big(\begin{pmatrix}
m & 0 \\
0 & 1 
\end{pmatrix}\big)\overline{W_{R(h)\varphi}\big(\begin{pmatrix}
n & 0 \\
0 & 1 
\end{pmatrix}\big)}$$

where $\pi$ runs over cuspidal subrepresentations in $L^2([G])$, $\mathscr{B}(\pi)$ is any choice of orthonormal basis of $\pi$ and  $W_{\varphi}$ is the Whittaker function of $\varphi$ as explained in the previous chapter.

The Eisenstein contribution can also be made explicit in terms of Fourier coefficients of Eisenstein series. We will not pursue this here since it will not be needed in the sequel. For more details on Fourier expansion of weight $0$ Eisenstein series, the reader may consult page 43, section 5.6 of \cite{KL13}. 
\section{Strategy and proof of the main theorem}

We begin this chapter by recalling the theorem that we will prove

\begin{mainthm*}
Let $\pi$ be a cuspidal $\GL_2$-automorphic representation over $F$ of level $\mf N$ an integral ideal. Let $\phi$ be an automorphic form in $\pi$ that is of level $\mf N$. Let $K$ be an isotypic trace function defined over the residue field $k(\mf p)$ with an associated Fourier-M\"obius group contained in the standard Borel subgroup. There exists an absolute constant $s>0$ s.t.

$$\sum_{m\in\Ft}K(m_{\mf p})W_{\phi,f}\begin{pmatrix}
m\pi_{\mf p } & 0 \\
0 & 1 
\end{pmatrix}V(m_{\infty})\ll_{\phi, F, V,\delta} \cond(K)^s\Nm(\mf p)^{\frac{1}{2}-\delta}$$

for any $\delta<\frac{1}{8}$ if $\mf N$ is coprime to $\mf p$ and for any $\delta<\frac{1}{12}$ if $0\leq v_{\mf p}(\mf N)\leq 1$ . (The implied constants being different in both the cases.)
\end{mainthm*}

We fix the following in our problem: Let $\mf p$ be a prime in $F$ and $k(\mf p)$ be the residue field (of size $\Nm(\mf p)$). 

 We fix a set $$\Lambda=\{\mf l\in \Spec(\mcO_F)\text{ | }L\leq \Nm(\mf l)\leq 2L\}$$ and $$x_{\mf l}\in\Cc\text{ for }\mf l\in\Lambda. $$ We assume that the elements of $\Lambda$ are coprime to $\mf p$ and $\mf N$. 
 
 Since we are interested in the situation when $L$ and $\Nm(\mf p)$ are very large, we will assume $\cond(\psi)$ has valuation 0 at places corresponding to  $ \mf p$ and to all $\mf l\in \Lambda$.

Let $K:k(\mf p)\xrightarrow[]{}\Cc$ be a function on the residue field at $\mf p$ and $\widehat{K}$ be its Fourier transform, defined by $$\whK(a)=\frac{1}{\sqrt{\Nm(\mf p)}}\sum_{x\in k(\mf p)}K(x)\psi_{\mf p} (-ax/\un p)$$ for $a\in k(\mf p)$. Note the abuse of notation: by $\psi_{\mf p} (-ax/\un p)$ we mean $\psi_{\mf p} (-a_0x_0/\un p)$ for any $a, x\in \mcO_{\mf p}$ whose residue is $a_0, x_0\in k(\mf p)$ respectively. Here we are using the relation between the standard local additive character and the standard additive character of the residue field. 
Let us define the operator $R(f)$ as follows:

 We now fix $$f=\sum_{\substack{\mf l_1,\mf l_2\in \Lambda\\ \mf l_1 \neq \mf l_2}} x_{\mf l_1}\ov{x_{\mf l_2}}h[\mf l_1]*h[\mf l_2]^*$$ where $$h[\mf l]=h_{\infty}[\mf l]\prod_{\mf q}h_{\mf q}[\mf l]$$ is defined such that
 
   \medskip
   
    \begin{enumerate}\label{defn1}
\item For the prime $\mf q=\mf p$
$$h_\mf q[\mf l]:=\frac{1}{\meas(K_0(\mf p))\sqrt{\Nm(\mf p)}}\sum_{x\in k(\mf p)}\widehat{K}(x)\ind_{Z(F_\mf q)\delta_\mf q\begin{pmatrix}
\un q & x \\
0 & 1 
\end{pmatrix}K_0(\mf p)_{\mf q}}$$ 

where $K_0(\mf p)=\{M\in G(\whO)| M_{21}\equiv 0 \mod \mf p\whO\}$.

\item For primes $\mf q=\mf l$ $$h_{\mf q}[\mf l]=\frac{1}{\sqrt{\Nm(\mf l)}}\ind_{Z(F_\mf q)\{M\in M_2(\mcO_\mf q)\text{ | }\text{det}(M)\in \mf l\mcO_\mf q^*\}}$$ 

\item  For primes $\mf q\nmid \mf l\mf p$ $$h_\mf q[\mf l]=\ind_{Z(F_\mf q) K_0(\mf N)_\mf q}$$ where $K_0(\mf N)=\{M\in G(\whO)| M_{21}\equiv 0 \mod \mf N\whO\}$. 

\item  We will define the archimedean components subsequently.
 \end{enumerate}
Compare with the definition in section \S \ref{Hecke}.
 We will also apply the relative trace formula to another choice of test function $$f=h[1]*h[1]^*$$ where $h[1]$ is defined such that

  \medskip
   
    \begin{enumerate}\label{defn2}
\item For the prime $\mf q=\mf p$
$$h_\mf q[1]:=\frac{1}{\meas(K_0(\mf p))\sqrt{\Nm(\mf p)}}\sum_{x\in k(\mf p)}\widehat{K}(x)\ind_{Z(F_\mf q)\delta_\mf q\begin{pmatrix}
\un q & x \\
0 & 1 
\end{pmatrix}K_0(\mf p)_{\mf q}}$$ 

where $K_0(\mf p)=\{M\in G(\whO)| M_{21}\equiv 0 \mod \mf p\whO\}$.

\item  For primes $\mf q\neq \mf p$ $$h_\mf q[1]=\ind_{Z(F_\mf q) K_0(\mf N)_\mf q}$$ where $K_0(\mf N)=\{M\in G(\whO)| M_{21}\equiv 0 \mod \mf N\whO\}$.

\item  We will define the archimedean components subsequently.
 \end{enumerate}
 \subsection{The spectral side of our relative trace formula}
 As we have mentioned in the introduction, we are interested in computing the amplified spectral average of the correlation sum we are interested in using the relative trace formula. To this end we set up some notations: 

 $$I_{m,n}:=\int\limits_{[N]\times [N]}K(x,y)\psi_m(x)\ov{\psi_n(y)}dxdy$$
    
    where $\psi_m(n(x))=\psi(mx)$, $\psi$ is the standard additive character on $$N(F)\setminus N(\Aa)\simeq F\setminus \Aa$$ and $m,n\in F$.  Here $K$ is the kernel of the twisted Hecke operator $R(f)$ we defined below. We will be using the same framework in turn for $$f=\sum_{\substack{\mf l_1,\mf l_2\in \Lambda\\ \mf l_1 \neq \mf l_2}} x_{\mf l_1}\ov{x_{\mf l_2}}h[\mf l_1]*h[\mf l_2]^*$$ and $$f=h[1]*h[1]^*.$$

Following the discussion in section \S\ref{summ} we can write,

$$I_{m,n}=I_{m,n,\text{cusp}}+I_{m,n,\text{Eis}}$$

where 

 $$I_{m,n,\text{cusp}}:=\int\limits_{[N]\times [N]}K_{\text{cusp}}(x,y)\psi_m(x)\ov{\psi_n(y)}dxdy$$

 and 
 $$I_{m,n,\text{Eis}}:=\int\limits_{[N]\times [N]}K_{\text{Eis}}(x,y)\psi_m(x)\ov{\psi_n(y)}dxdy$$
    
    Recall that the residual contribution is $0$. We will then sum $I_{m,n}$ over all $m,n\in \Ft$.

Define $$I:=\sum_{m,n\in \Ft}I_{m,n}.$$
Likewise we may define $I_{\text{cusp}}$ and $I_{\text{Eis}}$. Note that using the expression

$$I_{m,n,\text{cusp}}=\sum_{\pi}\sum_{\varphi\in \mathscr{B}(\pi)} W_{R(h)\varphi}\begin{pmatrix}
m & 0 \\
0 & 1 
\end{pmatrix}\overline{W_{R(h)\varphi}\begin{pmatrix}
n & 0 \\
0 & 1 
\end{pmatrix}}$$

discussed in section \S\ref{summ}, we have the following expression 

$$I_{\text{cusp}}=\sum_{\pi}\sum_{\varphi\in \mathscr{B}(\pi)} \sum_{m\in \Ft}\left|W_{R(h)\varphi}\begin{pmatrix}
m & 0 \\
0 & 1 
\end{pmatrix}\right|^2.$$

This allows us to deduce the non-negativity of the cuspidal contribution and likewise for the Eisenstein contribution. The non-negativity of individual summands in the cuspidal contribution will allow us to upper bound them by the whole spectral side and hence by the geometric side.

Here is a road map for what follows: 
\begin{itemize}
    \item section \S\ref{spectsidell} computes the $I_{m,n, \text{cusp}}$ for our two choices of $f$ to show that it indeed is the amplified spectral average we want
    \item section \S\ref{nonarchsupp} makes precise the non-archimedean support of the spectral average.
    \item section \S\ref{archrestat} formulates a restatement of the problem modifying the archimedean aspect.
    \item section \S\ref{proofmt} details the proof of the two main theorems assuming certain estimates for the geometric side.
    \item chapter \S\ref{geomside} computes the geometric side for the two choices of $f$ and proves the estimates mentioned in the previous item.
\end{itemize}

\subsubsection{What does the spectral side look like?} \label{spectsidell}
  Let us calculate the cuspidal contribution to the spectral side for our two choices of $f$. From section \S\ref{summ} we have
  $$I_{m,n, \text{cusp}}=\sum_{\pi}\sum_{\varphi\in \mathscr{B}(\pi)} W_{R(h)\varphi}\begin{pmatrix}
m & 0 \\
0 & 1 
\end{pmatrix}\overline{W_{R(h)\varphi}\begin{pmatrix}
n & 0 \\
0 & 1 
\end{pmatrix}}$$
where $\mathscr{B}(\pi)$ is an orthonormal basis of an irreducible cuspidal representation $\pi$. Let us choose the orthonormal basis to be made of pure tensors. 

Let us look at the contribution of one cuspidal irreducible representation. Let $\pi$ be an irreducible cuspidal representation, let $\pi=\bigotimes_v\pi_v$ (as explained in the preliminaries). Let $\vphi\in\pi$ be a pure tensor. Assume $\vphi_f$, the non-archimedean part is right $K_0(\mf N)$ invariant. 

As we saw in the preliminaries we can factorise the Whittaker function of a pure tensor automorphic form

Using the definition of $h_f[\mf l]$, a quick calculation shows:
\begin{align*}
W_{R(h_f[\mf l])\vphi_{f}}\begin{pmatrix}
m & 0 \\
0 & 1 
\end{pmatrix} &=\lambda_{\pi}(\mf l)\prod_{\substack{v<\infty\\v\neq \mf p}}W_{\vphi_{v}}\begin{pmatrix}
m & 0 \\
0 & 1 
\end{pmatrix}\frac{1}{\sqrt{\Nm(\mf p)}}\sum_{x\in k(\mf p)}\widehat{K}(x)W_{\vphi_{\mf p}}\begin{pmatrix}
m\un p & mx \\
0 & 1 
\end{pmatrix}\\
&=\lambda_{\pi}(\mf l)\prod_{\substack{v<\infty\\v\neq \mf p}}W_{\vphi_{v}}\begin{pmatrix}
m & 0 \\
0 & 1 
\end{pmatrix}W_{\vphi_{\mf p}}\begin{pmatrix}
m\un p & 0 \\
0 & 1 
\end{pmatrix}\frac{1}{\sqrt{\Nm(\mf p)}}\sum_{x\in k(\mf p)}\widehat{K}(x)\psi(mx).
\end{align*}

Note that since the conductor of $\pi$ is coprime to $\mf p$  $$W_{\vphi_{\mf p}}\begin{pmatrix}
m\un p & 0 \\
0 & 1 
\end{pmatrix}\neq 0$$ iff $$m\un p\in \mcO_p.$$
For $m\in \frac{1}{\un p}\mcO_p$,

$$W_{R(h_f[\mf l])\vphi_{f}}\begin{pmatrix}
m & 0 \\
0 & 1 
\end{pmatrix}=\lambda_{\pi}(\mf l)W_{\vphi_{f}}\begin{pmatrix}
m\un p & 0 \\
0 & 1 
\end{pmatrix}K(m_p)$$

where $m_p$ is the residue of $m\un p$. 

Likewise for $m\in \frac{1}{\un p}\mcO_p$,

$$W_{R(h_f[1])\vphi_{f}}\begin{pmatrix}
m & 0 \\
0 & 1 
\end{pmatrix}=W_{\vphi_{f}}\begin{pmatrix}
m\un p & 0 \\
0 & 1 
\end{pmatrix}K(m_p)$$

where $m_p$ is the residue of $m\un p$.

This computation is the first ingredient needed to carry out our proof strategy in section \S\ref{proofmt}.  We will examine in the next section the support of the finite part of $\cond(\pi)$.

 \subsubsection{Non-archimedean support of the spectral side}
  \label{nonarchsupp}
The non-archimedean support of the spectral side is controlled by the invariance of $h$. We have the following proposition:

\begin{prop}\label{invarianceprop} 
Let $K_a$ and $K_b$ be open, compact subgroups of $G(F_{\mf q})$ and $f\in L^1(\overline{G}(F_{\mf q}))$ which is left $K_a$-invariant and right $K_b$-invariant. Consider a unitary representation $(V, R)$ of $G(F_{\mf q})$. We have $$R(f)V\subseteq V^{K_a}$$ and $$R(f)((V^{K_b})^\perp)=\{0\}$$
\end{prop}

\proof

Let $\phi$ be a vector in a unitary representation $(V, R)$ of $G(F_{\mf q})$. 

$$R(f)\phi:=\int\limits_{\overline{G}(F_{\mf p})}f(h)R(h)\phi dh$$

Note that the above integral is defined to satisfy for any $\psi \in V$,

$$\langle R(f)\phi,\psi\rangle =\int\limits_{\overline{G}(F_{\mf p})}f(h)\langle R(h)\phi, \psi\rangle dh.$$
For any $k\in K_a$

\begin{align*}\langle R(k)R(f)\phi,\psi\rangle &= \langle R(f)\phi, R(k^{-1})\psi\rangle =\int\limits_{\overline{G}(F_{\mf p})}f(h)\langle R(h)\phi, R(k^{-1})\psi\rangle dh\\ &=\int\limits_{\overline{G}(F_{\mf p})}f(h)\langle R(kh)\phi, \psi\rangle dh=\int\limits_{\overline{G}(F_{\mf p})}f(k^{-1}h)\langle R(h)\phi, \psi\rangle dh \\ &=\langle R(f)\phi,\psi\rangle .
\end{align*}
In other words $R(f)\phi$ is $K_a$ invariant. Here we have used that $\overline{G}$ is unimodular and $f$ is left $K_a$ invariant.

We have that $$R(f)^{*}=R(f^*)$$ i.e. the adjoint of $R(f)$ is the right translation by $f^*$ which is defined by $f^*(g)=\overline{f(g^{-1})}$.  (Note that $f^{*}$ is left $K_b$ invariant.)

For any $\phi \in {V^{K_b}}^{\perp}$ and any $\psi\in V$: 

$$\langle R(f)\phi,\psi\rangle =\langle \phi,R(f)^{*}\psi\rangle =\langle \phi,R(f^{*})\psi\rangle =0.$$
Here we have used that $R(f^{*})\psi$ is $K_b$ invariant. Since $\psi$ is arbitrary we conclude $R(f)\phi=0$.
\qed

Let us examine the invariance of $h_{\mf p}[l]$ next since that is the non-trivial one:

\subsubsection{Invariance of $h_{\mf p}[\mf l]$}
We have that $h_{\mf p}[\mf l]=h_{\mf p}[1]$ is right invariant under the open compact $K_0(\mf p)$ and left invariant under the open compact $K_1(\mf p)$.: 

$$h_{\mf p}[\mf l]=\frac{1}{\meas(K_0(\mf p))\sqrt{\Nm(\mf p)}}\sum_{x\in k(\mf p)}\widehat{K}(x)\ind_{Z(F_{\mf p})\begin{pmatrix}
\un p & x \\
0 & 1 
\end{pmatrix}K_0(\mf p)}.$$
    
$h_{\mf p}[\mf l]$ is clearly right $K_0({\mf p})$-invariant. Let us investigate the left action. Let $$g=\begin{pmatrix}
a & b \\
c & d 
\end{pmatrix}\in K_{\mf p}.$$

$$^g h_{\mf p}[\mf l](u):=h_{\mf p}[\mf l](gu).$$

Let $$u\in Z(F_\mf p)\begin{pmatrix}
\un p & x \\
0 & 1 
\end{pmatrix}K_0(\mf p).$$
Note that $$Z(F_{\mf p})\begin{pmatrix}
\un p & x \\
0 & 1 
\end{pmatrix}K_0(\mf p)$$ for different $x\in k(\mf p)$ are disjoint. Let us investigate when

$$gu\in Z(F_\mf p)\begin{pmatrix}
\un p & y \\
0 & 1 
\end{pmatrix}K_0(\mf p).$$
This holds iff
$$g\begin{pmatrix}
\un p & x \\
0 & 1 
\end{pmatrix}\in \begin{pmatrix}
\un p & y \\
0 & 1 
\end{pmatrix}K_0(\mf p)\iff \begin{pmatrix}
\un p & y \\
0 & 1 
\end{pmatrix}^{-1}g\begin{pmatrix}
\un p & x \\
0 & 1 
\end{pmatrix}\in K_0(\mf p).$$

  By explicit calculation $$ \begin{pmatrix}
\un p & y \\
0 & 1 
\end{pmatrix}^{-1}g\begin{pmatrix}
\un p & x \\
0 & 1 
\end{pmatrix}=\begin{pmatrix}
a-cy & \frac{x(a-cy)+b-dy}{\un p} \\
c\un p & cx+d 
\end{pmatrix}.$$   

If $g\in (Id+\un p K_{\mf p})$, $$\begin{pmatrix}
a-cy & \frac{x(a-cy)+b-dy}{\un p} \\
c\un p & cx+d 
\end{pmatrix}\in K_0(\mf p)$$ for $x=y$ and therefore $$gu\in Z(F_{\mf p})\begin{pmatrix}
\un p & x \\
0 & 1 
\end{pmatrix}K_0(\mf p).$$
It follows that$$^g h_{\mf p}[\mf l](u)=h_{\mf p}[\mf l](u)$$ for $g\in (Id+\un pK_{\mf p})=:K_1(\mf p).$
\qed

\begin{cor}\label{invcor}
$R(h_{\mf p}[\mf l])=R(h_{\mf p}[1])$ is non-zero only on local representations $\pi_{\mf p}$ of conductor at most $\mf p$. 
\end{cor}
The corollary follows using proposition \ref{invarianceprop} by noting that for $f=h_{\mf p}[\mf l]$ one has $K_a=K_1(\mf p)$ and $K_b=K_0(\mf p)$.
\subsection{A restatement in the archimedean aspect}\label{archrestat}

Given a smooth compactly supported function $V:(\Rr^{\times})^{r_1}\times(\Cc^{\times})^{r_2}\xrightarrow[]{}\Cc$ there exists a smooth vector $\phi_{\infty}\in\pi_{\infty}$ s.t. 

$$W_{\phi,\infty}\begin{pmatrix}
x & 0 \\
0 & 1 
\end{pmatrix}=V(x)$$

for $x\in (\Rr^{\times})^{r_1}\times(\Cc^{\times})^{r_2}$. See \cite{GGPS69}.

By the theorem of Diximier-Malliavin, \cite{MD78} $\exists$ $h_1,\dots,h_n \in C_c^{\infty}(G(F_{\infty}))$ s.t. $$\phi_{\infty}=\sum_{i=1}^n R(h_i)\phi_i$$ for some smooth vectors $\phi_i\in\pi_{\infty}$. Combining the two statements:

\begin{multline*} 
\sum_{m\in F^{\times}}V(m) W_{\varphi,f}\begin{pmatrix}
m\un p & 0 \\
0 & 1 
\end{pmatrix} K(m_p)
\\=\sum_{i=1}^n\sum_{m\in F^{\times}}W_{R(h_i)\phi_i,\infty}\begin{pmatrix}
m & 0 \\
0 & 1 
\end{pmatrix}W_{\varphi,f}\begin{pmatrix}
m\un p & 0 \\
0 & 1 
\end{pmatrix} K(m_p).
\end{multline*}

In this way we can reduce the main theorem to the following statement:

\begin{thm*}
Let $\pi$ be a cuspidal $\GL_2$-automorphic representation over $F$ of level $\mf N$ integer. Let $\phi$ be an automorphic form in $\pi$ that is of level $\mf N$. Let $K$ be an isotypic trace function defined over the residue field $k(\mf p)$ with an associated Fourier-M\"obius group contained in the standard Borel subgroup. There exists an absolute constant $s>0$ s.t.
for every $h_{\infty}\in C_c(F_{\infty})$
$$\sum_{m\in\Ft}K(m_{\mf p})W_{\phi,f}\begin{pmatrix}
m\pi_{\mf p } & 0 \\
0 & 1 
\end{pmatrix}W_{R(h_{\infty})\phi_{\infty}}\begin{pmatrix}
m & 0 \\
0 & 1 
\end{pmatrix}\ll_{\phi, F, h_{\infty},\delta} \cond(K)^s\Nm(\mf p)^{\frac{1}{2}-\delta}$$

for the $\delta$ and implied constants as stated in the main theorems.
\end{thm*}

This restatement lends itself well to applying the relative trace formula as we will see in the next section.
\section{Proof of the main theorems }\label{proofmt}
The proof the main theorem is by estimating the amplified spectral average of the correlation sum:
\begin{multline*}
\sum_{\pi}\Big ( |\sum_{\mf l\in \Lambda}x_{\mf l}\lambda_{\pi}(\mf l)|^2\Big )\sum_{\varphi_{\infty}\in \mathscr{B}(\pi_{\infty})} \sum_{\varphi_f\in \mathscr{B}(\pi_f, \mf N\mf p)}  \left|\sum_{m\in F^{\times}}K(m_{\mf p})W_{\phi,f}\begin{pmatrix}
m\pi_{\mf p } & 0 \\
0 & 1 
\end{pmatrix}W_{R(h_{\infty})\phi_{\infty}}\begin{pmatrix}
m & 0 \\
0 & 1 
\end{pmatrix}\right|^2 .\end{multline*}
Here $\Lambda$ is the set of integral ideals of size $L$ and $x_{\mf l}\in \Cc$  (to be chosen later).Also $\mathscr{B}(\pi, \mf a)$ for an ideal $\mf a\lhd \mcO_F$ is an orthonormal basis of the subspace of right $K_0(\mf a)$-invariant vectors of an irreducible cuspidal representation $\pi$. (Recall that the subspace of $K_0(\mf a)$-invariant vectors of an irreducible cuspidal representation $\pi$ is finite dimensional.)

Using the computation of section \S\ref{spectsidell}, we have
\begin{multline*}
=\sum_{\pi}\Big ( |\sum_{\mf l\in \Lambda}x_{\mf l}\lambda_{\pi}(\mf l)|^2\Big )\sum_{\varphi_{\infty}\in \mathscr{B}(\pi_{\infty})} \sum_{\varphi_f\in \mathscr{B}(\pi_f, \mf N\mf p)}  \left|\sum_{m\in F^{\times}}W_{R(h_{\infty})\varphi_\infty}\begin{pmatrix}
m & 0 \\
0 & 1 
\end{pmatrix}W_{R(h[1]_f)\varphi_f}\begin{pmatrix}
m & 0 \\
0 & 1 
\end{pmatrix}\right|^2 .\end{multline*}

 We will split the $\mf l$ sum into diagonal and non diagonal pieces as follows: 
\begin{multline*}
=\sum_{\pi}  \sum_{\mf l\in \Lambda}|x_{\mf l}\lambda_{\pi}(\mf l)|^2 \sum_{\varphi_{\infty}\in \mathscr{B}(\pi_{\infty})} \sum_{\varphi_f\in \mathscr{B}(\pi_f, \mf N\mf p)}  \left|\sum_{m\in F^{\times}}W_{R(h_{\infty})\varphi_\infty}\begin{pmatrix}
m & 0 \\
0 & 1 
\end{pmatrix}W_{R(h[1]_f)\varphi_f}\begin{pmatrix}
m & 0 \\
0 & 1 
\end{pmatrix}\right|^2 \end{multline*}
\begin{multline*}+ \sum_{\pi}\Big ( \sum_{\substack{\mf l_1, \mf l_2\in \Lambda\\ \mf l_1\neq \mf l_2}}x_{\mf l_1}\ov{x_{\mf l_2}} \lambda_{\pi}(\mf l_1)\ov{\lambda_{\pi}(\mf l_2)}\Big )  \sum_{\varphi_{\infty}\in \mathscr{B}(\pi_{\infty})} \sum_{\varphi_f\in \mathscr{B}(\pi_f, \mf N\mf p)} \left|\sum_{m\in F^{\times}}W_{R(h_{\infty})\varphi_\infty}\begin{pmatrix}
m & 0 \\
0 & 1 
\end{pmatrix}W_{R(h[1]_f)\varphi_f}\begin{pmatrix}
m & 0 \\
0 & 1 
\end{pmatrix}\right|^2 .\end{multline*}
\subsection{Proof of main theorem 1}
\subsubsection{The term $\mf l_1=\mf l_2$}
Let us look at the term $\mf l_1=\mf l_2$ first. We have by Rankin-Selberg theory for $\pi$ of conductor dividing $\mf N\mf p$

$$\sum_{\mf l\in \Lambda}|\lambda_{\pi}(\mf l)|^2\ll_{\epsilon, \mf N}(\cond(\pi))^{\epsilon}L^{1+\eps}\ll_{\epsilon, \mf N} (\cond(\pi_{\infty}))^{\epsilon}(\Nm \mf p)^{\epsilon}L^{1+\eps}.$$

Fix $\eps>0$, Substituting the above bound we need to look at

\begin{multline*}
\sum_{\pi}(\cond(\pi_{\infty}))^{\epsilon}\sum_{\varphi_{\infty}\in \mathscr{B}(\pi_{\infty})} \sum_{\varphi_f\in \mathscr{B}(\pi_f, \mf N\mf p)} \left|\sum_{m\in F^{\times}}W_{R(h_{\infty})\varphi,\infty}\begin{pmatrix}
m & 0 \\
0 & 1 
\end{pmatrix}W_{R(h[1]_f)\varphi_f}\begin{pmatrix}
m & 0 \\
0 & 1 
\end{pmatrix}\right|^2 .
\end{multline*}

We split the sum into $\cond(\pi_{\infty})<\Nm(\mf p)$ and  $\cond(\pi_{\infty})\geq \Nm(\mf p)$. In the first case by positivity we get

\begin{multline*}
\sum_{\substack{\pi\\ \cond(\pi_{\infty})<\Nm(\mf p)}}(\cond(\pi_{\infty}))^{\epsilon}\sum_{\varphi_{\infty}\in \mathscr{B}(\pi_{\infty})} \sum_{\varphi_f\in \mathscr{B}(\pi_f, \mf N\mf p)} \left|\sum_{m\in F^{\times}}W_{R(h_{\infty})\varphi,\infty}\begin{pmatrix}
m & 0 \\
0 & 1 
\end{pmatrix}W_{R(h[1]_f)\varphi_f}\begin{pmatrix}
m & 0 \\
0 & 1 
\end{pmatrix}\right|^2 
\\ \leq (\Nm(\mf p))^{\epsilon} \sum_{\substack{\pi}} \sum_{\varphi_{\infty}\in \mathscr{B}(\pi_{\infty})} \sum_{\varphi_f\in \mathscr{B}(\pi_f, \mf N\mf p)} \left|\sum_{m\in F^{\times}}W_{R(h_{\infty})\varphi,\infty}\begin{pmatrix}
m & 0 \\
0 & 1 
\end{pmatrix}W_{R(h[1]_f)\varphi_f}\begin{pmatrix}
m & 0 \\
0 & 1 
\end{pmatrix}\right|^2 .
\end{multline*}

We will prove by computing the geometric side of the relative trace formula in propositions \ref{diagbound} and \ref{nondiagbound}, the following bound. This bound uses the fact that, by positivity the cuspidal contribution is bounded by the overall spectral average which is equal to the geometric side.
\begin{multline*}
\sum_{\pi}\sum_{\varphi_{\infty}\in \mathscr{B}(\pi_{\infty})} \sum_{\varphi_f\in \mathscr{B}(\pi_f)} \left|\sum_{m\in F^{\times}}W_{R(h_{\infty})\varphi,\infty}\begin{pmatrix}
m & 0 \\
0 & 1 
\end{pmatrix}W_{R(h[1]_f)\varphi_f}\begin{pmatrix}
m & 0 \\
0 & 1 
\end{pmatrix}\right|^2 \ll_{f_{\infty},F} \cond(K)^s \Nm(\mf p)
\end{multline*}
for some $s>0$.
Using the result of section \S\ref{nonarchsupp}, we have

\begin{multline*}
\sum_{\pi}\sum_{\varphi_{\infty}\in \mathscr{B}(\pi_{\infty})} \sum_{\varphi_f\in \mathscr{B}(\pi_f,\mf N\mf p)} \left|\sum_{m\in F^{\times}}W_{R(h_{\infty})\varphi,\infty}\begin{pmatrix}
m & 0 \\
0 & 1 
\end{pmatrix}W_{R(h[1]_f)\varphi_f}\begin{pmatrix}
m & 0 \\
0 & 1 
\end{pmatrix}\right|^2 \ll_{f_{\infty},F} \cond(K)^s \Nm(\mf p)
\end{multline*}
for some $s>0$.

In the second case, we use rapid decay of the archimedean part in analogy with section 3.3.2 (equation 3.6) in \cite{Nel17} and the trivial bound for the non archimedean average which is polynomial in $\Nm(\mf p)$ (see section \S \ref{triv}) to get 
\begin{multline*}
\sum_{\substack{\pi\\ \cond(\pi_{\infty})\geq \Nm(\mf p)}}(\cond(\pi_{\infty}))^{\epsilon}\sum_{\varphi_{\infty}\in \mathscr{B}(\pi_{\infty})} \sum_{\varphi_f\in \mathscr{B}(\pi_f, \mf N\mf p)} \left|\sum_{m\in F^{\times}}W_{R(h_{\infty})\varphi,\infty}\begin{pmatrix}
m & 0 \\
0 & 1 
\end{pmatrix}W_{R(h[1]_f)\varphi_f}\begin{pmatrix}
m & 0 \\
0 & 1 
\end{pmatrix}\right|^2  \\ \ll_A (\Nm(\mf p))^{-A}
\end{multline*}
for any $A>0$.

Putting together the results, we conclude
\begin{multline*}
\sum_{\pi} \sum_{\mf l\in \Lambda}|x_{\mf l}\lambda_{\pi}(\mf l)|^2\sum_{\varphi_{\infty}\in \mathscr{B}(\pi_{\infty})} \sum_{\varphi_f\in \mathscr{B}(\pi_f, \mf N\mf p)}  \left|\sum_{m\in F^{\times}}W_{R(h_{\infty})\varphi_\infty}\begin{pmatrix}
m & 0 \\
0 & 1 
\end{pmatrix}W_{R(h[1]_f)\varphi_f}\begin{pmatrix}
m & 0 \\
0 & 1 
\end{pmatrix}\right|^2  \\ \ll_{f_{\infty},F,\epsilon}\cond(K)^s \Nm(\mf p)^{1+\eps} L^{1+\eps}\end{multline*}

\subsubsection{The term $\mf l_1\neq\mf l_2$}
We see that the non-diagonal part of the sum can be rewritten as
\begin{multline*}
\sum_{\pi}\Big ( \sum_{\substack{\mf l_1, \mf l_2\in \Lambda\\ \mf l_1\neq \mf l_2}}x_{\mf l_1}\ov{x_{\mf l_2}}\lambda_{\pi}(\mf l_1)\ov{\lambda_{\pi}(\mf l_2)}\Big ) \sum_{\varphi_{\infty}\in \mathscr{B}(\pi_{\infty})} \sum_{\varphi_f\in \mathscr{B}(\pi_f, \mf N\mf p)} \left|\sum_{m\in F^{\times}}W_{R(h_{\infty})\varphi_\infty}\begin{pmatrix}
m & 0 \\
0 & 1 
\end{pmatrix}W_{R(h[1]_f)\varphi_f}\begin{pmatrix}
m & 0 \\
0 & 1 
\end{pmatrix}\right|^2\\
    =\sum_{\pi}\sum_{\varphi_{\infty}\in \mathscr{B}(\pi_{\infty})} \sum_{\varphi_f\in \mathscr{B}(\pi_f, \mf N\mf p)} \left|\sum_{m\in F^{\times}}W_{R(f)\varphi}\begin{pmatrix}
m & 0 \\
0 & 1 
\end{pmatrix}\right|^2 
\end{multline*}

where $$f=\sum_{\substack{\mf l_1,\mf l_2\in \Lambda\\ \mf l_1 \neq \mf l_2}} x_{\mf l_1}\ov{x_{\mf l_2}}h[\mf l_1]*h[\mf l_2]^*.$$

Under the assumption that the Fourier-M\"obius group of the trace function is contained in the Borel subgroup, we will show by applying the relative trace formula to $f$ and bounding the geometric side the following estimate. This bound uses the fact that, by positivity the cuspidal contribution is bounded by the overall spectral average which is equal to the geometric side. 
\begin{multline*}
   \sum_{\pi}\sum_{\varphi_{\infty}\in \mathscr{B}(\pi_{\infty})} \sum_{\varphi_f\in \mathscr{B}(\pi_f)} \left|\sum_{m\in F^{\times}}W_{R(f)\varphi}\begin{pmatrix}
m & 0 \\
0 & 1 
\end{pmatrix}\right|^2\\ \ll_{f_{\infty},F} \cond(K)^s\sqrt{\Nm(\mf p)}.L^{3+\eps}+\cond(K)^s\frac{\Nm(\mf p)}{L}\sum\limits_{ \substack{\mf l_1\neq \mf l_2\\\mf l_1,\mf l_2\in \Lambda}} |x_{\mf l_1}\overline{x_{\mf l_2}}|.
\end{multline*}
The above bounds are obtained in propositions \ref{diagboundv} and \ref{nondiagboundv}.

Using the result of section \S \ref{nonarchsupp}, we have
\begin{multline*}
   \sum_{\pi}\sum_{\varphi_{\infty}\in \mathscr{B}(\pi_{\infty})} \sum_{\varphi_f\in \mathscr{B}(\pi_f, \mf N\mf p)} \left|\sum_{m\in F^{\times}}W_{R(f)\varphi}\begin{pmatrix}
m & 0 \\
0 & 1 
\end{pmatrix}\right|^2\\ \ll_{f_{\infty},F} \cond(K)^s\sqrt{\Nm(\mf p)}.L^{3+\eps}+\cond(K)^s\frac{\Nm(\mf p)}{L}\sum\limits_{ \substack{\mf l_1\neq \mf l_2\\\mf l_1,\mf l_2\in \Lambda}} |x_{\mf l_1}\overline{x_{\mf l_2}}|.
\end{multline*}

\subsubsection{Conclusion}

Combining all the bounds together for weights of size 1 i.e. $|x_{\mf l}|=1$, we obtain
\begin{multline*}
\sum_{\pi} |\sum_{\mf l\in \Lambda}x_{\mf l}\lambda_{\pi}(\mf l)|^2\sum_{\varphi_{\infty}\in \mathscr{B}(\pi_{\infty})} \sum_{\varphi_f\in \mathscr{B}(\pi_f, \mf N\mf p)} \left|\sum_{m\in F^{\times}}W_{R(h_{\infty})\varphi,\infty}\begin{pmatrix}
m & 0 \\
0 & 1 
\end{pmatrix}W_{R(h[1]_f)\varphi_f}\begin{pmatrix}
m & 0 \\
0 & 1 
\end{pmatrix}\right|^2 
\\ \ll_{f_{\infty},F} \cond(K)^s(\Nm(\mf p))^{1+\eps}.L^{1+\eps}+\cond(K)^s\sqrt{\Nm(\mf p)}.L^{3+\eps}.\end{multline*}

Let $\phi$ be a cuspidal automorphic form in the representation $\pi$ of conductor $K_0(\mf N)$ that is a pure tensor. By positivity, we have the following bound 
\begin{multline*}
|\sum_{\mf l\in \Lambda}x_{\mf l}\lambda_{\pi}(\mf l)|^2\left|\sum_{m\in F^{\times}}W_{R(h_{\infty})\phi,\infty}\begin{pmatrix}
m & 0 \\
0 & 1 
\end{pmatrix}W_{R(h[1]_f)\phi_f}\begin{pmatrix}
m & 0 \\
0 & 1 
\end{pmatrix}\right|^2 
\\ \ll_{f_{\infty},F} \cond(K)^s(\Nm(\mf p))^{1+\eps}.L^{1+\eps}+\cond(K)^s\sqrt{\Nm(\mf p)}.L^{3+\eps}.\end{multline*}

Choosing
\[   
x_{\mf l}=  
     \begin{cases}
       \text{sign}(\lambda_{\pi}(\mf l)) &\text{if $\mf l\in \Lambda$ and $\lambda_{\pi}(\mf l)\neq 0$} \\
       0  &\quad\text{otherwise} 
     \end{cases}
\]

By the prime number theorem for Rankin-Selberg L functions (see \cite{AVe10} for this choice of amplifier) we have with the above choice,

$$|\sum_{\mf l\in \Lambda}x_{\mf l}\lambda_{\pi}(\mf l)|=\sum_{\mf l\in \Lambda}|\lambda_{\pi}(\mf l)|\gg_{\pi} \frac{L}{(\log L)^2}.$$

Putting all this information together and choosing $L=(\Nm(\mf p))^{\frac{1}{4}}$ we get

$$\left|\sum_{m\in F^{\times}}W_{R(h_{\infty})\phi,\infty}\begin{pmatrix}
m & 0 \\
0 & 1 
\end{pmatrix} W_{\phi}\begin{pmatrix}
m\un p & 0 \\
0 & 1 
\end{pmatrix} K(m_p)\right|\\ \ll_{f_{\infty},F,\pi} \cond(K)^s(\Nm(\mf p))^{\frac{3}{8}+\eps}.$$

 This completes the proof modulo computing the geometric side and proving the bounds that we stated in the course of the proof. This is done in the following chapter.

\subsection{Proof of main theorem 2}

\subsubsection{The term $\mf l_1=\mf l_2$}
Let us look at the term $\mf l_1=\mf l_2$ first. 
$$\sum_{\mf l\in \Lambda}|x_{\mf l}\lambda_{\pi}(\mf l)|^2=\sum_{\substack{\mf l \text{ prime }\\ \Nm(\mf l)\sim L}}|\ov{\lambda_{\pi_0}(\mf l)}\lambda_{\pi}(\mf l)|^2+\sum_{\substack{\mf l \text{ prime }\\ \Nm(\mf l)\sim L}}|(-1)\lambda_{\pi}(\mf l^2)|^2$$

For the first sum, we use Cauchy-Schwarz inequality and bounds for fourth power sums of Hecke eigenvalues. To obtain a bound for the fourth power sum of Hecke eigenvalues we proceed in the same way as G.L\"u (\cite{Lu09}), the o(1) power dependence on the conductor can be argued in the same way as Remark 1.11 in Philippe Michel's lecture in the Park City lecture series \cite{PhM07}.For the second sum the properties of the Rankin-Selberg L-function suffice to conclude:

$$\sum_{\mf l\in \Lambda}|x_{\mf l}\lambda_{\pi}(\mf l)|^2\ll_{\eps,\mf N} (\cond(\pi_{\infty}))^{\eps}(\cond(\pi_{0,\infty}))^{\eps}L^{1+\eps}$$
Fix $\eps>0$, Substituting the above bound we need to look at
\begin{multline*}
\sum_{\pi}(\cond(\pi_{\infty}))^{\epsilon}\sum_{\varphi_{\infty}\in \mathscr{B}(\pi_{\infty})} \sum_{\varphi_f\in \mathscr{B}(\pi_f, \mf N\mf p)} \left|\sum_{m\in F^{\times}}W_{R(h_{\infty})\varphi,\infty}\begin{pmatrix}
m & 0 \\
0 & 1 
\end{pmatrix}W_{R(h[1]_f)\varphi_f}\begin{pmatrix}
m & 0 \\
0 & 1 
\end{pmatrix}\right|^2 .
\end{multline*}

We split the sum into $\cond(\pi_{\infty})<\Nm(\mf p)$ and  $\cond(\pi_{\infty})\geq \Nm(\mf p)$. In the first case, by positivity we get
\begin{multline*}
\sum_{\substack{\pi\\ \cond(\pi_{\infty})<\Nm(\mf p)}}(\cond(\pi_{\infty}))^{\epsilon}\sum_{\varphi_{\infty}\in \mathscr{B}(\pi_{\infty})} \sum_{\varphi_f\in \mathscr{B}(\pi_f, \mf N\mf p)} \left|\sum_{m\in F^{\times}}W_{R(h_{\infty})\varphi,\infty}\begin{pmatrix}
m & 0 \\
0 & 1 
\end{pmatrix}W_{R(h[1]_f)\varphi_f}\begin{pmatrix}
m & 0 \\
0 & 1 
\end{pmatrix}\right|^2 
\\ \leq (\Nm(\mf p))^{\epsilon} \sum_{\substack{\pi}} \sum_{\varphi_{\infty}\in \mathscr{B}(\pi_{\infty})} \sum_{\varphi_f\in \mathscr{B}(\pi_f, \mf N\mf p)} \left|\sum_{m\in F^{\times}}W_{R(h_{\infty})\varphi,\infty}\begin{pmatrix}
m & 0 \\
0 & 1 
\end{pmatrix}W_{R(h[1]_f)\varphi_f}\begin{pmatrix}
m & 0 \\
0 & 1 
\end{pmatrix}\right|^2 .
\end{multline*}

We will prove by computing the geometric side of the relative trace formula in propositions \ref{diagbound} and \ref{nondiagbound}, the following bound. This bound uses the fact that, by positivity the cuspidal contribution is bounded by the overall spectral average which is equal to the geometric side.
\begin{multline*}
\sum_{\pi}\sum_{\varphi_{\infty}\in \mathscr{B}(\pi_{\infty})} \sum_{\varphi_f\in \mathscr{B}(\pi_f)} \left|\sum_{m\in F^{\times}}W_{R(h_{\infty})\varphi,\infty}\begin{pmatrix}
m & 0 \\
0 & 1 
\end{pmatrix}W_{R(h[1]_f)\varphi_f}\begin{pmatrix}
m & 0 \\
0 & 1 
\end{pmatrix}\right|^2 \ll_{f_{\infty},F} \cond(K)^s \Nm(\mf p)
\end{multline*}
for some $s>0$.
Using the result of section \S\ref{nonarchsupp}, we have
\begin{multline*}
\sum_{\pi}\sum_{\varphi_{\infty}\in \mathscr{B}(\pi_{\infty})} \sum_{\varphi_f\in \mathscr{B}(\pi_f,\mf N\mf p)} \left|\sum_{m\in F^{\times}}W_{R(h_{\infty})\varphi,\infty}\begin{pmatrix}
m & 0 \\
0 & 1 
\end{pmatrix}W_{R(h[1]_f)\varphi_f}\begin{pmatrix}
m & 0 \\
0 & 1 
\end{pmatrix}\right|^2 \ll_{f_{\infty},F} \cond(K)^s \Nm(\mf p)
\end{multline*}
for some $s>0$.

In the second case, we use rapid decay of the archimedean part in analogy with section 3.3.2 (equation 3.6) in \cite{Nel17} and the trivial bound for the non archimedean average which is polynomial in $\Nm(\mf p)$ (see section \S \ref{triv}) to get 
\begin{multline*}
\sum_{\substack{\pi\\ \cond(\pi_{\infty})\geq \Nm(\mf p)}}(\cond(\pi_{\infty}))^{\epsilon}\sum_{\varphi_{\infty}\in \mathscr{B}(\pi_{\infty})} \sum_{\varphi_f\in \mathscr{B}(\pi_f, \mf N\mf p)} \left|\sum_{m\in F^{\times}}W_{R(h_{\infty})\varphi,\infty}\begin{pmatrix}
m & 0 \\
0 & 1 
\end{pmatrix}W_{R(h[1]_f)\varphi_f}\begin{pmatrix}
m & 0 \\
0 & 1 
\end{pmatrix}\right|^2  \\ \ll_A (\Nm(\mf p))^{-A}
\end{multline*}
for any $A>0$.

Putting together the results, we conclude
\begin{multline*}
\sum_{\pi} \sum_{\mf l\in \Lambda}|x_{\mf l}\lambda_{\pi}(\mf l)|^2\sum_{\varphi_{\infty}\in \mathscr{B}(\pi_{\infty})} \sum_{\varphi_f\in \mathscr{B}(\pi_f, \mf N\mf p)}  \left|\sum_{m\in F^{\times}}W_{R(h_{\infty})\varphi_\infty}\begin{pmatrix}
m & 0 \\
0 & 1 
\end{pmatrix}W_{R(h[1]_f)\varphi_f}\begin{pmatrix}
m & 0 \\
0 & 1 
\end{pmatrix}\right|^2  \\ \ll_{f_{\infty},F,\epsilon}\cond(K)^s \Nm(\mf p)^{1+\eps} L^{1+\eps}\end{multline*}

\subsubsection{The term $\mf l_1\neq\mf l_2$}
We see that the non-diagonal part of the sum can be rewritten as
\begin{multline*}
\sum_{\pi}\Big ( \sum_{\substack{\mf l_1, \mf l_2\in \Lambda\\ \mf l_1\neq \mf l_2}}x_{\mf l_1}\ov{x_{\mf l_2}}\lambda_{\pi}(\mf l_1)\ov{\lambda_{\pi}(\mf l_2)}\Big ) \sum_{\varphi_{\infty}\in \mathscr{B}(\pi_{\infty})} \sum_{\varphi_f\in \mathscr{B}(\pi_f, \mf N\mf p)} \left|\sum_{m\in F^{\times}}W_{R(h_{\infty})\varphi_\infty}\begin{pmatrix}
m & 0 \\
0 & 1 
\end{pmatrix}W_{R(h[1]_f)\varphi_f}\begin{pmatrix}
m & 0 \\
0 & 1 
\end{pmatrix}\right|^2\\
    =\sum_{\pi}\sum_{\varphi_{\infty}\in \mathscr{B}(\pi_{\infty})} \sum_{\varphi_f\in \mathscr{B}(\pi_f, \mf N\mf p)} \left|\sum_{m\in F^{\times}}W_{R(f)\varphi}\begin{pmatrix}
m & 0 \\
0 & 1 
\end{pmatrix}\right|^2 
\end{multline*}

where $$f=\sum_{\substack{\mf l_1,\mf l_2\in \Lambda\\ \mf l_1 \neq \mf l_2}} x_{\mf l_1}\ov{x_{\mf l_2}}h[\mf l_1]*h[\mf l_2]^*.$$

Under the assumption that the Fourier-M\"obius group of the trace function is contained in the Borel subgroup, we will show by applying the relative trace formula to $f$ and bounding the geometric side the following estimate. This bound uses the fact that, by positivity the cuspidal contribution is bounded by the overall spectral average which is equal to the geometric side. 
\begin{multline*}
   \sum_{\pi}\sum_{\varphi_{\infty}\in \mathscr{B}(\pi_{\infty})} \sum_{\varphi_f\in \mathscr{B}(\pi_f)} \left|\sum_{m\in F^{\times}}W_{R(f)\varphi}\begin{pmatrix}
m & 0 \\
0 & 1 
\end{pmatrix}\right|^2\\ \ll_{f_{\infty},F} \cond(K)^s\sqrt{\Nm(\mf p)}.L^{2+\eps}\sum\limits_{ \substack{\mf l_1\neq \mf l_2\\\mf l_1,\mf l_2\in \Lambda}} |x_{\mf l_1}\overline{x_{\mf l_2}}|+\cond(K)^s\frac{\Nm(\mf p)}{L}\sum\limits_{ \substack{\mf l_1\neq \mf l_2\\\mf l_1,\mf l_2\in \Lambda}} |x_{\mf l_1}\overline{x_{\mf l_2}}|.
\end{multline*}
The above bounds are obtained in propositions \ref{diagbounddfi} and \ref{nondiagbounddfi}. (The above estimate works for any $s\geq 4$. )

Using the result of section \S \ref{nonarchsupp}, we have
\begin{multline*}
   \sum_{\pi}\sum_{\varphi_{\infty}\in \mathscr{B}(\pi_{\infty})} \sum_{\varphi_f\in \mathscr{B}(\pi_f, \mf N\mf p)} \left|\sum_{m\in F^{\times}}W_{R(f)\varphi}\begin{pmatrix}
m & 0 \\
0 & 1 
\end{pmatrix}\right|^2\\ \ll_{f_{\infty},F} L^{2+\eps}\sum\limits_{ \substack{\mf l_1\neq \mf l_2\\\mf l_1,\mf l_2\in \Lambda}} |x_{\mf l_1}\overline{x_{\mf l_2}}|+\cond(K)^s\frac{\Nm(\mf p)}{L}\sum\limits_{ \substack{\mf l_1\neq \mf l_2\\\mf l_1,\mf l_2\in \Lambda}} |x_{\mf l_1}\overline{x_{\mf l_2}}|.
\end{multline*}

\subsubsection{Conclusion}

Combining all the bounds together and using Cauchy-Schwartz and Rankin-Selberg bounds for the bilinear sums in the amplifier weights, we obtain
\begin{multline*}
\sum_{\pi} |\sum_{\mf l\in \Lambda}x_{\mf l}\lambda_{\pi}(\mf l)|^2\sum_{\varphi_{\infty}\in \mathscr{B}(\pi_{\infty})} \sum_{\varphi_f\in \mathscr{B}(\pi_f, \mf N\mf p)} \left|\sum_{m\in F^{\times}}W_{R(h_{\infty})\varphi,\infty}\begin{pmatrix}
m & 0 \\
0 & 1 
\end{pmatrix}W_{R(h[1]_f)\varphi_f}\begin{pmatrix}
m & 0 \\
0 & 1 
\end{pmatrix}\right|^2 
\\ \ll_{f_{\infty},F} \cond(K)^s(\Nm(\mf p))^{1+\eps}.L^{1+\eps}+\cond(K)^s\sqrt{\Nm(\mf p)}.L^{4+\eps}.\end{multline*}

Let $\phi$ be a cuspidal automorphic form in the representation $\pi_0$ of conductor a divisor of $\mf N \mf p$ that is a pure tensor. By positivity, we have the following bound 
\begin{multline*}
|\sum_{\mf l\in \Lambda}x_{\mf l}\lambda_{\pi_0}(\mf l)|^2\left|\sum_{m\in F^{\times}}W_{R(h_{\infty})\phi,\infty}\begin{pmatrix}
m & 0 \\
0 & 1 
\end{pmatrix}W_{R(h[1]_f)\phi_f}\begin{pmatrix}
m & 0 \\
0 & 1 
\end{pmatrix}\right|^2 
\\ \ll_{f_{\infty},F} \cond(K)^s(\Nm(\mf p))^{1+\eps}.L^{1+\eps}+\cond(K)^s\sqrt{\Nm(\mf p)}.L^{4+\eps}.\end{multline*}

Choosing
\[   
x_{\mf l}=  
     \begin{cases}
       \ov{\lambda_{\pi_0}(\mf l)} &\text{if $\mf l \sim L$ prime and $(\mf l,\mf N\mf p)=\mcO_{F}$} \\
       -1 &\text{if $\mf l=\mf l'^2$, $\mf l' \sim L$ prime and $(\mf l',\mf N\mf p)=\mcO_{F}$} \\
       0  &\quad\text{otherwise} 
     \end{cases}
\]

By the prime number theorem for the Dedekind zeta function (see \cite{DFI94} for this choice of amplifier) we have with the above choice,

$$|\sum_{\mf l\in \Lambda}x_{\mf l}\lambda_{\pi_0}(\mf l)|=\sum_{\substack{\mf l\sim L \text{  prime}\\ (\mf l,\mf N\mf p)=(1)}}(|\lambda_{\pi_0}(\mf l)|^2-\lambda_{\pi_0}(\mf l^2))=\sum_{\substack{\mf l\sim L \text{  prime}\\ (\mf l,\mf N\mf p)=(1)}} 1 \gg_{F} L^{1+o(1)}.$$

Putting all this information together and choosing $L=(\Nm(\mf p))^{\frac{1}{6}}$ we get

$$\left|\sum_{m\in F^{\times}}W_{R(h_{\infty})\phi,\infty}\begin{pmatrix}
m & 0 \\
0 & 1 
\end{pmatrix} W_{\phi}\begin{pmatrix}
m\un p & 0 \\
0 & 1 
\end{pmatrix} K(m_p)\right|\\ \ll_{f_{\infty},F,\pi} \cond(K)^s(\Nm(\mf p))^{\frac{5}{12}+\eps}.$$

 Note that the implied constant in this case is polynomial in the archimedean part of the conductor of $\pi$. This completes the proof modulo computing the geometric side and proving the bounds that we stated in the course of the proof.  This is done in the following chapter.

\section{The geometric side}\label{geomside}

The aim of this chapter is to complete the proof of the main theorem by computing the geometric side and bounding it. We obtain estimates in the case of $f=h[1]*h[1]^*$ and in the case $f=\sum_{\mf l_1\neq\mf l_2\in \Lambda} x_{\mf l_1}\ov{x_{\mf l_2}}h[\mf l_1]*h[\mf l_2]^*$. The bounds for the case of $f=h[1]*h[1]^*$ are in propositions \ref{diagbound} and \ref{nondiagbound}, for the Venkatesh amplifier in propositions \ref{diagboundv} and \ref{nondiagboundv} and for the DFI amplifier in propositions  \ref{diagbounddfi} and \ref{nondiagbounddfi}.
  We will now compute the integral using the geometric expansion of the kernel
     \begin{align*}      
     I_{m,n} &=\int\limits_{[N]\times [N]}K(x,y)\psi_m(x)\ov{\psi_n(y)}dxdy \\
             &=\int\limits_{[N]\times [N]}\sum_{\gamma\in \overline{G}(F)}f(x^{-1}\gamma y)\psi_m(x)\ov{\psi_n(y)}dxdy.
     \end{align*}
     
  Unfolding the integral and rearranging according to the orbits we get(see \cite{KL06} for details)
    
     $$I_{m,n}=\sum_{[\delta]}I_{\delta}$$
     where 
       $$I_{\delta}=\int\limits_{H_{\delta}(F)/H(\Aa)}f(x^{-1}\delta y)\psi_m(x)\ov{\psi_n(y)}dxdy$$
     where $$[\delta]=\{x^{-1}\delta y|x,y\in N(F)\}\text{ for }\delta\in \overline{G}(F)$$ 
     
   and    $$H=N\times N\text{ with    }H_{\delta}=\{(x,y)\in H|x^{-1}\delta y=\delta\}\text{ for }\delta\in\overline{G}.$$ The $I_{\delta}$ depend on $m,n$ but the dependence is dropped for simplicity.

     Further $I_{\delta}=0$ unless $\psi_m(x)\ov{\psi_n(y)}=1$ for $(x,y)\in H_{\delta}(\Aa)$. If $\psi_m(x)\ov{\psi_n(y)}=1$ for $(x,y)\in H_{\delta}(\Aa)$, we will say $[\delta]$ is relevant. We will now calculate $I_{\delta}$ for relevant orbits.
\begin{obs}

 Let $\mf q\neq \mf p$ be a finite place, we know that $f$ is $K_0(\mf N)_\mf q$-bi invariant. In particular
    $$f(g)=f(g\delta_{\mf q}\begin{pmatrix}
1 & t \\
0 & 1 
\end{pmatrix}=f(\delta_{\mf q}\begin{pmatrix}
1 & t \\
0 & 1 
\end{pmatrix}g)$$ for $t\in \mcO_\mf q$

Using this in $I_{\delta}$ we get  $I_{\delta}\neq 0$ for some $\delta$ only if  $m\mcO_\mf q\subseteq \cond(\psi)_{\mf q}$ and $n\mcO_\mf q\subseteq \cond(\psi)_{\mf q}$ for all finite place $\mf q\neq \mf p$. 
\end{obs}

As discussed in the section before, we will apply the relative trace formula to the case $f=h[1]*h[1]^*$ and to the case $f=\sum_{\mf l_1\neq\mf l_2\in \Lambda} x_{\mf l_1}\ov{x_{\mf l_2}}h[\mf l_1]*h[\mf l_2]^*$. In both these cases the above observation is applicable.

\subsection{The diagonal contribution}

Recall that this corresponds to the case:
$$\delta= \begin{pmatrix}
\lambda & 0 \\
0 & 1 
\end{pmatrix}. $$

Recall that $\delta$ is relevant only if $$n=m\lambda.$$

    \begin{align*}  
    I_{\delta}&= \int\limits_{\{t(\lambda ,1)|t\in F\}\backslash \Aa \times\Aa}f\begin{pmatrix}
1 & -t_1 \\
0 & 1 
\end{pmatrix} \begin{pmatrix}
\lambda & 0 \\
0 & 1 
\end{pmatrix} \begin{pmatrix}
1 & t_2 \\
0 & 1 
\end{pmatrix}\psi(mt_1-nt_2)dt_1dt_2
\\ &= \int\limits_{\{t(\lambda ,1)|t\in F\}\backslash \Aa \times\Aa}f\begin{pmatrix}
\lambda & \lambda t_2-t_1 \\
0 & 1 
\end{pmatrix}\psi(mt_1-nt_2)dt_1dt_2.
\end{align*}

Set $x=\lambda t_2-t_1$ and $y=t_2$ to get

$$I_{\delta}= \int\limits_{\{0\}\times F\backslash \Aa \times\Aa}f\begin{pmatrix}
\lambda & x \\
0 & 1 
\end{pmatrix}\psi(-mx)dxdy=c\int\limits_{ \Aa }f\begin{pmatrix}
\lambda & x \\
0 & 1 
\end{pmatrix}\psi(-mx)dx$$

with $c=meas(F\backslash\Aa)$.

Recall that we will work with two cases $$f=\sum_{\substack{\mf l_1,\mf l_2\in \Lambda \\ \mf l_1 \neq \mf l_2}} x_{\mf l_1}\ov{x_{\mf l_2}}h[\mf l_1]*h[\mf l_2]^*:=\sum_{\substack{\mf l_1,\mf l_2\in \Lambda \\ \mf l_1 \neq \mf l_2}} x_{\mf l_1}\ov{x_{\mf l_2}}f[\mf l_1,\mf l_2]$$ and $$f=h[1]*h[1]^*.$$

Since $f$ and $\psi$ can be factored, the integral factors into local integrals.

\subsubsection{Local computation at $\mf q|\mf l_1\mf l_2$, $\mf l_1\neq \mf l_2$} \label{diag2}
 Let $v_{\mf q}(\mf l_1)=d_1$ and $v_{\mf q}(\mf l_2)=d_2$ :
\begin{multline*}
    \int\limits_{ F_{\mf q}}f_{\mf q}[\mf l_1,\mf l_2]\begin{pmatrix}
\lambda & -x \\
0 & 1 
\end{pmatrix}\psi_{\mf q}(mx)dx\\ =\frac{1}{\sqrt{\Nm(\mf q^{d_1+d_2})}}\int\limits_{ F_{\mf q} } \int\limits_{\overline{G(F_{\mf q})}} \ind_{Z(F_{\mf q})M(\mf l_1)_{\mf q}} \begin{pmatrix}
\lambda & -x \\
0 & 1 
\end{pmatrix}g^{-1}) \ind_{Z(F_{\mf q}) M(\mf l_2)_{\mf q}} (g^{-1} ) \psi_{\mf q}(m x) dg dx
\end{multline*}
where $$M(\mf l)_{\mf q}= \{M\in M_2(\mcO_{\mf q})|\text{det}(M)\in \mf l\mcO_{\mf q}^*\}.$$
\begin{multline*}
    \int\limits_{ F_{\mf q}}f_{\mf q}[\mf l_1,\mf l_2]\begin{pmatrix}
\lambda & -x \\
0 & 1 
\end{pmatrix}\psi_{\mf q}(mx)dx =\frac{1}{\sqrt{\Nm(\mf q^{d_1+d_2})}}\sum_{a=0}^{d_2}\sum_{b\text{ mod }\un q^a}\int\limits_{ F_{\mf q} } \ind_{Z(F_{\mf q})M(\mf l_1)_{\mf q}} \left(\begin{pmatrix}
\lambda & -x \\
0 & 1 
\end{pmatrix}\begin{pmatrix}
\un q^a & b \\
0 & \un q^{d_2-a} 
\end{pmatrix}\right) \psi_{\mf q}(m x) \\\int\limits_{\overline{G(F_{\mf q})}}  \ind_{Z(F_{\mf q})\begin{pmatrix}
\un q^a & b \\
0 & \un q^{d_2-a} 
\end{pmatrix}K_{\mf q}} (g^{-1} )  dg dx
\end{multline*}
Simplifying the above we get
\begin{align*}
\int\limits_{ F_{\mf q}}f_{\mf q}[\mf l_1,\mf l_2]\begin{pmatrix}
\lambda & -x \\
0 & 1 
\end{pmatrix}\psi_{\mf q}(mx)dx &=\frac{1}{\sqrt{\Nm(\mf q^{d_1+d_2})}}\sum_{a=0}^{d_2}\sum_{b\text{ mod }\un q^a}\int\limits_{ F_{\mf q} } \ind_{Z(F_{\mf q})M(\mf l_1)_{\mf q}} \begin{pmatrix}
\lambda\un q^a & b\lambda-x\un q^{d_2-a} \\
0 & \un q^{d_2-a} 
\end{pmatrix} 
 \psi_{\mf q}(m x)\\&=\frac{1}{\sqrt{\Nm(\mf q^{d_1+d_2})}}\sum_{a=0}^{d_2}\sum_{b\text{ mod }\un q^a}\sum_{\alpha\in\Zz}\ind_{v_{\mf q}(\lambda)=d_1-d_2+2\alpha}\\& \int\limits_{ F_{\mf q} } \ind_{M(\mf l_1)_{\mf q}} \begin{pmatrix}
\lambda\un q^{a-\alpha} & b\lambda\un q^{-\alpha}-x\un q^{d_2-a-\alpha} \\
0 & \un q^{d_2-a-\alpha} 
\end{pmatrix} 
 \psi_{\mf q}(m x)dx \\ &=\frac{1}{\sqrt{\Nm(\mf q^{d_1+d_2})}}\sum_{a=0}^{d_2}\sum_{b\text{ mod }\un q^a}\psi_{\mf q}(nb\un q^{a-d_2})\sum_{\alpha=d_2-d_1-a}^{d_2-a}\ind_{v_{\mf q}(\lambda)=d_1-d_2+2\alpha}\\&
 \int\limits_{ \mf q^{a-d_2+\alpha}\mcO_{\mf q}} \psi_{\mf q}(mx) dx\\ &=\frac{1}{\sqrt{\Nm(\mf q^{d_1+d_2})}}\sum_{a=0}^{d_2}\sum_{\alpha=d_2-d_1-a}^{d_2-a}\ind_{v_{\mf q}(\lambda)=d_1-d_2+2\alpha}\ind_{v_{\mf q}(m)\geq M(\mf q, a,\alpha)}(\Nm(\mf q))^{(\frac{d_1+d_2}{2}-\frac{v_{\mf q}(\lambda)}{2})}.
\end{align*}

Here $$M(\mf q, a,\alpha)=\max(d_2-a-\alpha,d_2-a-\alpha+d_2-d_1-\alpha)+v_{\mf q}(\cond(\psi_{\mf q}))$$

\subsubsection{Local computation at $\mf p$} \label{diag3}
\begin{multline*}
\int\limits_{ F_{\mf p}}f_{\mf p}[l_1,l_2]\begin{pmatrix}
\lambda & -x \\
0 & 1 
\end{pmatrix}\psi_{\mf p}(mx)dx =\frac{1}{\meas(K_0(\mf p))^2 \Nm(\mf p)}\sum_{a,b\in k(\mf p)}\whK(a)\overline{\whK(b)}\\
\int\limits_{ F_q } \int\limits_{\overline{G(F_{\mf p})}} \ind_{Z(F_{\mf p})\begin{pmatrix}
p & a \\
0 & 1 
\end{pmatrix}K_0(\mf p)} \begin{pmatrix}
\lambda & -x \\
0 & 1 
\end{pmatrix}g^{-1} )  \ind_{Z(F_{\mf p})\begin{pmatrix}
\un p & b \\
0 & 1 
\end{pmatrix}K_0(\mf p)} (g^{-1} ) \psi_{\mf p}(m x) dg dx.
\end{multline*}
  The integrand is non zero if and only if 
  
  $$   \begin{pmatrix}
\lambda & -x \\
0 & 1 
\end{pmatrix}\begin{pmatrix}
\un p & b \\
0 & 1 
\end{pmatrix}K_0(\mf p)  \in Z(F_{\mf p})\begin{pmatrix}
\un p & a \\
0 & 1 
\end{pmatrix}K_0(\mf p).$$ 
Considering determinants, this implies:

$$\lambda\in (\Ft_{\mf p})^2\mcOt_{\mf p}.$$

Let $v_{\mf p}(\lambda)=2d$. The above condition is equivalent to 
 $$   \begin{pmatrix}
\lambda \un p^{-d+1} & -\un p^{-d}(-b\lambda +x) \\
0 & \un p^{-d} 
\end{pmatrix}  \in  \begin{pmatrix}
\un p & a \\
0 & 1 
\end{pmatrix}K_0(\mf p).$$ 

Looking at the bottom elements we conclude that in fact 
$$d=0.$$ 

The conclusion is that,
 \begin{multline*}
 \ind_{Z(F_{\mf p})\begin{pmatrix}
\un p & a \\
0 & 1 
\end{pmatrix}K_0(\mf p)} \begin{pmatrix}
\lambda & -x \\
0 & 1 
\end{pmatrix}g^{-1} )  \ind_{Z(F_{\mf p})\begin{pmatrix}
\un p & b \\
0 & 1 
\end{pmatrix}K_0(\mf p)} (g^{-1} )\\ =\ind_{\lambda\in \mcO_{\mf p}^*}\ind_{x\in -(a+b\lambda)+\un p\mcO_{\mf p}} \ind_{Z(F_{\mf p})K_0(\mf p)\begin{pmatrix}
\un p & b \\
0 & 1 
\end{pmatrix}^{-1}} (g ).
 \end{multline*}
It follows that
\begin{multline*}
\int\limits_{ F_{\mf p}}f_{\mf p}[\mf l_1,\mf l_2]\begin{pmatrix}
\lambda & -x \\
0 & 1 
\end{pmatrix}\psi_{\mf p}(mx)dx
\\ =\ind_{\lambda\in \mcO_{\mf p}^*}\frac{1}{\Nm(\mf p)^2\meas(K_0(\mf p))}\ind_{v_{\mf p}(m)\geq -1} \sum_{a,b\in k(\mf p)}\whK(a)\overline{\whK(b)} \psi_{\mf p}(-(ma+bn)). \end{multline*}
Since $\lambda\in\mcO_{\mf p}^*$ satisfies $n=m\lambda$, we have $v_{\mf p}(n)\geq -1$. Suppose $$m=\frac{m_{\mf p}}{\un p}+\mcO_{\mf p}$$ and $$n=\frac{n_{\mf p}}{\un p}+\mcO_{\mf p}$$ The above simplifies to 
\begin{multline*}
\int\limits_{ F_{\mf p}}f_{\mf p}[\mf l_1,\mf l_2]\begin{pmatrix}
\lambda & -x \\
0 & 1 
\end{pmatrix}\psi_{\mf p}(mx)dx
=\ind_{\lambda\in \mcO_{\mf p}^*}\frac{1}{\Nm(\mf p)\meas(K_0(\mf p))}\ind_{v_{\mf p}(m)\geq -1} \ind_{v_{\mf p}(n)\geq -1} K(m_{\mf p}) \overline{K(-n_{\mf p})}.      
\end{multline*}
\subsubsection{Local computation for $\mf q\nmid \mf l\mf p$} \label{diag4}

Similarly for $\mf q\nmid \mf l\mf p$ we have

 $$\int\limits_{ F_{\mf q}}f_{\mf q}[\mf l_1,\mf l_2]\begin{pmatrix}
\lambda & -x \\
0 & 1 
\end{pmatrix}\psi_{\mf q}(mx)dx=\meas(K_0(\mf N)_{\mf q})\ind_{\lambda\in\mcO_{\mf q}^*}\ind_{m\in\cond(\psi)_{\mf q}}.$$

\subsubsection{Archimedean local computations} \label{diagarch}
Let us now consider the archimedean part 
$$\phi(m,\lambda)= \int\limits_{ F_{\infty}}f_{\infty}[\mf l_1,\mf l_2]\begin{pmatrix}
\lambda & -x \\
0 & 1 
\end{pmatrix}\psi(mx)dx.$$ Note that this does not depend on $\mf l_1$ and $\mf l_2$.

Since $f_{\infty}[\mf l_1,\mf l_2]$ is compactly supported on $\PGL_2(F_{\infty})$, we conclude that $\phi(m,\lambda)=0$ except for $\lambda$ whose archimedean embedding lies in a compact subset of $F_{\infty}^{\times}$ (the compact set depends only on $f_{\infty}$ and not on $m$.). Further since $f_{\infty}[\mf l_1,\mf l_2]$ is smooth, $\phi(m,\lambda)$ is a Schwartz function in $m$.
\subsection{The diagonal contribution and bounds}
Combining the results of sections  \S\ref{diag2} , \S\ref{diag3} and \S\ref{diag4}  we get the following propositions
\subsubsection{Diagonal contribution for $f=h[1]*h[1]^*$}
\begin{prop}
The diagonal contribution for $f=h[1]*h[1]^*$ i.e. the sum over $I_\delta$, $\delta$ of the form $\begin{pmatrix}
\lambda & 0 \\
0 & 1 
\end{pmatrix}$, $\lambda\in \Ft$ reduces to a single term $\lambda=\frac{n}{m}$: 
\begin{align*}
A_1(m,n)&:=\sum_{\delta} I_{\delta} \\&=c'\ind_{v_{\mf p}(m)\geq -1} \ind_{v_{\mf p}(n)\geq -1}\ind_{m,n\in \cond(\psi)^{(\mf p)}}  K(m_p) \overline{K(-n_p)} \phi(m,\lambda)
\end{align*}

where $$c'=\frac{\meas(K_0(\mf N))\meas(\Aa/F)}{\Nm(\mf p)\meas(K_0(\mf p))} \text{ and } \lambda:=\frac{n}{m}. $$
\end{prop}

Now we average over $m,n\in F^{\times}$,

\begin{prop} \label{diagbound}
$$\left|\sum_{m\in F^{\times}}\sum_{n\in F^{\times}}A_1(m,n)\right|\ll_{f_{\infty},\mf N}\norm{K}_{\infty}^2\Nm(\mf p).$$

\end{prop}
\proof
\begin{align*}
    \left |\sum_{m\in F^{\times}}\sum_{n\in F^{\times}}A(m,n)\right |  &=\left |\sum\limits_{m, n \in \cond(\psi)\mf P^{-1}}c'  K(m_p) \overline{K(-n_p)} \phi(m,\lambda)\right |\\ 
     &\ll_{F} \norm{K}_{\infty}^2\sum\limits_{m, n \in \cond(\psi)\mf P^{-1}} \ind_{ \substack{\lambda\in\mcOt \\ n=m\lambda}} \left | \phi(m,\lambda)\right |\\ &= \norm{K}_{\infty}^2 \sum\limits_{m \in \cond(\psi)\mf P^{-1}}\sum_{\lambda\in \mcOt}\left | \phi(m,\lambda)\right |.
\end{align*}
The $\lambda$ sum is finite by \ref{diagarch} and we estimate the $m$ sum using lemma \ref{lpc} to get the result.
\subsubsection{Diagonal contribution for $f=\sum\limits_{\mf l_1 \neq \mf l_2} x_{\mf l_1}\ov{x_{\mf l_2}}h[\mf l_1]*h[\mf l_2]^*$}
\begin{prop}\label{diagl}
The diagonal contribution for $$f=\sum_{\substack{\mf l_1,\mf l_2\in \Lambda\\ \mf l_1 \neq \mf l_2}} x_{\mf l_1}\ov{x_{\mf l_2}}h[\mf l_1]*h[\mf l_2]^*$$ i.e. the sum over $I_\delta$, $\delta$ of the form $\begin{pmatrix}
\lambda & 0 \\
0 & 1 
\end{pmatrix}$, $\lambda\in \Ft$ reduces to a single term $\lambda=\frac{n}{m}$,
\begin{multline*}
    A_2(m,n):=\sum_{\delta} I_{\delta}=c'\ind_{v_{\mf p}(m)\geq -1} \ind_{v_{\mf p}(n)\geq -1}\ind_{m,n\in \cond(\psi)^{(\mf p)}}\\ \Bigg{(} \sum\limits_{ \substack{\mf l_1\neq \mf l_2\\\mf l_1,\mf l_2\in \Lambda}} x_{\mf l_1}\overline{x_{\mf l_2}} \ind_{\lambda\in(\whO^{\times})^{(\mf l_1\mf l_2)}} K(m_p) \overline{K(-n_p)}\\  \prod_{\mf q|\mf l_1\mf l_2}\Big{(}\sum_{a=0}^{v_{\mf q}(\mf l_2)}\sum_{\alpha=v_{\mf q}(\mf l_2)-v_{\mf q}(\mf l_1)-a}^{v_{\mf q}(\mf l_2)-a}\ind_{v_{\mf q}(\lambda)=v_{\mf q}(\mf l_1)-v_{\mf q}(\mf l_2)+2\alpha}\ind_{v_{\mf q}(m)\geq M(\mf q, a,\alpha)}(\Nm(\mf q))^{-\frac{v_{\mf q}(\lambda)}{2}}\Big{)}\Bigg{)}\phi(m,\lambda) 
\end{multline*}

where $$c'=\frac{\meas(K_0(\mf N))\meas(\Aa/F)}{\Nm(\mf p)\meas(K_0(\mf p))}\text{ and }\lambda=\frac{n}{m} $$
$$ M(\mf q, a,\alpha)=\max(v_{\mf q }(\mf l_2)-a-\alpha,v_{\mf q }(\mf l_2)-a-\alpha+v_{\mf q }(\mf l_2)-v_{\mf q }(\mf l_1)-\alpha)+v_{\mf q}(\cond(\psi_{\mf q}))$$
\end{prop}

As discussed in \ref{diagarch}, $\phi(m,\lambda)=0$ unless $\lambda$ lies in a compact subset of $F_{\infty}^{\times}$. For $L$ large enough we will have for the $\lambda$ occurring in $A_2(m,n)$  $$\prod_{\mf q|\mf l_1\mf l_2}\Nm(\mf q)^{v_{\mf q}(\lambda)}\sim 1.$$ Averaging over $m,n$ we have for $L$ large enough, we have the following bounds:
\begin{prop}[Bounds for the Venkatesh amplifier]\label{diagboundv}

$$\left|\sum_{m\in F^{\times}}\sum_{n\in F^{\times}}A_2(m,n)\right|\ll_{f_{\infty},\mf N}\norm{K}_{\infty}^2\frac{\Nm(\mf p)}{L}\sum\limits_{ \substack{\mf l_1\neq \mf l_2\\\mf l_1,\mf l_2\in \Lambda}} |x_{\mf l_1}\overline{x_{\mf l_2}}|.$$

\end{prop}

\proof

 We first note that $c'$ is of size 1 and $$\Lambda=\{\mf l| \mf l\text{ prime and }\mf l\sim L \text{ coprime to }\mf N\mf p\}.$$ 
 So $\mf l_1\neq \mf l_2$ means for each summand in the diagonal contribution we have only the two possibilities $\mf q=\mf l_1$ or $\mf q =\mf l_2$. Moreover the archimedean condition and the remark above implies the only possibilities for $\lambda$ is $v_{l_1}(\lambda)=1$ and $v_{l_1}(\lambda)=-1$ or vice versa. Putting this all together,
\begin{multline*}
    |A_2(m,n)|\leq c'\ind_{v_{\mf p}(m)\geq -1} \ind_{v_{\mf p}(n)\geq -1}\ind_{m,n\in (\cond{\psi})^{(\mf p)}}\ind_{n=m\lambda}\frac{1}{\sqrt{\Nm(\mf l_1\mf l_2)}}\\ \Bigg{(} \sum\limits_{ \substack{ l_1\neq  l_2\\ l_1, l_2\in \Lambda}} \Big| x_{ \mf l_1}\overline{x_{ \mf l_2}} \Big( \Nm(\mf l_2)\ind_{v_{\mf l_2}(m)\geq 1} \ind_{\lambda\in \frac{\mf l_1}{\mf l_2}\mcOt}+ \Nm(\mf l_1) \ind_{v_{\mf l_1}(m)\geq 1} \ind_{\lambda\in \frac{\mf l_2}{\mf l_1}\mcOt}\Big ) K(m_p) \overline{K(-n_p)} \Big| \Bigg{)}\times|\phi(m,\lambda)|. 
\end{multline*}

This implies 
\begin{multline*}
   \sum_{m\in F}\sum_{n\in F} |A_2(m,n)|\leq c' \sum\limits_{ \substack{ l_1\neq  l_2\\ l_1, l_2\in \Lambda}} \left |x_{ \mf l_1}\overline{x_{ \mf l_2}}\right |\ind_{\substack{[l_1]=[l_2]\\ [l_1],[l_2]\in \Cl(\mcO)}}\frac{1}{\sqrt{\Nm(\mf l_1\mf l_2)}}\\ \Big(\sum_{m\in \Ft } \sum_{u\in \mcOt }\ind_{v_{\mf p}(m)\geq -1} \ind_{v_{\mf l_2}(m)\geq 1} \ind_{m\in (\cond{\psi})^{(\mf p)}}  \Nm(\mf l_2)\left|K(m_p) \overline{K(-m_p\lambda_0u)}\phi(m,\lambda_0 u)\right| \\  + \sum_{m\in \Ft } \sum_{u\in \mcOt } \ind_{v_{\mf p}(m)\geq -1}\ind_{v_{\mf l_1}(m)\geq 1} \ind_{m\in (\cond{\psi})^{(\mf p)}}\Nm(\mf l_1) \left|K(m_p) \overline{K(-m_p\frac{u}{\lambda_0})}\phi(m,\frac{u}{\lambda_0})\right| \Big).   
\end{multline*}

The sum over units is finite and of absolutely bounded length (depending only on $F,f_{\infty}$). The $m$ sum is over the lattice $\frac{\mf l_2}{\mf p}\cond(\psi)$ in the first case and over $\frac{\mf l_1}{\mf p}\cond(\psi)$ in the second case. The result follows by applying the lattice point counting lemma \ref{lpc}.
\qed

\begin{prop}[Bounds for the DFI amplifier]\label{diagbounddfi}

$$\left|\sum_{m\in F^{\times}}\sum_{n\in F^{\times}}A_2(m,n)\right|\ll_{f_{\infty},\mf N}\norm{K}_{\infty}^2\frac{\Nm(\mf p)}{L}\sum\limits_{ \substack{\mf l_1\neq \mf l_2\\\mf l_1,\mf l_2\sim L\text{ prime}}} |x_{\mf l_1}\overline{x_{\mf l_2}}|+ \norm{K}_{\infty}^2\frac{\Nm(\mf p)}{L^2}\sum\limits_{ \substack{\mf l_1\neq \mf l_2\\\mf l_1,\mf l_2\sim L\text{ prime}}} |x_{\mf l_1^2}\overline{x_{\mf l_2^2}}|.$$

\end{prop}

\proof

We have,

\begin{multline*}
 \sum_{m,n} A_2(m,n)=\sum\limits_{ \substack{\mf l_1\neq \mf l_2\\\mf l_1,\mf l_2\in \Lambda}} x_{\mf l_1}\overline{x_{\mf l_2}} \sum_{m,n} c'\ind_{v_{\mf p}(m)\geq -1} \ind_{v_{\mf p}(n)\geq -1}\ind_{m,n\in \cond(\psi)^{(\mf p)}}\\ \Bigg{(}  \ind_{\lambda\in(\whO^{\times})^{(\mf l_1\mf l_2)}} K(m_p) \overline{K(-n_p)}\\  \prod_{\mf q|\mf l_1\mf l_2}\Big{(}\sum_{a=0}^{v_{\mf q}(\mf l_2)}\sum_{\alpha=v_{\mf q}(\mf l_2)-v_{\mf q}(\mf l_1)-a}^{v_{\mf q}(\mf l_2)-a}\ind_{v_{\mf q}(\lambda)=v_{\mf q}(\mf l_1)-v_{\mf q}(\mf l_2)+2\alpha}\ind_{v_{\mf q}(m)\geq M(\mf q, a,\alpha)}(\Nm(\mf q))^{-\frac{v_{\mf q}(\lambda)}{2}}\Big{)}\Bigg{)}\phi(m,\lambda)   
\end{multline*}

We first note that $c'$ is of size 1 and $$\Lambda=\{\mf l| \mf l\text{ prime and }\mf l\sim L \text{ coprime to }\mf N\mf p\}\cup \{\mf l^2| \mf l\text{ prime and }\mf l\sim L\text{ coprime to }\mf N\mf p\}.$$

We can rewrite this sum as follows since $n=m\lambda$:

\begin{multline*}
 \sum_{m,n} A_2(m,n)=\sum\limits_{ \substack{\mf l_1\neq \mf l_2\\\mf l_1,\mf l_2\in \Lambda}} x_{\mf l_1}\overline{x_{\mf l_2}} \sum_{m,\lambda} c'\ind_{v_{\mf p}(m)\geq -1} \ind_{m\in \cond(\psi)^{(\mf p)}}\\ \Bigg{(}  \ind_{\lambda\in(\whO^{\times})^{(\mf l_1\mf l_2)}} K(m_p) \overline{K(-n_p)}\\  \prod_{\mf q|\mf l_1\mf l_2}\Big{(}\sum_{a=0}^{v_{\mf q}(\mf l_2)}\sum_{\alpha=v_{\mf q}(\mf l_2)-v_{\mf q}(\mf l_1)-a}^{v_{\mf q}(\mf l_2)-a}\ind_{v_{\mf q}(\lambda)=v_{\mf q}(\mf l_1)-v_{\mf q}(\mf l_2)+2\alpha}\ind_{v_{\mf q}(m)\geq M(\mf q, a,\alpha)}(\Nm(\mf q))^{-\frac{v_{\mf q}(\lambda)}{2}}\Big{)}\Bigg{)}\phi(m,\lambda)   
\end{multline*}

Note that $m$ is summed over a lattice  weighted by a Schwartz function and the $\lambda$ is constrained by the cutoff function to lie in a compact subset of $\Ft_{\infty}$. In each summand the valuation of $\lambda$ gets fixed at all places so it is essentially a sum over units.

We will carry out bounds for each distinct pair of elements in $\Lambda$ separately:

\textbf{Case 1:} A pair of primes $\mf l_1, \mf l_2$ ($\mf l_1\neq \mf l_2$)

By the remark following proposition \ref{diagl}, the only surviving terms are those that satisfy $$v_{\mf l_1}(\lambda)+v_{\mf l_2}(\lambda)=0$$
\begin{multline*}
     \sum_{m,n} A_2(m,n)= c'\sum\limits_{ \substack{ l_1\neq  l_2\\ l_1, l_2\in \Lambda}}  x_{ \mf l_1}\overline{x_{ \mf l_2}}\sum_{m,\lambda} \ind_{v_{\mf p}(m)\geq -1} \ind_{m\in (\cond{\psi})^{(\mf p)}}\\ \Bigg{(}  \Big( \sqrt{\frac{\Nm(\mf l_2)}{\Nm(\mf l_1)}}\ind_{v_{\mf l_2}(m)\geq 1} \ind_{\lambda\in \frac{\mf l_1}{\mf l_2}\mcOt}+ \sqrt{\frac{\Nm(\mf l_1)}{\Nm(\mf l_2)}} \ind_{v_{\mf l_1}(m)\geq 1} \ind_{\lambda\in \frac{\mf l_2}{\mf l_1}\mcOt}\Big ) K(m_p) \overline{K(-n_p)}  \Bigg{)}\phi(m,\lambda). 
\end{multline*}

This implies 
\begin{multline*}
   \sum_{m\in F}\sum_{n\in F} |A_2(m,n)|\ll_{F,\mf N} \sum\limits_{ \substack{ l_1\neq  l_2\\ l_1, l_2\in \Lambda}} \left |x_{ \mf l_1}\overline{x_{ \mf l_2}}\right |\ind_{\substack{[l_1]=[l_2]\\ [l_1],[l_2]\in \Cl(\mcO)}}\frac{1}{\sqrt{\Nm(\mf l_1\mf l_2)}}\\ \Big(\sum_{m\in \Ft } \sum_{u\in \mcOt }\ind_{v_{\mf p}(m)\geq -1} \ind_{v_{\mf l_2}(m)\geq 1} \ind_{m\in (\cond{\psi})^{(\mf p)}}  \Nm(\mf l_2)\left|K(m_p) \overline{K(-m_p\lambda_0u)}\phi(m,\lambda_0 u)\right| \\  + \sum_{m\in \Ft } \sum_{u\in \mcOt } \ind_{v_{\mf p}(m)\geq -1}\ind_{v_{\mf l_1}(m)\geq 1} \ind_{m\in (\cond{\psi})^{(\mf p)}}\Nm(\mf l_1) \left|K(m_p) \overline{K(-m_p\frac{u}{\lambda_0})}\phi(m,\frac{u}{\lambda_0})\right| \Big).   
\end{multline*}

The sum over units is finite and of absolutely bounded length (depending only on $F,f_{\infty}$). The $m$ sum is over the lattice $\frac{\mf l_2}{\mf p}\cond(\psi)$ in the first case and over $\frac{\mf l_1}{\mf p}\cond(\psi)$ in the second case. The result follows by applying the lattice point counting lemma \ref{lpc}.

\textbf{Case 2:} A prime and a square of prime $\mf l_1, \mf l_2^2$ (The primes $\mf l_1$ and $\mf l_2$ may or may not be equal.)

By the remark following proposition \ref{diagl}, the only surviving terms are those that satisfy $$v_{\mf l_1}(\lambda)+v_{\mf l_2}(\lambda)=0.$$ The reader can check that none of the terms satisfy this condition and the contribution of this case is $0$.

\textbf{Case 3:} Two prime squares $\mf l_1^2, \mf l_2^2$ (The primes $\mf l_1\neq \mf l_2$.)

By the remark following proposition \ref{diagl}, the only surviving terms are those that satisfy $$v_{\mf l_1}(\lambda)+v_{\mf l_2}(\lambda)=0.$$ 
\begin{multline*}
     \sum_{m,n} A_2(m,n)= c'\sum\limits_{ \substack{ l_1\neq  l_2\\ l_1, l_2\text{ prime }\sim L}}  x_{ \mf l_1^2}\overline{x_{ \mf l_2^2}}\sum_{m,\lambda} \ind_{v_{\mf p}(m)\geq -1} \ind_{m\in (\cond{\psi})^{(\mf p)}}\\ \Bigg{(}  \Big( \frac{\Nm(\mf l_2)}{\Nm(\mf l_1)} \ind_{v_{\mf l_2}(m)\geq 2} \ind_{\lambda\in \frac{\mf l_1^2}{\mf l_2^2}\mcOt} + \ind_{v_{\mf l_1}(m)\geq 1}\ind_{v_{\mf l_1}(m)\geq 1} \ind_{\lambda\in \mcOt} + \frac{\Nm(\mf l_1)}{\Nm(\mf l_2)}\ind_{v_{\mf l_1}(m)\geq 2} \ind_{\lambda\in \frac{\mf l_2^2}{\mf l_1^2}\mcOt}\Big ) K(m_p) \overline{K(-n_p)}  \Bigg{)}\phi(m,\lambda). 
\end{multline*}

This implies 
\begin{multline*}
   \sum_{m\in F}\sum_{n\in F} |A_2(m,n)|\ll_{F,\mf N} \sum\limits_{ \substack{ l_1\neq  l_2\\ l_1, l_2\text{ prime }\sim L}}  \left |x_{ \mf l_1^2}\overline{x_{ \mf l_2^2}}\right |\ind_{\substack{[l_1^2]=[l_2^2]\\ [l_1],[l_2]\in \Cl(\mcO)}}\\ \Big(\sum_{m\in \Ft } \sum_{u\in \mcOt }\ind_{v_{\mf p}(m)\geq -1} \ind_{v_{\mf l_2}(m)\geq 2} \ind_{m\in (\cond{\psi})^{(\mf p)}}  \frac{\Nm(\mf l_2)}{\Nm(\mf l_1)}\left|K(m_p) \overline{K(-m_p\lambda_0u)}\phi(m,\lambda_0 u)\right| \\  + \sum_{m\in \Ft } \sum_{u\in \mcOt }\ind_{v_{\mf p}(m)\geq -1} \ind_{v_{\mf l_1}(m)\geq 2} \ind_{m\in (\cond{\psi})^{(\mf p)}}  \frac{\Nm(\mf l_1)}{\Nm(\mf l_2)}\left|K(m_p) \overline{K(-m_p\lambda_0u)}\phi(m,\lambda_0 u)\right|\\ + \sum_{m\in \Ft } \sum_{u\in \mcOt }\ind_{v_{\mf p}(m)\geq -1} \ind_{v_{\mf l_1}(m)\geq 1} \ind_{v_{\mf l_2}(m)\geq 1} \ind_{m\in (\cond{\psi})^{(\mf p)}}  \left|K(m_p) \overline{K(-m_p\lambda_0u)}\phi(m,\lambda_0 u)\right|\Big).   
\end{multline*}

The sum over units is finite and of absolutely bounded length (depending only on $F,f_{\infty}$). The $m$ sum is over the lattice $\frac{\mf l_2^2}{\mf p}\cond(\psi)$ in the first case, over $\frac{\mf l_1^2}{\mf p}\cond(\psi)$ in the second case and over $\frac{\mf l_1\mf l_2}{\mf p}\cond(\psi)$ in the third case. The result follows by applying the lattice point counting lemma \ref{lpc}.

\qed

\subsection{The non-diagonal contribution}
The next set of orbits is parametrized by $$\delta= \begin{pmatrix}
0 & \mu \\
1 & 0 
\end{pmatrix} $$ for $\mu\in F^{\times}$.

By explicit calculation, the stabilizer is
$$H_{\delta}(\Aa)={(e,e)}.$$

So all $\delta$'s are relevant.

 With  $\mu\in F^{\times}$ we have
   \begin{align*}    
   I_{\delta}&=\int\limits_{ N(\Aa)\times N(\Aa)}f(n(t_1)^{-1}\delta n(t_2))\psi(mt_1)\psi(nt_2)dn_1dn_2
   \\ &=\int\limits_{ \Aa\times \Aa}f\begin{pmatrix}
-t_1 & \mu-t_1t_2 \\
1 & t_2 
\end{pmatrix}\psi(mt_1)\psi(nt_2)dt_1dt_2.
\end{align*}

Like in the diagonal case, the orbital integral is factorizable. We write down the non-archimedean parts below:

If $f=h[1]*h[1]^*$, 

 $$I_{\delta,f}=\int\limits_{ \Aa_f\times \Aa_f}\int\limits_{\overline{G(\Aa_f)}}h[1]\begin{pmatrix}
-t_1' & \mu-t_1't_2' \\
1 & t_2' 
\end{pmatrix}g^{-1})\overline{h[1](g^{-1})}\psi(mt_1')\psi(nt_2')dt_1'dt_2'dg.$$

Considering the support of $h[1]$ and computing determinants, we conclude that there exists $c\in \Aft$ (unique upto multiplication by an element of $\whOt$) s.t. $$c^2\mu\in \whO^{\times}.$$ 

For any such choice of $c\in \Aft$, let us make a change of variables $t_1=ct_1'$ and $t_2=ct_2'$ to get

$$I_{\delta,f}=\frac{1}{||c||^2_{\Af}}\int\limits_{ \Aa_f\times \Aa_f}\int\limits_{\overline{G(\Aa_f)}}h[1]\begin{pmatrix}
-t_1 & c\mu-\frac{t_1t_2}{c} \\
c & t_2 
\end{pmatrix}g^{-1})\overline{h[1](g^{-1})}\psi(\frac{mt_1-nt_2}{c}) dt_1dt_2dg.$$

If $f=h[\mf l_1]*h[\mf l_2]^*$

 $$I_{\delta,f}=\int\limits_{ \Aa_f\times \Aa_f}\int\limits_{\overline{G(\Aa_f)}}h_f[\mf l_1]\begin{pmatrix}
-t_1' & \mu-t_1't_2' \\
1 & t_2' 
\end{pmatrix}g^{-1})\overline{h_f[\mf l_2](g^{-1})}\psi(mt_1')\psi(nt_2')dt_1'dt_2'dg.$$

Considering the support of $h_f[\mf l]$ and computing determinants, we conclude that there exists $c\in \Aft$ (unique upto multiplication by an element of $\whO^{\times}$) s.t. $$c^2\mu\in \frac{\mf l_1}{\mf l_2}\whO^{\times}.$$ 

For any such choice of $c\in \Aft$, let us make a change of variables $t_1=ct_1'$ and $t_2=ct_2'$ to get

$$I_{\delta,f}=\frac{1}{||c||^2_{\Af}}\int\limits_{ \Aa_f\times \Aa_f}\int\limits_{\overline{G(\Aa_f)}}h_f[\mf l_1]\begin{pmatrix}
-t_1 & c\mu-\frac{t_1t_2}{c} \\
c & t_2 
\end{pmatrix}g^{-1})\overline{h_f[\mf l_2](g^{-1})}\psi(\frac{mt_1-nt_2}{c}) dt_1dt_2dg.$$

Note that $c$ is an id\`ele as defined above, but in what follows we will also denote the components of $c$ by $c$ so as not to overload notation. The reader may be able to understand which local component is meant, from context. 

\subsubsection{Local computation for $\mf q\nmid \mf l_1 \mf l_2\mf p$} \label{nondiag1}
Let $S$ be the set of finite places coprime to $\mf l_1 \mf l_2\mf p$. Let us calculate:

$$\int\limits_{ \Aa_S\times \Aa_S}\int\limits_{\overline{G(\Aa_S)}}h_1\begin{pmatrix}
-t_1 & c\mu-\frac{t_1t_2}{c} \\
c & t_2 
\end{pmatrix}g^{-1})\overline{h_2(g^{-1})}\psi(\frac{mt_1-nt_2}{c})dt_1dt_2dg$$

where $(h_1,h_2)=(h_S[\mf l_1],h_S[\mf l_2])$ or $(h_1,h_2)=(h_{\mf p,S},h_{\mf p,S})$. Recall that $h_S=\ind_{K_0(\mf N)_S Z(\Aa_S)}$ in both these cases, we see that the above equals:

$$\meas(K_0(\mf N))\int\limits_{ \whO_S\times \whO_S}\ind_{c\in\whO_S}\ind_{c\in \mf N\whO_S}\ind_ {t_1t_2-c^2\mu \in c\whO_S} \psi(\frac{mt_1-nt_2}{c})dt_1dt_2.$$

Since the integrand is constant on cosets of $c\whO_s$, it follows that

$$\int\limits_{ \Aa_S\times \Aa_S}\int\limits_{\overline{G(\Aa_S)}}h_1\begin{pmatrix}
-t_1 & c\mu-\frac{t_1t_2}{c} \\
c & t_2 
\end{pmatrix}g^{-1})\overline{h_2(g^{-1})}\psi(\frac{mt_1-nt_2}{c})dt_1dt_2dg$$

$$=\meas(K_0(\mf N))(\meas(c\whO_s))^2\ind_{c\in \mf N\whO_S}\ind_{m\in\cond(\psi)_S}\ind_{n\in\cond(\psi)_S}\sum_{\substack{s_1,s_2\in\whO_s/c\whO_s \\ s_1s_2\equiv c^2\mu \mod c\whO_s}}\psi(\frac{ms_1-ns_2}{c})$$

\subsubsection{Local computation at $\mf p$}\label{nondiag2}

Let us calculate the local integral:
\begin{multline*}
 \frac{1}{\meas(K_0(\mf p))^2 \Nm(\mf p)} \sum_{a,b\in k(\mf p)}\whK(a)\overline{\whK(b)} \int\limits_{ F_{\mf p}\times F_{\mf p} }
 \int\limits_{\overline{G}(F_{\mf p})} \ind_{Z(F_{\mf p})\begin{pmatrix}
\un p & a \\
0 & 1 
\end{pmatrix}K_0(\mf p)} \begin{pmatrix}
-ct_1' & c\mu-ct_1't_2' \\
c & ct_2' 
\end{pmatrix}g^{-1} ) \\ \ind_{Z(F_{\mf p})\begin{pmatrix}
\un p & b \\
0 & 1 
\end{pmatrix}K_0(\mf p)} (g^{-1} ) \psi_{\mf p}(m t_1'+nt_2') dg dt_1'dt_2'.
\end{multline*}
We make the change of variables:

$$t_1=c(t_1'+a)$$ $$t_2=c(t_2'+b)$$ to get
\begin{multline*}
 \frac{1}{\meas(K_0(\mf p))^2 \Nm(\mf p)}\frac{1}{|c|_{\mf p}^2}\sum_{a,b\in k(\mf p)}\whK(a)\overline{\whK(b)}\psi_{\mf p}(-(ma+nb))\\ \int\limits_{ F_{\mf p}\times F_{\mf p} } \int\limits_{\overline{G}(F_{\mf p})} \ind_{Z(F_{\mf p})\begin{pmatrix}
\un p & a \\
0 & 1 
\end{pmatrix}K_0(\mf p)} \begin{pmatrix}
-t_1+ac & c\mu-c(\frac{t_1}{c}-a)(\frac{t_2}{c}-b) \\
c & t_2-bc 
\end{pmatrix}g^{-1} ) \\ \ind_{Z(F_{\mf p})\begin{pmatrix}
\un p & b \\
0 & 1 
\end{pmatrix}K_0(\mf p)} (g^{-1} ) \psi_{\mf p}(\frac{mt_1-nt_2}{c}) dg dt_1dt_2.
\end{multline*}
We want to understand when we have 
$$\begin{pmatrix}
-\un pt_1+ac\un p & c\mu-c(\frac{t_1}{c}-a)(\frac{t_2}{c}-b)-bt_1+abc \\
c\un p & t_2
\end{pmatrix}\in \begin{pmatrix}
\un p & a \\
0 & 1 
\end{pmatrix}K_0(\mf p).$$

This is equivalent to demanding:

$$\begin{pmatrix}
-t_1 & \frac{(c^2\mu-t_1t_2)}{c\un p} \\
c\un p & t_2
\end{pmatrix}\in K_0(\mf p).$$

Therefore
\begin{multline*} 
\ind_{Z(F_{\mf p})\begin{pmatrix}
\un p & a \\
0 & 1 
\end{pmatrix}K_0(\mf p)} \begin{pmatrix}
-t_1+ac & c\mu-c(\frac{t_1}{c}-a)(\frac{t_2}{c}-b) \\
c & t_2-bc 
\end{pmatrix}g^{-1}   \ind_{Z(F_{\mf p})\begin{pmatrix}
\un p & b \\
0 & 1 
\end{pmatrix}K_0(\mf p)} (g^{-1} ) \\ =\ind_{v_{\mf p}(c)\geq 0}\ind_{\mcO_{\mf p}}(t_1)\ind_{\mcO_{\mf p}}(t_2)\ind_{t_1t_2- c^2\mu \in c\un p } \ind_{Z(F_{\mf p})K_0(\mf p)\begin{pmatrix}
\un p & b \\
0 & 1 
\end{pmatrix}^{-1}} (g ).
\end{multline*} 
The integral evaluates to 
\begin{multline*}
 \frac{1}{\meas(K_0(\mf p)) \Nm(\mf p)}\frac{1}{|c|_{\mf p}^2}\ind_{v_{\mf p}(c)\geq 0} \sum_{a,b\in \Fp}\whK(a)\overline{\whK(b)}\psi_{\mf p}(-(ma+nb))\\ \int\limits_{ \mcO_{\mf p}\times \mcO_{\mf p} }  \ind_{t_1t_2\equiv c^2\mu \mod c\un p}   \psi_{\mf p}(\frac{mt_1-nt_2}{c})  dt_1dt_2.
\end{multline*}
The integrand is constant on cosets modulo $c\un p$, so we get
\begin{multline*}
\ind_{v_{\mf p}(c)\geq 0} \frac{1}{\meas(K_0(\mf p)) \Nm(\mf p)}\sum_{a,b\in \Fp}\whK(a)\overline{\whK(b)}\psi_{\mf p}(-(ma+nb))\\
\sum_{\substack{s_1,s_2\in\mcO_{\mf p}/c\un p\mcO_{\mf p} \\ s_1s_2\equiv c^2\mu}}\psi(\frac{ms_1-ns_2}{c})\int\limits_{ \un p\mcO_{\mf p}\times \un p\mcO_{\mf p} }    \psi_{\mf p}(mt_1-nt_2)  dt_1dt_2.
\end{multline*}
The local integral at $\mf p$ is
\begin{multline*}
    \ind_{v_{\mf p}(c)\geq 0} \frac{1}{\meas(K_0(\mf p)) (\Nm(\mf p))^2}\ind_{v_{\mf p}(m)\geq -1}\ind_{v_{\mf p}(n)\geq -1}K(m_{\mf p})\overline{K(-n_{\mf p})}\\
    \sum_{\substack{s_1,s_2\in\mcO_{\mf p}/c\un p\mcO_{\mf p} \\ s_1s_2\equiv c^2\mu}}\psi_{\mf p}(\frac{m\un ps_1-n\un ps_2}{c\un p}).
\end{multline*}

Recall our assumption that $\mf p$ is large enough so that $\cond(\psi)$ and $\mf p$ are coprime.

\subsubsection{Local computation at $\mf q|\mf l_1\mf l_2$, $\mf l_1\neq \mf l_2$}\label{nondiag4}

 Let $v_{\mf q}(\mf l_1)=d_1$ and $v_{\mf q}(\mf l_2)=d_2$. Looking at the determinants we can choose $c\in \Aft$ in such a way that $c^2_{\mf q}\mu\in \mf q^{d_1-d_2}\mcO^{\times}_{\mf q}$. We need to calculate the local integral:
$$\frac{1}{\sqrt{\Nm(\mf l_1)}}\int\limits_{ F_{\mf q}\times F_{\mf q} } \int\limits_{\overline{G(F_{\mf q})}} \ind_{Z(F_{\mf q})M(\mf q^{d_1})_{\mf q}} \begin{pmatrix}
-ct_1 & c\mu-ct_1t_2 \\
c & ct_2 
\end{pmatrix}g^{-1} \ind_{Z(F_{\mf q}) M(\mf q^{d_2})_{\mf q}} (g^{-1} ) \psi_{\mf q}(m t_1-n t_2) dg dt_1dt_2$$
where $$M(\mf l)_{\mf q}= \{M\in M_2(\mcO_{\mf q})|\text{det}(M)\in \mf l\mcO_{\mf q}^{\times}\}.$$

$$=\frac{1}{\sqrt{\Nm(\mf q^{d_1+d_2})}}\sum_{a=0}^{d_2}\sum_{b \text{ mod }\mf q^a}\int\limits_{ F_{\mf q}\times F_{\mf q} }  \ind_{M(\mf q^{d_1})_{\mf q}} \left(\begin{pmatrix}
-ct_1 & c\mu-ct_1t_2 \\
c & ct_2 
\end{pmatrix}\begin{pmatrix}
\un q^a & b \\
0 & \un q^{d_2-a} 
\end{pmatrix}\right )  \psi_{\mf q}(m t_1-n t_2)  dt_1dt_2.$$
Changing variables $t_1$ by $ct_1\un q^a$ and $t_2$ by $bc+ct_2\un q^{d_2-a}$ we need to look at
$$\frac{1}{|c|^2|\un q|^{d_2}}\frac{1}{\sqrt{\Nm(\mf q^{d_1+d_2})}}\sum_{a=0}^{d_2}\sum_{b \text{ mod }\mf q^a}\psi_{\mf q}\left(\frac{-nb}{\un q^{d_2-a}}\right)\int\limits_{ F_{\mf q}\times F_{\mf q} }  \ind_{M(\mf q^{d_1})_{\mf q}} \begin{pmatrix}
-t_1 & c\mu\un q^{d_2-a}-\frac{t_1t_2}{c\un q^a} \\
c\un q^a & t_2 
\end{pmatrix}  \psi_{\mf q}\left(\frac{m t_1-\frac{n\un q^a}{\un q^{d_2-a}} t_2}{c\un q^a}\right)  dt_1dt_2$$
$$=\frac{1}{|c|^2|\un q|^{d_2}}\frac{1}{\sqrt{\Nm(\mf q^{d_1+d_2})}}\sum_{a=0}^{d_2}\ind_{n \in \un q^{d_2-a}\cond(\psi_{\mf q})}\int\limits_{ \mcO_{\mf q}\times \mcO_{\mf q} } \ind_{c^2\mu\un q^{d_2}- t_1t_2 \in c\un q^a\mcO_q} \ind_{c\un q^a\in \mcO_{\mf q}}   \psi_{\mf q}\left(\frac{m t_1-\frac{n\un q^a}{\un q^{d_2-a}} t_2}{c\un q^a}\right)  dt_1dt_2$$
$$=\frac{1}{|\un q|^{d_2}}\frac{1}{\sqrt{\Nm(\mf q^{d_1+d_2})}}\sum_{a=0}^{d_2}\ind_{n \in \un q^{d_2-a}\cond(\psi_{\mf q})}\ind_{m\in \cond(\psi_{\mf q})}\ind_{c\un q^a\in \mcO_{\mf q}}  \sum_{\substack{s_1,s_2 \mod c\un q^a\\ s_1s_2\equiv c^2\mu\un q^{d_2} \mod c\un q^a\mcO_{\mf q}}} \psi_{\mf q}\left(\frac{m s_1-\frac{n\un q^a}{\un q^{d_2-a}} s_2}{c\un q^a}\right).$$

The local integral at this place $\mf q$ is 

$$\frac{1}{|\un q|^{d_2}}\frac{1}{\sqrt{\Nm(\mf q^{d_1+d_2})}}\sum_{a=0}^{d_2}\ind_{n \in \un q^{d_2-a}\cond(\psi_{\mf q})}\ind_{m\in \cond(\psi_{\mf q})}\ind_{c\un q^a\in \mcO_{\mf q}}  \sum_{\substack{s_1,s_2 \mod c\un q^a\\ s_1s_2\equiv c^2\mu\un q^{d_2} \mod c\un q^a\mcO_{\mf q}}} \psi_{\mf q}\left(\frac{m s_1-\frac{n\un q^a}{\un q^{d_2-a}} s_2}{c\un q^a}\right).$$

\subsubsection{Archimedean local computation} \label{nondiagarch}
The archimedean local component is given by

$$\phi(m,n,\mu)= \int\limits_{ F_{\infty}\times F_{\infty}}f_{\infty}[\mf l_1, \mf l_2]\begin{pmatrix}
-t_1 & \mu-t_1t_2 \\
1 & t_2 
\end{pmatrix}\psi(mt_1)\psi(nt_2)dt_1dt_2.$$

Since $f_{\infty}[\mf l_1, \mf l_2]$ is compactly supported on $\PGL_2(F_{\infty})$, we conclude that $\phi(m,n,\mu)=0$ except for $\mu$ that satisfies

$$|\mu|_v\gg_{f,F}1\text{ for every }v|\infty.$$

Further since $f_{\infty}[\mf l_1, \mf l_2]$ is smooth, $\phi(m,n,\mu)$ is a Schwartz function on $m,n$.

For later computations we will also need to compute the Fourier transform. It has the following simple shape

$$\hat{\phi}(m,n,\mu)= f_{\infty}[\mf l_1, \mf l_2]\begin{pmatrix}
m & \mu-mn \\
1 & -n 
\end{pmatrix}.$$

\subsection{The non-diagonal contribution and bounds}
First, we would like to define some Kloostermann sums: Let $S$ be a set of finite places and let $a,b,c,d\in \whO_S$. Recall that $\whO_S:=\prod_{\mf q\in S}\mcO_{\mf q}.$ 

We denote $$\Kl_S(a,b,d;c)=\sum_{\substack{s_1,s_2 \in\whO_S/c\whO_S\\ s_1s_2\equiv d}} \psi_S( \frac{as_1+bs_2}{c})$$

and by  $$\Kl_S(a,b;c)=\sum_{\substack{s_1,s_2 \in\whO_S/c\whO_S\\ s_1s_2\equiv 1}} \psi_S( \frac{as_1+bs_2}{c}).$$

We will write $\Kl(a,b,d;c)$ instead of $\Kl_S(a,b,d;c)$ and $\Kl(a,b;c)$ instead of $\Kl_S(a,b;c)$ when the set of places is clear from context, in order to simplify notation.

Note if $d\in\whOt_S$, then $\Kl(a,b,d;c)=\Kl(a,db;c)=\Kl(ad,b;c).$

 \subsubsection{Non-diagonal contribution for $f=h[1]*h[1]^*$ }
  Let $$C_0(\mf p)=\frac{1}{\meas(K_0(\mf p)) \meas(K_0(\mf N))(\Nm(\mf p))^2}.$$ We will drop the dependence on $\mf p$ in the sequel for simplicity.

Combining the results of sections \S\ref{nondiag1} , \S\ref{nondiag2} we get the following proposition
\begin{prop}
The non-diagonal contribution to the geometric side is given by
\begin{multline*}  
B(m,n)=C_0 \sum\limits_{\mu \in \mathcal{C}}\ind_{c\in \mf N\whO}\ind_{m\un p\in\whO}\ind_{n\un p\in\whO} \Kl(m\un p,n\un p,c_{\mf p}^2\mu;c_{\mf p}\un p) \\ \Kl(m,n,c_{S}^2\mu;c_{S})K(m_{\mf p})\overline{K(-n_{\mf p})} \phi(m,n,\mu) 
\end{multline*}

where $$\mathcal{C}=F^{\times}\cap (\Aft)^2\whOt.$$  This is exactly the set of rational elements that have an even valuation at all the places.
$c\in \Aft$ is any idele satisfying $$c^2\mu\in\whOt.$$
The choice of $c$ for a given $\mu$ is not unique. Two possible choices of $c$ differ by multiplying by an element of $\whOt$. Each summand however is independent of this choice. 

\end{prop}

\begin{obs}\label{obs1}
As discussed in \S\ref{nondiagarch} $|\mu|_{\infty}\gg_{f,F} 1$ , $$|c|_{\Af}\gg_{f,F}1.$$ The choice of $c$ is only up to multiplication by an element of $\whOt$, so we may think of the ideal that corresponds to $c$, and the above inequality implies the norm of that ideal is bounded.

Observe that in the sum the $c$ that appears has to be integral. We have in conclusion $$1\ll_{f,F} |c|_{\Af}\ll 1.$$

So for sufficiently large $\Nm(\mf p)$ , we have $$|c|_{\mf p}=1. $$ 
\end{obs}

We wish to analyse and bound the sum $$\sum_{m,n\in F}B(m,n) $$ using the adelic Poisson summation formula. Note that $B(m,n)$ is indeed a Schwartz-Bruhat function of $m,n$. 
Recall that $$C_0=\frac{1}{\meas(K_0(\mf p)) \meas(K_0(\mf N))(\Nm(\mf p))^2}.$$
\begin{prop}
The adelic Fourier transform of $B(m,n)$ is given by \begin{multline*}
\wh{B}(m,n)
=C_0 . \Nm(\mf p)\sum\limits_{\mu \in \mathcal{C}}\sum_{\substack{t_1t_2\equiv c^2\mu \\ \mod c\un p}} \ind_{c\in \mf N\whO} \ind_{m+\frac{t_1}{c}\in \cond(\psi)}\ind_{n+\frac{t_2}{c}\in \cond(\psi)} \\ \whK\bigg(m+\frac{t_1}{c}\bigg)\ov{\whK\bigg(n+\frac{t_2}{c}\bigg)}
\wh{\phi}(m,n,\mu).   \end{multline*}

\end{prop}
\begin{prop} \label{nondiagbound}
      For $\Nm(\mf p)$ large enough as in \ref{obs1}, we have the following bound 
       $$\left|\sum_{m,n}\wh{B}(m,n)\right|\ll_{f_{\infty}} \cond(K)^4 \Nm(\mf p).$$
    where $\cond(K)$ is the smallest conductor among all sheaves whose trace function is $K$.
    \end{prop}
    
\proof

First note that we have 

$$C_0 . \Nm(\mf p)\ll_{F, \mf N} 1.$$

Since $K$ is a trace function which is not an additive character, $\whK$ is bounded independent of $\Nm(\mf p)$. (Recall that the Fourier transform on $k(\mf p)$ is unitarily normalized.) So bounding the summands trivially, we get
\begin{multline*}
    \left|\sum_{m,n}\wh{B}(m,n)\right|=\bigg |C_0 . \Nm(\mf p) \sum\limits_{m,n\in \Ft}\sum\limits_{\mu \in \mathcal{C}}\sum\limits_{\substack{t_1t_2\equiv c^2\mu \\ \mod c\un p}} \ind_{c\in \mf N\whO}\ind_{m+\frac{t_1}{c}\in \cond(\psi)}\ind_{n+\frac{t_2}{c}\in \cond(\psi)} \\ \whK\bigg(m+\frac{t_1}{c}\bigg)\ov{\whK\bigg(n+\frac{t_2}{c}\bigg)}
\wh{\phi}(m,n,\mu) \bigg | \\ \ll_{F} \norm{\wh{K}}_{\infty}^2\sum\limits_{\mu \in \mathcal{C}}\ind_{c\in \mf N\whO} \sum\limits_{\substack{t_1t_2\equiv c^2\mu \\ \mod c\un p}} \sum_{m,n} \left |\ind_{m+\frac{t_1}{c}\in \cond(\psi)}\ind_{n+\frac{t_2}{c}\in \cond(\psi)} 
\wh{\phi}(m,n,\mu) \right |.
\end{multline*}
The innermost sum in terms of $m,n$ is absolutely bounded uniformly in $t_1,t_2$. The bound depends on the test function $f_{\infty}$. So we get 

$$    
 \left|\sum_{m,n}\wh{B}(m,n)\right| \ll_{F, f_{\infty}} \norm{\wh{K}}_{\infty}^2\sum\limits_{\mu \in \mathcal{C}} \ind_{c\in \mf N\whO} \sum\limits_{\substack{t_1t_2\equiv c^2\mu \\ \mod c\un p}} 1.  $$

Proceeding further,

$$  
 \left|\sum_{m,n}\wh{B}(m,n)\right|\ll_{F, f_{\infty}} \norm{\wh{K}}_{\infty}^2|c|_{\Af}^{-1}\Nm(\mf p) \ll_{F, f_{\infty}} \norm{\wh{K}}_{\infty}^2 \Nm(\mf p)$$

In the first inequality we use the fact that the $\mu$-sum is of absolutely bounded length (depending on the test function). From observation \ref{obs1} we know that $|c|_{\Af}\gg 1$ so the last inequality follows. Finally since we have assumed the sheaf $\mcF$ underlying the trace function $K$ is Fourier, the correlation sum of the trace function with the additive character satisfies square root cancellation i.e.

$$\norm{\wh{K}}_{\infty}\ll \cond(K)^2.$$

We refer the reader to section \S \ref{corrsum}.

\subsubsection{Non-diagonal contribution for $f=\sum\limits_{\mf l_1 \neq \mf l_2} x_{\mf l_1}\ov{x_{\mf l_2}}h[\mf l_1]*h[\mf l_2]^*$}
Let $S[\mf l_1,\mf l_2]$ denote the set of finite places not dividing $\mf l_1\mf l_2$ and $S'[\mf l_1,\mf l_2]$ denote the set of finite places not dividing $\mf l_1\mf l_2\mf p$. For every $\mf l_1,\mf l_2$ and $\mu$ we choose $c\in\Aft$ s.t. $c^2\mu\in \frac{\un {l_1}}{\un {l_2}}\whOt$. We don't denote explicitly the dependence of c on $\mf l_1,\mf l_2,\mu$ to not overload the notation.
\begin{prop}
The non-diagonal contribution to the geometric side is given by
\begin{multline*}
B_2(m,n)=C_0 . \Nm(\mf p)\frac{1}{\sqrt{\Nm(\mf l_2)\Nm(\mf l_1)}}\sum\limits_{ \substack{\mf l_1\neq \mf l_2\\\mf l_1,\mf l_2\in \Lambda}}\sum\limits_{\mu \in \mathcal{C}_{\mf l_1,\mf l_2}}\ind_{c_{S[\mf l_1,\mf l_2]}\in \mf N\whO_{S[\mf l_1,\mf l_2]}}\ind_{m\un p\in\cond(\psi)\whO}\ind_{n\un p\in\cond(\psi)\whO^{(\mf l_1\mf l_2)}}\\ \Kl(m\un p,n\un p,c_{\mf p}^2\mu;c_{\mf p}\un p)\Kl(m,n,c_{S'[\mf l_1,\mf l_2]}^2\mu;c_{S'[\mf l_1,\mf l_2]})K(m_{\mf p})\overline{K(-n_{\mf p})}\\\prod_{\mf q|\mf l_1\mf l_2} \left(\sum_{a=0}^{v_{\mf q}(\mf l_2)}\Nm(\mf q)^{2a-v_{\mf q}(\mf l_2)}\ind_{n \in \un q^{v_{\mf q}(\mf l_2)-a}\cond(\psi_{\mf q})}\ind_{c\un q^a\in \mcO_{\mf q}}\Kl(m,n\un q^{2a-v_{\mf q}(\mf l_2)},c^2\mu\un q^{v_{\mf q}(\mf l_2)},c\un q^a)\right)\phi(m,n,\mu).\end{multline*} 

Here $$\mathcal{C}_{\mf l_1,\mf l_2}=F^{\times}\cap (\Aft)^2\frac{\un {l_1}}{\un {l_2}}\whOt.$$ 
For $\mu\in \mathcal{C}_{\mf l_1,\mf l_2}$, $c\in \Aft$ is any idele satisfying $$c^2\mu\in\frac{\un {l_1}}{\un {l_2}}\whOt.$$

The choice of $c$ for a given $\mu$ is not unique. Two possible choices of $c$ differ by multiplying by an element of $\whOt$. Each summand however is independent of this choice. 

\end{prop}

\begin{obs}[Archimedean constraint for the Venkatesh amplifier]\label{obsv}
Since $|\mu|_{\infty}\gg_{f_{\infty},F} 1$  and elements in $\Lambda$ are of the same size, we have from $c^2\mu\in\frac{\un {l_1}}{\un {l_2}}\whOt$, 
$$|c|_{\Af}\gg_{f_{\infty},F}1.$$ The choice of $c$ is only up to multiplication by an element of $\whOt$, so we may think of the ideal that corresponds to $c$, and the above inequality implies the norm of that ideal is bounded.

Observe that in the sum $c\un {l_2}$ has to be integral. We have in conclusion $$1\ll_{f,F} |c|_{\Af}\ll L$$ or more precisely $$1\ll_{f,F} |c|_{\Af}$$ and $$|c\un {l_2}|_{\Af}\leq 1.$$ 

 We will take all the places in $\Lambda$ to have the same size $L$ (i.e. $L\leq \Nm(\mf l) \leq 2L$) and we will let $\Nm(\mf p)\xrightarrow[]{}\infty$  and places in $\Lambda$ go to infinity s.t. $\dfrac{\Nm(\mf p)}{L}\xrightarrow[]{}\infty$ .

So for sufficiently large (depending on $f_{\infty}$) $\mf p$ and $ L$, we have $$|c_{\mf p}|_{\mf p}=1 \text{ and  } |c_{\mf l_2}|_{\mf l_2}\geq 1 $$ in the sum. 

\end{obs}
\begin{obs}[Archimedean constraint for the DFI amplifier case]\label{obsdfi} 

Now the elements of $\Lambda$ are no longer the same size.

Since $|\mu|_{\infty}\gg_{f_{\infty},F} 1$  we have $$|c|_{\Af}\gg_{f_{\infty},F}\sqrt{\left|\frac{\mf l_1}{\mf l_2}\right|_{\Af}}.$$ 

We have 

$$|c|_{\Af}\gg_{f_{\infty},F}  
     \begin{cases}
       1 &\text{if $\mf l_1,\mf l_2\in \Lambda$ are both primes or squares of primes} \\
        \frac{1}{\sqrt{L}}&\quad\text{if $\mf l_1$ is a square and $\mf l_2$ is prime}
       \\
         \sqrt{L}&\quad\text{if $\mf l_1$ is a prime and $\mf l_2$ is square}
       \end{cases}
$$
The choice of $c$ is only up to multiplication by an element of $\whOt$, so we may think of the ideal that corresponds to $c$, and the above inequality implies the norm of that ideal is bounded.

Observe that in the summand corresponding to a particular value of $(a_{i})_{i\in \Spec(\mf l_2)}$ $$c\prod_{i\in \Spec(\mf l_2)}\un {q}^{a_i}\in \whO.$$ We have in conclusion the ideal $J$ corresponding to the idele $c\prod_{i\in \Spec(\mf l_2)}\un {q}^{a_i}$ is integral and satisfies

$$|\Nm(J)|\ll_{f_{\infty},F}  
     \begin{cases}
       \Nm(\mf l) &\text{if $\mf l_1,\mf l_2\in \Lambda$ are both primes or squares of primes} \\
       \Nm(\mf l)\sqrt{L}  &\quad\text{if $\mf l_1$ is a square and $\mf l_2$ is prime}
       \\
        \frac{\Nm(\mf l)}{\sqrt{L}} &\quad\text{if $\mf l_1$ is a prime and $\mf l_2$ is square}
       \end{cases}
$$ 

where $\mf l=\prod_{\mf q\in \Spec(\mf l_2)}(\mf q)^{a_{\mf q}}$

 We will let $\Nm(\mf p)\xrightarrow[]{}\infty$  and places in $\Lambda$ go to infinity s.t. $\dfrac{\Nm(\mf p)}{L^2}\xrightarrow[]{}\infty$ .

So for sufficiently large (depending on $f_{\infty}$) $\mf p$ and $ L$, we have $$|c_{\mf p}|_{\mf p}=1$$ in the sum. 

\end{obs}
We wish to now analyse and bound the sum  $$\sum_{m,n\in F}B_2(m,n)$$
    using the adelic Poisson summation formula. Note that  $B_2(m,n)$ is indeed a Schwartz-Bruhat function of $m,n$. 
    \begin{prop}
    The adelic Fourier transform of $B_2(m,n)$ is given by
\begin{multline*}
 \wh{B_2}(m,n)=C_0 . \Nm(\mf p)\sum\limits_{ \substack{\mf l_1\neq\mf  l_2\\\mf l_1,\mf l_2\in \Lambda}}\frac{1}{\sqrt{\Nm(\mf l_1)\Nm(\mf l_2)}}x_{\mf l_1}\ov{x_{\mf l_2}}\sum\limits_{\mu \in \mathcal{C}_{\mf l_1,\mf l_2}}\ind_{ c\in \mf N\whO^{(\mf l_1\mf l_2)}} \ind_{ cm \in \whO^{(\mf l_1\mf l_2)}}\ind_{cn \in \whO^{(\mf l_1\mf l_2)}}\\
 \prod_{\mf q|\mf l_1\mf l_2}\Bigg(\sum_{a=0}^{v_{\mf q}(\mf l_2)} \ind_{\un {q}^a c\in \mf N\mcO_{\mf q}} \ind_{\un {q}^a cm \in \mcO_{\mf q}}\ind_{\un {q}^{v_{\mf q}(\mf l_2)}cn \in \mcO_{\mf q}} \\ 
 \sum_{\substack{s_1,s_2\text{ mod } c\un q^a\\ s_1s_2 \equiv c^2 \un q^{v_{\mf q}(\mf l_2)}\mu }}\ind_{\un {q}^a cm\equiv s_1 \text{ mod } (c\un q^a \cond(\psi))}\ind_{\un {q}^{v_{\mf q}(\mf l_2)} cn\equiv \un q^a s_2 \text{ mod } (c\un q^a \cond(\psi))}\Bigg)\\ \Cor(K,\gamma_{m,n}(\mu)) \wh{\phi}(m,n,\mu).\end{multline*}
  
    \end{prop}
Note that $$\Cor(K,\gamma_{m,n}(\mu))=\sum_{z\in k(\mf p)}\wh{K}(\gamma_{m,n}(\mu)z)\ov{\wh{K}(z)}$$ is the correlation sum with $$\gamma_{m,n}(\mu)=\begin{pmatrix}
m & \mu-mn \\
1 & -n 
\end{pmatrix}\in \PGL_2(k(\mf p)).$$

Recall that 
$$\wh{\phi}(m,n,\mu)=f_{\infty}\begin{pmatrix}
m & \mu-mn \\
1 & -n 
\end{pmatrix}.$$

\begin{prop}[\textbf{Bounds for the non-diagonal contribution for the Venkatesh amplifier}] \label{nondiagboundv} 
     Let us assume that the trace function $K$ has Fourier-M\"obius group contained in the standard Borel subgroup. We have the following bound for $L$ and $\Nm(\mf p)$ large enough (as in \ref{obsv})
       $$\left|\sum_{m,n\in F}\wh{B_2}(m,n)\right|\ll_{f_{\infty}}\cond(K)^4 L^{1+\text{o}(1)}\sqrt{\Nm(\mf p)}\sum\limits_{ \substack{\mf l_1\neq \mf  l_2\\\mf l_1,\mf l_2\in \Lambda}}|x_{\mf l_1}\ov{x_{\mf l_2}}| $$
       where $\cond(K)$ is the smallest conductor of a sheaf $\mcF$ whose trace function is $K$.
\end{prop}
    
\proof
For the Venkatesh amplifier, the adelic Fourier transform of $B_2(m,n)$ is given by
\begin{multline*}
 \wh{B_2}(m,n)=C_0 . \frac{\Nm(\mf p)}{\sqrt{\Nm(\mf l_1\mf l_2)}}\sum\limits_{ \substack{\mf l_1\neq \mf  l_2\\\mf l_1,\mf l_2\in \Lambda}}x_{\mf l_1}\ov{x_{\mf l_2}}  \sum\limits_{\mu \in \mathcal{C}_{\mf l_1,\mf l_2}}\Bigg(\ind_{\un {l_2}c\in \mf N\whO}\ind_{\un {l_2}cm\in\whO}\ind_{\un {l_2}cn\in\whO}\\ \Cor(K,\gamma_{m,n}(\mu))\ind_{mn\in \mu+\frac{1}{\un {l_2}c}\whO} \wh{\phi}(m,n,\mu)\Bigg).\end{multline*}
 
Recall that
$$C_0 . \Nm(\mf p)\ll_{F, \mf N} 1.$$

We have 
\begin{multline*}
\sum_{m,n\in F}|B_2(m,n)|= \Bigg|C_0 . \frac{\Nm(\mf p)}{\sqrt{\Nm(\mf l_1\mf l_2)}}\sum_{m,n\in F}\sum\limits_{ \substack{\mf l_1\neq \mf  l_2\\ \mf l_1,\mf l_2\in \Lambda}}x_{\mf l_1}\ov{x_{\mf l_2}}\sum_{\mu \in \mathcal{C}_{\mf l_1,\mf l_2}} \ind_{\un {l_2}c\in \mf N\whO} \ind_{\un {l_2}cm\in\whO} \ind_{\un {l_2}cn\in\whO}\\ \ind_{mn\in \mu+\frac{1}{\un {l_2}c}\whO} \Cor(K,\gamma_{m,n}(\mu))  \wh{\phi}(m,n,\mu)\Bigg|
\end{multline*}

Since $\gamma_{m,n}(\mu)$ is clearly not contained in the standard Borel subgroup and hence in the automorphism group attached to $\wh{K}$, using corollary \ref{corrsumsize}, we have square root cancellation, $$|\Cor(K,\gamma_{m,n}(\mu))|\leq (\cond(K))^4\sqrt{\Nm(\mf p)} .$$
\begin{multline*}
   \sum_{m,n\in F}|B_2(m,n)| \ll_{\cond(K)} \frac{\sqrt{\Nm(\mf p)}}{\sqrt{\Nm(\mf l_1\mf l_2)}} \sum\limits_{ \substack{\mf l_1\neq \mf  l_2\\\mf l_1,\mf l_2\in \Lambda}}|x_{\mf l_1}\ov{x_{\mf l_2}}|\sum_{\mu \in \mathcal{C}_{\mf l_1,\mf l_2}}\ind_{\un {l_2}c\in \mf N\whO}  \sum_{m,n\in F}\ind_{m\in I}\ind_{n\in I} \ind_{mn\in \mu+I} \left |\wh{\phi}(m,n,\mu)\right| 
\end{multline*}

where $I:=F\cap \frac{1}{\un {l_2}c}\whO$ is a fractional ideal satisfying $$\Nm(I)= |\un {l_2}c|_{\Af}\gg_{F,f_{\infty}} L^{-1}.$$

We may think of the $\mu$-sum as a sum over ideals $J$ s.t. $\frac{\mf l_2}{\mf l_1}J^2$ is a principal ideal. This is by defining the ideal $J$ to satisfy
$$(\mu)= \frac{\mf l_2}{\mf l_1}.J^2.$$

Moreover we saw in observation \ref{obsv} that $\wh{\phi}(m,n,\mu)\neq 0$ only if $|\Nm(J)|\ll_{f,F}1$ The non-archimedean condition implies $J\mf l_2$ is integral. Let us write $J=\frac{1}{\mf l_2}J'$ with $J'$ integral. Since $\mf l_2$ is of size $L$, we have $\wh{\phi}(m,n,\mu)\neq 0$ only if $|\Nm(J')|\ll_{f,F} L$. So the $\mu$-sum is finite and of  length $L$. Note that we have
$$(\mu)= \frac{1}{\mf l_2 \mf l_1}.(J')^2$$
and
$$I=\frac{J'}{\mf l_1\mf l_2}$$
under the new parametrization.

Also recall that 
$$\wh{\phi}(m,n,\mu)=f_{\infty}\begin{pmatrix}
m & \mu-mn \\
1 & -n 
\end{pmatrix}.$$

Hence (since $f_{\infty}$ is compactly supported modulo the center)

$$ |\wh{\phi}(m,n,\mu)| \ll_{f_{\infty}} \ind_{|m|_{\infty}\ll \sqrt{|\mu|_{\infty}}}\ind_{|n|_{\infty}\ll \sqrt{|\mu|_{\infty}}}\ind_{|mn-\mu|_{\infty}\ll \sqrt{|\mu|_{\infty}}}.$$

The inequalities in the indicators also depend on $f_{\infty}$ or more precisely its support.

 Now we bound each term of the $\mu$-sum

$$
    \sum_{m,n\in F}\ind_{m\in I}\ind_{n\in I} \ind_{mn\in \mu+I} \left |\wh{\phi}(m,n,\mu)\right|  \ll \sum_{|k-\mu|_{\infty}\ll \sqrt{|\mu|_{\infty}}}  \ind_{k\in \mu +I} \sum_{\substack{m,n\in I\\|m|_{\infty}\ll \sqrt{|\mu|_{\infty}}\\|n|_{\infty}\ll \sqrt{|\mu|_{\infty}}}}\ind_{mn=k}.
$$

We bound the innermost sum using the divisor counting lemma (lemma \ref{divc}). Also note that $\Nm(k)\ll \Nm(\mu)$.

$$
\sum_{|k-\mu|_{\infty}\ll \sqrt{|\mu|_{\infty}}}  \ind_{k\in \mu +I} \sum_{\substack{m,n\in I\\|m|_{\infty}\ll \sqrt{|\mu|_{\infty}}\\|n|_{\infty}\ll \sqrt{|\mu|_{\infty}}}}\ind_{mn=k} \ll (\Nm(\mu))^{\text{o}(1)} \sum_{|k-\mu|_{\infty}\ll \sqrt{|\mu|_{\infty}}}  \ind_{k\in \mu +I}.
$$
  By the lattice point counting lemma \ref{lpc}  we have
\begin{multline*}
     \sum_{m,n\in F}\ind_{m\in I}\ind_{n\in I} \ind_{mn\in \mu+I} \left |\wh{\phi}(m,n,\mu)\right|  \ll (\Nm(\mu))^{\text{o}(1)} \frac{(\Nm(\mu))^{1/2}}{\Nm(I)}\\=(\Nm(\mu))^{\text{o}(1)}\sqrt{\Nm(\mf l_1 \mf l_2)}\ll L(\Nm(\mu))^{\text{o}(1)}.
\end{multline*}
 
The last equality follows by substituting the expression for $(\mu)$ and $I$ in terms of $J'$ that we noted earlier.

Now rewriting the $\mu$ sum under the new parametrization using $J'$, we have
\begin{multline*}
     \frac{\sqrt{\Nm(\mf p)}}{\sqrt{\Nm(\mf l_1\mf l_2)}} \sum\limits_{ \substack{\mf l_1\neq \mf  l_2\\\mf l_1,\mf l_2\in \Lambda}}|x_{\mf l_1}\ov{x_{\mf l_2}}|\sum\limits_{ \substack{J'\lhd \mcO_F\\ \mf N\mid J'\\ \Nm(J')<L\\ [(J')^2/\mf l_1 \mf  l_2]=1}}\sum_{\substack{\frac{(J')^2}{\mf l_1 \mf  l_2}=(\mu)\\|\mu|_{\infty}<L^2}}  \sum_{m,n\in F}\ind_{m\in I}\ind_{n\in I} \ind_{mn\in \mu+I} \left |\wh{\phi}(m,n,\mu)\right|. 
\end{multline*}

 Therefore putting everything together and using the unit counting lemma (lemma \ref{unic}) we get
 $$\ll_{F, \mf N, f_{\infty},\cond(K)} L^{1+\text{o}(1)}\sqrt{\Nm(\mf p)}\sum\limits_{ \substack{\mf l_1\neq \mf  l_2\\\mf l_1,\mf l_2\in \Lambda}}|x_{\mf l_1}\ov{x_{\mf l_2}}|. $$

\begin{prop}[\textbf{Bounds for the non-diagonal contribution for the DFI amplifier}] \label{nondiagbounddfi} 
     Let us assume that the trace function $K$ has Fourier-M\"obius group contained in the standard Borel subgroup. We have the following bound for $L$ and $\Nm(\mf p)$ large enough (as in \ref{obsdfi})
       \begin{multline*}
       \left|\sum_{m,n\in F}\wh{B_2}(m,n)\right|\ll_{f_{\infty}}\cond(K)^4\sqrt{\Nm(\mf p)}\Big( L^{1+\text{o}(1)}\sum\limits_{ \substack{\mf l_1\neq \mf  l_2\\\mf l_1,\mf l_2\in \Lambda\\ \mf l_1,\mf l_2 \text{ prime }}}|x_{\mf l_1}\ov{x_{\mf l_2}}|+
       L^{2+\text{o}(1)}\sum\limits_{ \substack{\mf l_1\neq \mf  l_2\\\mf l_1^2,\mf l_2^2\in \Lambda\\ \mf l_1,\mf l_2 \text{ prime }}}|x_{\mf l_1^2}\ov{x_{\mf l_2^2}}|\\ +L^{3/2+\text{o}(1)}\sum\limits_{ \substack{\mf l_1,\mf l_2^2\in \Lambda\\ \mf l_1,\mf l_2 \text{ prime }}}|x_{\mf l_1}\ov{x_{\mf l_2^2}}|+
       L^{3/2+\text{o}(1)}\sum\limits_{ \substack{\mf l_1^2,\mf l_2\in \Lambda\\ \mf l_1,\mf l_2 \text{ prime }}}|x_{\mf l_1^2}\ov{x_{\mf l_2}}|\Big) 
       \end{multline*}
       where $\cond(K)$ is the smallest conductor of a sheaf $\mcF$ whose trace function is $K$.
\end{prop}
    
\proof
Recall that
$$C_0 . \Nm(\mf p)\ll_{F, \mf N} 1.$$

We have 
$$\Lambda=\{\mf l| \mf l\text{ prime and }\mf l\sim L \text{ coprime to }\mf N\mf p\}\cup \{\mf l^2| \mf l\text{ prime and }\mf l\sim L\text{ coprime to }\mf N\mf p\}.$$

We carry out the proof in six cases according to the elements of $\Lambda$:

\textbf{Case 1:} $\mf l_1$, $\mf l_2$ prime ($\mf l_1\neq \mf l_2$)

\begin{multline*}
\sum_{m,n\in F}|B_2(m,n)|\leq \Bigg|C_0 . \Nm(\mf p)\frac{1}{\sqrt{\Nm(\mf l_1\mf l_2)}}\sum\limits_{ \substack{\mf l_1\neq\mf  l_2\\\mf l_1,\mf l_2\in \Lambda \\ \mf l_1,\mf l_2 \text{ prime }}}x_{\mf l_1}\ov{x_{\mf l_2}}\sum\limits_{\mu \in \mathcal{C}_{\mf l_1,\mf l_2}}\ind_{ c\in \mf N\whO^{(\mf l_1\mf l_2)}} \ind_{ cm \in \whO^{(\mf l_1\mf l_2)}}\ind_{cn \in \whO^{(\mf l_1\mf l_2)}}\\
 \Bigg( \ind_{c\in \mf N\mcO_{\mf l_1}} \ind_{cm \in \mcO_{\mf l_1}}\ind_{cn \in \mcO_{\mf l_1}}
 \ind_{ mn\in \mu+\frac{1}{c}\mcO_{\mf l_1}}\Bigg) \Bigg(\sum_{a=0}^{1} \ind_{\un {l_2}^a c\in \mf N\mcO_{\mf l_2}} \ind_{\un {l_2}^a cm \in \mcO_{\mf l_2}}\ind_{\un {l_2}cn \in \mcO_{\mf l_2}} \ind_{ mn\in \mu+\frac{1}{c\un{l_2}}\mcO_{\mf l_2}}\Bigg)\\ \Cor(K,\gamma_{m,n}(\mu)) \wh{\phi}(m,n,\mu).\Bigg|
\end{multline*}

Recall that $$\mathcal{C}_{\mf l_1,\mf l_2}=F^{\times}\cap (\Aft)^2\frac{\un {l_1}}{\un {l_2}}\whOt$$ and for $\mu\in \mathcal{C}_{\mf l_1,\mf l_2}$, $c\in \Aft$ is any idele satisfying $$c^2\mu\in\frac{\un {l_1}}{\un {l_2}}\whOt.$$

Since $\gamma_{m,n}(\mu)$ is clearly not contained in the standard Borel subgroup and hence in the automorphism group attached to $\wh{K}$, using corollary \ref{corrsumsize}, we have square root cancellation, $$|\Cor(K,\gamma_{m,n}(\mu))|\leq (\cond(K))^4\sqrt{\Nm(\mf p)} .$$

\begin{multline*}
   \sum_{m,n\in F}|B_2(m,n)| \ll(\cond(K))^4 \frac{\sqrt{\Nm(\mf p)}}{\sqrt{\Nm(\mf l_1\mf l_2)}} \sum\limits_{ \substack{\mf l_1\neq \mf  l_2\\\mf l_1,\mf l_2\in \Lambda\\ \mf l_1,\mf l_2 \text{ prime }}}|x_{\mf l_1}\ov{x_{\mf l_2}}|\sum_{a=0}^1\sum_{\mu \in \mathcal{C}_{\mf l_1,\mf l_2}}\ind_{\un {l_2}^a c\in \mf N\whO} \\  \sum_{m,n\in F}\ind_{m\in I_{a}}\ind_{n\in I} \ind_{mn\in \mu+I} \left |\wh{\phi}(m,n,\mu)\right| 
\end{multline*}

where $I:=F\cap \frac{1}{\un {l_2}c}\whO$ and $I_a:=F\cap \frac{1}{\un {l_2}^a c}\whO$ are fractional ideals satisfying $$\Nm(I)= |\un {l_2}c|_{\Af}\gg_{F,f_{\infty}} L^{-1}$$ and $$\Nm(I_a)= |\un {l_2}^a c|_{\Af}\gg_{F,f_{\infty}} L^{-a}$$

We may think of the $\mu$-sum as a sum over ideals $J$ s.t. $\frac{\mf l_1}{\mf l_2}J^2$ is a principal ideal. This is by defining the ideal $J$ to satisfy
$$(\mu)= \frac{\mf l_1}{\mf l_2}.J^2.$$

Moreover we saw in observation \ref{obsdfi} that in the case where $\mf l_1, \mf l_2$ are prime, $\wh{\phi}(m,n,\mu)\neq 0$ only if $|\Nm(J)|\ll_{f,F}1$ The non-archimedean condition implies $J\mf l_2^a$ is integral. Let us write $J'=\mf l_2^a J$ then $J'$ is integral. Since $\mf l_2^a$ is of size $L^a$, we have $\wh{\phi}(m,n,\mu)\neq 0$ only if $|\Nm(J')|\ll_{f,F} L^a$. So the $\mu$-sum is finite and of length $L^a$.

Note that we have
$$I=\frac{J}{\mf l_2}$$
under the new parametrization.

Also recall that 
$$\wh{\phi}(m,n,\mu)=f_{\infty}\begin{pmatrix}
m & \mu-mn \\
1 & -n 
\end{pmatrix}.$$

Hence (since $f_{\infty}$ is compactly supported modulo the center)

$$ |\wh{\phi}(m,n,\mu)| \ll_{f_{\infty}} \ind_{|m|_{\infty}\ll \sqrt{|\mu|_{\infty}}}\ind_{|n|_{\infty}\ll \sqrt{|\mu|_{\infty}}}\ind_{|mn-\mu|_{\infty}\ll \sqrt{|\mu|_{\infty}}}.$$

The inequalities in the indicators also depend on $f_{\infty}$ or more precisely its support.

 Now we bound each term of the $\mu$-sum

$$
    \sum_{m,n\in F}\ind_{m\in I_a}\ind_{n\in I} \ind_{mn\in \mu+I} \left |\wh{\phi}(m,n,\mu)\right|  \ll \sum_{|k-\mu|_{\infty}\ll \sqrt{|\mu|_{\infty}}}  \ind_{k\in \mu +I} \sum_{\substack{m\in I_a\\ n\in I\\ |m|_{\infty}\ll \sqrt{|\mu|_{\infty}}\\|n|_{\infty}\ll \sqrt{|\mu|_{\infty}}}}\ind_{mn=k}.
$$

We bound the innermost sum using the divisor counting lemma (lemma \ref{divc}). Also note that $\Nm(k)\ll \Nm(\mu)$.

$$
\sum_{|k-\mu|_{\infty}\ll \sqrt{|\mu|_{\infty}}}  \ind_{k\in \mu +I} \sum_{\substack{m\in I_a\\ n\in I\\|m|_{\infty}\ll_{f_{\infty}} \sqrt{|\mu|_{\infty}}\\|n|_{\infty}\ll \sqrt{|\mu|_{\infty}}}}\ind_{mn=k} \ll (\Nm(\mu))^{\text{o}(1)} \sum_{|k-\mu|_{\infty}\ll \sqrt{|\mu|_{\infty}}}  \ind_{k\in \mu +I}.
$$
  By the lattice point counting lemma \ref{lpc}  we have
\begin{multline*}
     \sum_{m,n\in F}\ind_{m\in I_a}\ind_{n\in I} \ind_{mn\in \mu+I} \left |\wh{\phi}(m,n,\mu)\right|  \ll (\Nm(\mu))^{\text{o}(1)} \frac{(\Nm(\mu))^{1/2}}{\Nm(I)}\\=(\Nm(\mu))^{\text{o}(1)}\sqrt{\Nm(\mf l_1 \mf l_2)}.
\end{multline*}
 
The last equality follows by substituting the expression for $(\mu)$ and $I$ in terms of $J$ that we noted earlier.

Now rewriting the $\mu$ sum under the new parametrization using $J'$ we have

\begin{multline*}
     \frac{\sqrt{\Nm(\mf p)}}{\sqrt{\Nm(\mf l_1\mf l_2)}} \sum\limits_{ \substack{\mf l_1\neq \mf  l_2\\\mf l_1,\mf l_2\in \Lambda\\ \mf l_1,\mf l_2 \text{ prime }}}|x_{\mf l_1}\ov{x_{\mf l_2}}|\sum_{a=0}^1\sum\limits_{ \substack{J'\lhd \mcO_F\\ \mf N\mid J'\\ \Nm(J')<L^a\\ [(J')^2.\mf l_1 /\mf  l_2^{1+2a}]=1}}\sum_{\substack{(J')^2.\mf l_1 /\mf  l_2^{1+2a}=(\mu)\\|\mu|_{\infty}<L^{2a}}}  \sum_{m,n\in F}\ind_{m\in I_a}\ind_{n\in I} \ind_{mn\in \mu+I} \left |\wh{\phi}(m,n,\mu)\right|. 
\end{multline*}

 Therefore putting everything together and using the unit counting lemma (lemma \ref{unic}) we get 
 $$\ll_{F, \mf N, f_{\infty}} \cond(K)^4 \sqrt{\Nm(\mf p)}\left(\sum_{a=0}^1 L^{a+\text{o}(1)}\right)\sum\limits_{ \substack{\mf l_1\neq \mf  l_2\\\mf l_1,\mf l_2\in \Lambda\\ \mf l_1,\mf l_2 \text{ prime }}}|x_{\mf l_1}\ov{x_{\mf l_2}}|. $$

\textbf{Case 2:} $\mf l_1^2$ square of prime and $\mf l_2^2$ square of prime ($\mf l_1\neq \mf l_2$)

\begin{multline*}
\sum_{m,n\in F}|B_2(m,n)|\leq \Bigg|C_0 . \Nm(\mf p)\frac{1}{\sqrt{\Nm(\mf l_1^2\mf l_2^2)}}\sum\limits_{ \substack{\mf l_1\neq\mf  l_2\\\mf l_1^2,\mf l_2^2\in \Lambda \\ \mf l_1,\mf l_2 \text{ prime }}}x_{\mf l_1^2}\ov{x_{\mf l_2^2}}\sum\limits_{\mu \in \mathcal{C}_{\mf l_1,\mf l_2}}\ind_{ c\in \mf N\whO^{(\mf l_1\mf l_2)}} \ind_{ cm \in \whO^{(\mf l_1\mf l_2)}}\ind_{cn \in \whO^{(\mf l_1\mf l_2)}}\\
 \Bigg( \ind_{c\in \mf N\mcO_{\mf l_1}} \ind_{cm \in \mcO_{\mf l_1}}\ind_{cn \in \mcO_{\mf l_1}}
 \ind_{ mn\in \mu+\frac{1}{c}\mcO_{\mf l_1}}\Bigg) \Bigg(\sum_{a=0}^{2} \ind_{\un {l_2}^a c\in \mf N\mcO_{\mf l_2}} \ind_{\un {l_2}^a cm \in \mcO_{\mf l_2}}\ind_{\un {l_2}^2cn \in \mcO_{\mf l_2}} \ind_{ mn\in \mu+\frac{1}{c\un{l_2}^2}\mcO_{\mf l_2}}\Bigg)\\ \Cor(K,\gamma_{m,n}(\mu)) \wh{\phi}(m,n,\mu).\Bigg|
\end{multline*}

Recall that $$\mathcal{C}_{\mf l_1,\mf l_2}=F^{\times}\cap (\Aft)^2\frac{\un {l_1}^2}{\un {l_2}^2}\whOt$$ and for $\mu\in \mathcal{C}_{\mf l_1,\mf l_2}$, $c\in \Aft$ is any idele satisfying $$c^2\mu\in\frac{\un {l_1}^2}{\un {l_2}^2}\whOt.$$

Since $\gamma_{m,n}(\mu)$ is clearly not contained in the standard Borel subgroup and hence in the automorphism group attached to $\wh{K}$, using corollary \ref{corrsumsize}, we have square root cancellation, $$|\Cor(K,\gamma_{m,n}(\mu))|\leq (\cond(K))^4\sqrt{\Nm(\mf p)} .$$

\begin{multline*}
   \sum_{m,n\in F}|B_2(m,n)| \ll(\cond(K))^4 \frac{\sqrt{\Nm(\mf p)}}{\Nm(\mf l_1\mf l_2)} \sum\limits_{ \substack{\mf l_1\neq \mf  l_2\\\mf l_1^2,\mf l_2^2\in \Lambda\\ \mf l_1,\mf l_2 \text{ prime }}}|x_{\mf l_1^2}\ov{x_{\mf l_2^2}}|\\ \sum_{a=0}^2\sum_{\mu \in \mathcal{C}_{\mf l_1,\mf l_2}}\ind_{\un {l_2}^a c\in \mf N\whO}  \sum_{m,n\in F}\ind_{m\in I_{a}}\ind_{n\in I} \ind_{mn\in \mu+I} \left |\wh{\phi}(m,n,\mu)\right| 
\end{multline*}

where $I:=F\cap \frac{1}{\un {l_2}^2c}\whO$ and $I_a:=F\cap \frac{1}{\un {l_2}^a c}\whO$ are fractional ideals satisfying $$\Nm(I)= |\un {l_2}^2c|_{\Af}\gg_{F,f_{\infty}} L^{-2}$$ and $$\Nm(I_a)= |\un {l_2}^a c|_{\Af}\gg_{F,f_{\infty}} L^{-a}$$

We may think of the $\mu$-sum as a sum over ideals $J$ s.t. $\frac{\mf l_1^2}{\mf l_2^2}J^2$ is a principal ideal. This is by defining the ideal $J$ to satisfy
$$(\mu)= \frac{\mf l_1}{\mf l_2}.J^2.$$

Moreover we saw in observation \ref{obsdfi} that in the case where $\mf l_1^2, \mf l_2^2$ are squares of primes, $\wh{\phi}(m,n,\mu)\neq 0$ only if $|\Nm(J)|\ll_{f,F}1$ The non-archimedean condition implies $J\mf l_2^a$ is integral. Let us write $J'=\mf l_2^a J$ then $J'$ is integral. Since $\mf l_2^a$ is of size $L^a$, we have $\wh{\phi}(m,n,\mu)\neq 0$ only if $|\Nm(J')|\ll_{f,F} L^a$. So the $\mu$-sum is finite and of length $L^a$.

Note that we have
$$I=\frac{J}{\mf l_2^2}$$
under the new parametrization.

Also recall that 
$$\wh{\phi}(m,n,\mu)=f_{\infty}\begin{pmatrix}
m & \mu-mn \\
1 & -n 
\end{pmatrix}.$$

Hence (since $f_{\infty}$ is compactly supported modulo the center)

$$ |\wh{\phi}(m,n,\mu)| \ll_{f_{\infty}} \ind_{|m|_{\infty}\ll \sqrt{|\mu|_{\infty}}}\ind_{|n|_{\infty}\ll \sqrt{|\mu|_{\infty}}}\ind_{|mn-\mu|_{\infty}\ll \sqrt{|\mu|_{\infty}}}.$$

The inequalities in the indicators also depend on $f_{\infty}$ or more precisely its support.

 Now we bound each term of the $\mu$-sum

$$
    \sum_{m,n\in F}\ind_{m\in I_a}\ind_{n\in I} \ind_{mn\in \mu+I} \left |\wh{\phi}(m,n,\mu)\right|  \ll \sum_{|k-\mu|_{\infty}\ll \sqrt{|\mu|_{\infty}}}  \ind_{k\in \mu +I} \sum_{\substack{m\in I_a\\ n\in I\\ |m|_{\infty}\ll \sqrt{|\mu|_{\infty}}\\|n|_{\infty}\ll \sqrt{|\mu|_{\infty}}}}\ind_{mn=k}.
$$

We bound the innermost sum using the divisor counting lemma (lemma \ref{divc}). Also note that $\Nm(k)\ll \Nm(\mu)$.

$$
\sum_{|k-\mu|_{\infty}\ll \sqrt{|\mu|_{\infty}}}  \ind_{k\in \mu +I} \sum_{\substack{m\in I_a\\ n\in I\\|m|_{\infty}\ll_{f_{\infty}} \sqrt{|\mu|_{\infty}}\\|n|_{\infty}\ll \sqrt{|\mu|_{\infty}}}}\ind_{mn=k} \ll (\Nm(\mu))^{\text{o}(1)} \sum_{|k-\mu|_{\infty}\ll \sqrt{|\mu|_{\infty}}}  \ind_{k\in \mu +I}.
$$
  By the lattice point counting lemma \ref{lpc}  we have
\begin{multline*}
     \sum_{m,n\in F}\ind_{m\in I_a}\ind_{n\in I} \ind_{mn\in \mu+I} \left |\wh{\phi}(m,n,\mu)\right|  \ll (\Nm(\mu))^{\text{o}(1)} \frac{(\Nm(\mu))^{1/2}}{\Nm(I)}=(\Nm(\mu))^{\text{o}(1)}\Nm(\mf l_1 \mf l_2).
\end{multline*}
 
The last equality follows by substituting the expression for $(\mu)$ and $I$ in terms of $J$ that we noted earlier.

Now rewriting the $\mu$ sum under the new parametrization using $J'$ we have

\begin{multline*}
     \frac{\sqrt{\Nm(\mf p)}}{\Nm(\mf l_1\mf l_2)} \sum\limits_{ \substack{\mf l_1\neq \mf  l_2\\\mf l_1^2,\mf l_2^2\in \Lambda\\ \mf l_1,\mf l_2 \text{ prime }}}|x_{\mf l_1^2}\ov{x_{\mf l_2^2}}|\sum_{a=0}^2\sum\limits_{ \substack{J'\lhd \mcO_F\\ \mf N\mid J'\\ \Nm(J')<L^a\\ [(J')^2.\mf l_1 /\mf  l_2^{1+2a}]=1}}\sum_{\substack{(J')^2.\mf l_1 /\mf  l_2^{1+2a}=(\mu)\\|\mu|_{\infty}<L^{2a}}}  \sum_{m,n\in F}\ind_{m\in I_a}\ind_{n\in I} \ind_{mn\in \mu+I} \left |\wh{\phi}(m,n,\mu)\right|. 
\end{multline*}

 Therefore putting everything together and using the unit counting lemma (lemma \ref{unic}) we get 
 $$\ll_{F, \mf N, f_{\infty}} \cond(K)^4 \sqrt{\Nm(\mf p)}\left(\sum_{a=0}^2 L^{a+\text{o}(1)}\right)\sum\limits_{ \substack{\mf l_1\neq \mf  l_2\\\mf l_1^2,\mf l_2^2\in \Lambda\\ \mf l_1,\mf l_2 \text{ prime }}}|x_{\mf l_1^2}\ov{x_{\mf l_2^2}}|. $$

\textbf{Case 3:} $\mf l_1^2$ square of prime and $\mf l_2$ prime ($\mf l_1\neq \mf l_2$)

\begin{multline*}
\sum_{m,n\in F}|B_2(m,n)|\leq \Bigg|C_0 . \Nm(\mf p)\frac{1}{\sqrt{\Nm(\mf l_1^2\mf l_2)}}\sum\limits_{ \substack{\mf l_1\neq\mf  l_2\\\mf l_1^2,\mf l_2\in \Lambda \\ \mf l_1,\mf l_2 \text{ prime }}}x_{\mf l_1^2}\ov{x_{\mf l_2}}\sum\limits_{\mu \in \mathcal{C}_{\mf l_1,\mf l_2}}\ind_{ c\in \mf N\whO^{(\mf l_1\mf l_2)}} \ind_{ cm \in \whO^{(\mf l_1\mf l_2)}}\ind_{cn \in \whO^{(\mf l_1\mf l_2)}}\\
 \Bigg( \ind_{c\in \mf N\mcO_{\mf l_1}} \ind_{cm \in \mcO_{\mf l_1}}\ind_{cn \in \mcO_{\mf l_1}}
 \ind_{ mn\in \mu+\frac{1}{c}\mcO_{\mf l_1}}\Bigg) \Bigg(\sum_{a=0}^{1} \ind_{\un {l_2}^a c\in \mf N\mcO_{\mf l_2}} \ind_{\un {l_2}^a cm \in \mcO_{\mf l_2}}\ind_{\un{l_2} cn \in \mcO_{\mf l_2}} \ind_{ mn\in \mu+\frac{1}{c\un{l_2}}\mcO_{\mf l_2}}\Bigg)\\ \Cor(K,\gamma_{m,n}(\mu)) \wh{\phi}(m,n,\mu).\Bigg|
\end{multline*}

Recall that $$\mathcal{C}_{\mf l_1,\mf l_2}=F^{\times}\cap (\Aft)^2\frac{\un {l_1}^2}{\un {l_2}}\whOt$$ and for $\mu\in \mathcal{C}_{\mf l_1,\mf l_2}$, $c\in \Aft$ is any idele satisfying $$c^2\mu\in\frac{\un {l_1}^2}{\un {l_2}}\whOt.$$

Since $\gamma_{m,n}(\mu)$ is clearly not contained in the standard Borel subgroup and hence in the automorphism group attached to $\wh{K}$, using corollary \ref{corrsumsize}, we have square root cancellation, $$|\Cor(K,\gamma_{m,n}(\mu))|\leq (\cond(K))^4\sqrt{\Nm(\mf p)} .$$

\begin{multline*}
   \sum_{m,n\in F}|B_2(m,n)| \ll(\cond(K))^4 \sqrt{\frac{\Nm(\mf p)}{\Nm(\mf l_1^2\mf l_2)}} \sum\limits_{ \substack{\mf l_1\neq \mf  l_2\\\mf l_1^2,\mf l_2\in \Lambda\\ \mf l_1,\mf l_2 \text{ prime }}}|x_{\mf l_1^2}\ov{x_{\mf l_2}}|\\ \sum_{a=0}^1\sum_{\mu \in \mathcal{C}_{\mf l_1,\mf l_2}}\ind_{\un {l_2}^a c\in \mf N\whO}  \sum_{m,n\in F}\ind_{m\in I_{a}}\ind_{n\in I} \ind_{mn\in \mu+I} \left |\wh{\phi}(m,n,\mu)\right| 
\end{multline*}

where $I:=F\cap \frac{1}{\un {l_2}c}\whO$ and $I_a:=F\cap \frac{1}{\un {l_2}^a c}\whO$ are fractional ideals satisfying $$\Nm(I)= |\un {l_2}c|_{\Af}\gg_{F,f_{\infty}} L^{-\frac{3}{2}}$$ and $$\Nm(I_a)= |\un {l_2}^a c|_{\Af}\gg_{F,f_{\infty}} L^{-a-\frac{1}{2}}$$

We may think of the $\mu$-sum as a sum over ideals $J$ s.t. $\frac{\mf l_1^2}{\mf l_2}J^2$ is a principal ideal. This is by defining the ideal $J$ to satisfy
$$(\mu)= \frac{\mf l_1^2}{\mf l_2}.J^2.$$

Moreover we saw in observation \ref{obsdfi} that in this case , $\wh{\phi}(m,n,\mu)\neq 0$ only if $|\Nm(J)|\ll_{f,F}L^{\frac{1}{2}}$ The non-archimedean condition implies $J\mf l_2^a$ is integral. Let us write $J'=\mf l_2^a J$ then $J'$ is integral. Since $\mf l_2^a$ is of size $L^a$, we have $\wh{\phi}(m,n,\mu)\neq 0$ only if $|\Nm(J')|\ll_{f,F} L^{a+\frac{1}{2}}$. So the $\mu$-sum is finite and of length $L^{a+\frac{1}{2}}$.

Note that we have
$$I=\frac{J}{\mf l_2}$$
under the new parametrization.

Also recall that 
$$\wh{\phi}(m,n,\mu)=f_{\infty}\begin{pmatrix}
m & \mu-mn \\
1 & -n 
\end{pmatrix}.$$

Hence (since $f_{\infty}$ is compactly supported modulo the center)

$$ |\wh{\phi}(m,n,\mu)| \ll_{f_{\infty}} \ind_{|m|_{\infty}\ll \sqrt{|\mu|_{\infty}}}\ind_{|n|_{\infty}\ll \sqrt{|\mu|_{\infty}}}\ind_{|mn-\mu|_{\infty}\ll \sqrt{|\mu|_{\infty}}}.$$

The inequalities in the indicators also depend on $f_{\infty}$ or more precisely its support.

 Now we bound each term of the $\mu$-sum

$$
    \sum_{m,n\in F}\ind_{m\in I_a}\ind_{n\in I} \ind_{mn\in \mu+I} \left |\wh{\phi}(m,n,\mu)\right|  \ll \sum_{|k-\mu|_{\infty}\ll \sqrt{|\mu|_{\infty}}}  \ind_{k\in \mu +I} \sum_{\substack{m\in I_a\\ n\in I\\ |m|_{\infty}\ll \sqrt{|\mu|_{\infty}}\\|n|_{\infty}\ll \sqrt{|\mu|_{\infty}}}}\ind_{mn=k}.
$$

We bound the innermost sum using the divisor counting lemma (lemma \ref{divc}). Also note that $\Nm(k)\ll \Nm(\mu)$.

$$
\sum_{|k-\mu|_{\infty}\ll \sqrt{|\mu|_{\infty}}}  \ind_{k\in \mu +I} \sum_{\substack{m\in I_a\\ n\in I\\|m|_{\infty}\ll_{f_{\infty}} \sqrt{|\mu|_{\infty}}\\|n|_{\infty}\ll \sqrt{|\mu|_{\infty}}}}\ind_{mn=k} \ll (\Nm(\mu))^{\text{o}(1)} \sum_{|k-\mu|_{\infty}\ll \sqrt{|\mu|_{\infty}}}  \ind_{k\in \mu +I}.
$$
  By the lattice point counting lemma \ref{lpc}  we have
\begin{multline*}
     \sum_{m,n\in F}\ind_{m\in I_a}\ind_{n\in I} \ind_{mn\in \mu+I} \left |\wh{\phi}(m,n,\mu)\right|  \ll (\Nm(\mu))^{\text{o}(1)} \frac{(\Nm(\mu))^{1/2}}{\Nm(I)}=(\Nm(\mu))^{\text{o}(1)}\Nm(\mf l_1)\sqrt{\Nm(\mf l_2)}.
\end{multline*}
 
The last equality follows by substituting the expression for $(\mu)$ and $I$ in terms of $J$ that we noted earlier.

Now rewriting the $\mu$ sum under the new parametrization using $J'$ we have

\begin{multline*}
     \sqrt{\frac{\Nm(\mf p)}{\Nm(\mf l_1^2\mf l_2)}} \sum\limits_{ \substack{\mf l_1\neq \mf  l_2\\\mf l_1^2,\mf l_2\in \Lambda\\ \mf l_1,\mf l_2 \text{ prime }}}|x_{\mf l_1^2}\ov{x_{\mf l_2}}|\sum_{a=0}^1\sum\limits_{ \substack{J'\lhd \mcO_F\\ \mf N\mid J'\\ \Nm(J')<L^a\\ [(J')^2.\mf l_1 /\mf  l_2^{1+2a}]=1}}\sum_{\substack{(J')^2.\mf l_1 /\mf  l_2^{1+2a}=(\mu)\\|\mu|_{\infty}<L^{2a+1}}}  \sum_{m,n\in F}\ind_{m\in I_a}\ind_{n\in I} \ind_{mn\in \mu+I} \left |\wh{\phi}(m,n,\mu)\right|. 
\end{multline*}

 Therefore putting everything together and using the unit counting lemma (lemma \ref{unic}) we get 
 $$\ll_{F, \mf N, f_{\infty}} \cond(K)^4 \sqrt{\Nm(\mf p)}\left(\sum_{a=0}^1 L^{a+\frac{1}{2}+\text{o}(1)}\right)\sum\limits_{ \substack{\mf l_1\neq \mf  l_2\\\mf l_1^2,\mf l_2\in \Lambda\\ \mf l_1,\mf l_2 \text{ prime }}}|x_{\mf l_1^2}\ov{x_{\mf l_2}}|. $$

\textbf{Case 4:} $\mf l_1^2$ square of prime and $\mf l_2$ prime ($\mf l_1=\mf l_2$)

\begin{multline*}
\sum_{m,n\in F}|B_2(m,n)|\leq \Bigg|C_0 . \Nm(\mf p)\frac{1}{\sqrt{\Nm(\mf l_1^3)}}\sum\limits_{ \substack{\mf l_1^2,\mf l_1\in \Lambda \\ \mf l_1, \text{ prime }}}x_{\mf l_1^2}\ov{x_{\mf l_1}}\sum\limits_{\mu \in \mathcal{C}_{\mf l_1,\mf l_2}}\ind_{ c\in \mf N\whO^{(\mf l_1)}} \ind_{ cm \in \whO^{(\mf l_1)}}\ind_{cn \in \whO^{(\mf l_1)}}\\
 \Bigg( \ind_{c\in \mf N\mcO_{\mf l_1}} \ind_{cm \in \mcO_{\mf l_1}}\ind_{cn \in \mcO_{\mf l_1}}
 \ind_{ mn\in \mu+\frac{1}{c}\mcO_{\mf l_1}}\Bigg) \ \Cor(K,\gamma_{m,n}(\mu)) \wh{\phi}(m,n,\mu).\Bigg|
\end{multline*}

Recall that $$\mathcal{C}_{\mf l_1,\mf l_2}=F^{\times}\cap (\Aft)^2\un {l_1}\whOt$$ and for $\mu\in \mathcal{C}_{\mf l_1,\mf l_2}$, $c\in \Aft$ is any idele satisfying $$c^2\mu\in \un {l_1}\whOt.$$

Since $\gamma_{m,n}(\mu)$ is clearly not contained in the standard Borel subgroup and hence in the automorphism group attached to $\wh{K}$, using corollary \ref{corrsumsize}, we have square root cancellation, $$|\Cor(K,\gamma_{m,n}(\mu))|\leq (\cond(K))^4\sqrt{\Nm(\mf p)} .$$

\begin{multline*}
   \sum_{m,n\in F}|B_2(m,n)| \ll(\cond(K))^4 \sqrt{\frac{\Nm(\mf p)}{\Nm(\mf l_1^3)}} \sum\limits_{ \substack{\mf l_1^2,\mf l_1\in \Lambda\\ \mf l_1 \text{ prime }}}|x_{\mf l_1^2}\ov{x_{\mf l_1}}|\\ \sum_{a=0}^1\sum_{\mu \in \mathcal{C}_{\mf l_1}}\ind_{\un {l_1}^a c\in \mf N\whO}  \sum_{m,n\in F}\ind_{m\in I_{a}}\ind_{n\in I} \ind_{mn\in \mu+I} \left |\wh{\phi}(m,n,\mu)\right| 
\end{multline*}

where $I:=F\cap \frac{1}{\un {l_1}c}\whO$ and $I_a:=F\cap \frac{1}{\un {l_1}^a c}\whO$ are fractional ideals satisfying $$\Nm(I)= |\un {l_1}c|_{\Af}\gg_{F,f_{\infty}} L^{-\frac{3}{2}}$$ and $$\Nm(I_a)= |\un {l_1}^a c|_{\Af}\gg_{F,f_{\infty}} L^{-a-\frac{1}{2}}$$

We may think of the $\mu$-sum as a sum over ideals $J$ s.t. $\mf l_1 J^2$ is a principal ideal. This is by defining the ideal $J$ to satisfy
$$(\mu)= \mf l_1.J^2.$$

Moreover we saw in observation \ref{obsdfi} that in this case , $\wh{\phi}(m,n,\mu)\neq 0$ only if $|\Nm(J)|\ll_{f,F}L^{\frac{1}{2}}$ The non-archimedean condition implies $J\mf l_2^a$ is integral. Let us write $J'=\mf l_1^a J$ then $J'$ is integral. Since $\mf l_1^a$ is of size $L^a$, we have $\wh{\phi}(m,n,\mu)\neq 0$ only if $|\Nm(J')|\ll_{f,F} L^{a+\frac{1}{2}}$. So the $\mu$-sum is finite and of length $L^{a+\frac{1}{2}}$.

Note that we have
$$I=\frac{J}{\mf l_1}$$
under the new parametrization.

Also recall that 
$$\wh{\phi}(m,n,\mu)=f_{\infty}\begin{pmatrix}
m & \mu-mn \\
1 & -n 
\end{pmatrix}.$$

Hence (since $f_{\infty}$ is compactly supported modulo the center)

$$ |\wh{\phi}(m,n,\mu)| \ll_{f_{\infty}} \ind_{|m|_{\infty}\ll \sqrt{|\mu|_{\infty}}}\ind_{|n|_{\infty}\ll \sqrt{|\mu|_{\infty}}}\ind_{|mn-\mu|_{\infty}\ll \sqrt{|\mu|_{\infty}}}.$$

The inequalities in the indicators also depend on $f_{\infty}$ or more precisely its support.

 Now we bound each term of the $\mu$-sum

$$
    \sum_{m,n\in F}\ind_{m\in I_a}\ind_{n\in I} \ind_{mn\in \mu+I} \left |\wh{\phi}(m,n,\mu)\right|  \ll \sum_{|k-\mu|_{\infty}\ll \sqrt{|\mu|_{\infty}}}  \ind_{k\in \mu +I} \sum_{\substack{m\in I_a\\ n\in I\\ |m|_{\infty}\ll \sqrt{|\mu|_{\infty}}\\|n|_{\infty}\ll \sqrt{|\mu|_{\infty}}}}\ind_{mn=k}.
$$

We bound the innermost sum using the divisor counting lemma (lemma \ref{divc}). Also note that $\Nm(k)\ll \Nm(\mu)$.

$$
\sum_{|k-\mu|_{\infty}\ll \sqrt{|\mu|_{\infty}}}  \ind_{k\in \mu +I} \sum_{\substack{m\in I_a\\ n\in I\\|m|_{\infty}\ll_{f_{\infty}} \sqrt{|\mu|_{\infty}}\\|n|_{\infty}\ll \sqrt{|\mu|_{\infty}}}}\ind_{mn=k} \ll (\Nm(\mu))^{\text{o}(1)} \sum_{|k-\mu|_{\infty}\ll \sqrt{|\mu|_{\infty}}}  \ind_{k\in \mu +I}.
$$
  By the lattice point counting lemma \ref{lpc}  we have
\begin{multline*}
     \sum_{m,n\in F}\ind_{m\in I_a}\ind_{n\in I} \ind_{mn\in \mu+I} \left |\wh{\phi}(m,n,\mu)\right|  \ll (\Nm(\mu))^{\text{o}(1)} \frac{(\Nm(\mu))^{1/2}}{\Nm(I)}=(\Nm(\mu))^{\text{o}(1)}\Nm(\mf l_1)^{\frac{3}{2}}.
\end{multline*}
 
The last equality follows by substituting the expression for $(\mu)$ and $I$ in terms of $J$ that we noted earlier.

Now rewriting the $\mu$ sum under the new parametrization using $J'$ we have

\begin{multline*}
     \sqrt{\frac{\Nm(\mf p)}{\Nm(\mf l_1^3)}} \sum\limits_{ \substack{\mf l_1^2,\mf l_1\in \Lambda\\ \mf l_1 \text{ prime }}}|x_{\mf l_1^2}\ov{x_{\mf l_1}}|\sum_{a=0}^1\sum\limits_{ \substack{J'\lhd \mcO_F\\ \mf N\mid J'\\ \Nm(J')<L^a\\ [(J')^2. /\mf  l_1^{2a}]=1}}\sum_{\substack{(J')^2. /\mf  l_1^{2a}=(\mu)\\|\mu|_{\infty}<L^{2a+1}}}  \sum_{m,n\in F}\ind_{m\in I_a}\ind_{n\in I} \ind_{mn\in \mu+I} \left |\wh{\phi}(m,n,\mu)\right|. 
\end{multline*}

 Therefore putting everything together and using the unit counting lemma (lemma \ref{unic}) we get 
 $$\ll_{F, \mf N, f_{\infty}} \cond(K)^4 \sqrt{\Nm(\mf p)}\left(\sum_{a=0}^1 L^{a+\frac{1}{2}+\text{o}(1)}\right)\sum\limits_{ \substack{\mf l_1^2,\mf l_1\in \Lambda\\ \mf l_1 \text{ prime }}}|x_{\mf l_1^2}\ov{x_{\mf l_1}}|. $$

\textbf{Case 5:} $\mf l_1$  prime and $\mf l_2^2$ square of prime ($\mf l_1\neq \mf l_2$)

\begin{multline*}
\sum_{m,n\in F}|B_2(m,n)|\leq \Bigg|C_0 . \Nm(\mf p)\frac{1}{\sqrt{\Nm(\mf l_1\mf l_2^2)}}\sum\limits_{ \substack{\mf l_1\neq\mf  l_2\\\mf l_1,\mf l_2^2\in \Lambda \\ \mf l_1,\mf l_2 \text{ prime }}}x_{\mf l_1}\ov{x_{\mf l_2}^2}\sum\limits_{\mu \in \mathcal{C}_{\mf l_1,\mf l_2}}\ind_{ c\in \mf N\whO^{(\mf l_1\mf l_2)}} \ind_{ cm \in \whO^{(\mf l_1\mf l_2)}}\ind_{cn \in \whO^{(\mf l_1\mf l_2)}}\\
 \Bigg( \ind_{c\in \mf N\mcO_{\mf l_1}} \ind_{cm \in \mcO_{\mf l_1}}\ind_{cn \in \mcO_{\mf l_1}}
 \ind_{ mn \in \mu+\frac{1}{c}\mcO_{\mf l_1}}\Bigg) \Bigg(\sum_{a=0}^{2} \ind_{\un {l_2}^a c\in \mf N\mcO_{\mf l_2}} \ind_{\un {l_2}^a cm \in \mcO_{\mf l_2}}\ind_{\un{l_2}^2 cn \in \mcO_{\mf l_2}} \ind_{ mn\in \mu+\frac{1}{c\un{l_2}^2}\mcO_{\mf l_2}}\Bigg)\\ \Cor(K,\gamma_{m,n}(\mu)) \wh{\phi}(m,n,\mu).\Bigg|
\end{multline*}

Recall that $$\mathcal{C}_{\mf l_1,\mf l_2}=F^{\times}\cap (\Aft)^2\frac{\un {l_1}}{\un {l_2}^2}\whOt$$ and for $\mu\in \mathcal{C}_{\mf l_1,\mf l_2}$, $c\in \Aft$ is any idele satisfying $$c^2\mu\in\frac{\un {l_1}}{\un {l_2}^2}\whOt.$$

Since $\gamma_{m,n}(\mu)$ is clearly not contained in the standard Borel subgroup and hence in the automorphism group attached to $\wh{K}$, using corollary \ref{corrsumsize}, we have square root cancellation, $$|\Cor(K,\gamma_{m,n}(\mu))|\leq (\cond(K))^4\sqrt{\Nm(\mf p)} .$$

\begin{multline*}
   \sum_{m,n\in F}|B_2(m,n)| \ll(\cond(K))^4 \sqrt{\frac{\Nm(\mf p)}{\Nm(\mf l_1\mf l_2^2)}} \sum\limits_{ \substack{\mf l_1\neq \mf  l_2\\\mf l_1,\mf l_2^2\in \Lambda\\ \mf l_1,\mf l_2 \text{ prime }}}|x_{\mf l_1}\ov{x_{\mf l_2^2}}|\\ \sum_{a=0}^2\sum_{\mu \in \mathcal{C}_{\mf l_1,\mf l_2}}\ind_{\un {l_2}^a c\in \mf N\whO}  \sum_{m,n\in F}\ind_{m\in I_{a}}\ind_{n\in I} \ind_{mn\in \mu+I} \left |\wh{\phi}(m,n,\mu)\right| 
\end{multline*}

where $I:=F\cap \frac{1}{\un {l_2}^2c}\whO$ and $I_a:=F\cap \frac{1}{\un {l_2}^a c}\whO$ are fractional ideals satisfying $$\Nm(I)= |\un {l_2}^2c|_{\Af}\gg_{F,f_{\infty}} L^{-\frac{3}{2}}$$ and $$\Nm(I_a)= |\un {l_2}^a c|_{\Af}\gg_{F,f_{\infty}} L^{-a+\frac{1}{2}}$$

We may think of the $\mu$-sum as a sum over ideals $J$ s.t. $\frac{\mf l_1}{\mf l_2^2}J^2$ is a principal ideal. This is by defining the ideal $J$ to satisfy
$$(\mu)= \frac{\mf l_1}{\mf l_2^2}.J^2.$$

Moreover we saw in observation \ref{obsdfi} that in this case , $\wh{\phi}(m,n,\mu)\neq 0$ only if $|\Nm(J)|\ll_{f,F}L^{-\frac{1}{2}}$ The non-archimedean condition implies $J\mf l_2^a$ is integral. Let us write $J'=\mf l_2^a J$ then $J'$ is integral. Since $\mf l_2^a$ is of size $L^a$, we have $\wh{\phi}(m,n,\mu)\neq 0$ only if $|\Nm(J')|\ll_{f,F} L^{a-\frac{1}{2}}$. So the $\mu$-sum is finite and of length $L^{a-\frac{1}{2}}$. (The case $a=0$ yields the empty sum.)

Note that we have
$$I=\frac{J}{\mf l_2^2}$$
under the new parametrization.

Also recall that 
$$\wh{\phi}(m,n,\mu)=f_{\infty}\begin{pmatrix}
m & \mu-mn \\
1 & -n 
\end{pmatrix}.$$

Hence (since $f_{\infty}$ is compactly supported modulo the center)

$$ |\wh{\phi}(m,n,\mu)| \ll_{f_{\infty}} \ind_{|m|_{\infty}\ll \sqrt{|\mu|_{\infty}}}\ind_{|n|_{\infty}\ll \sqrt{|\mu|_{\infty}}}\ind_{|mn-\mu|_{\infty}\ll \sqrt{|\mu|_{\infty}}}.$$

The inequalities in the indicators also depend on $f_{\infty}$ or more precisely its support.

 Now we bound each term of the $\mu$-sum

$$
    \sum_{m,n\in F}\ind_{m\in I_a}\ind_{n\in I} \ind_{mn\in \mu+I} \left |\wh{\phi}(m,n,\mu)\right|  \ll \sum_{|k-\mu|_{\infty}\ll \sqrt{|\mu|_{\infty}}}  \ind_{k\in \mu +I} \sum_{\substack{m\in I_a\\ n\in I\\ |m|_{\infty}\ll \sqrt{|\mu|_{\infty}}\\|n|_{\infty}\ll \sqrt{|\mu|_{\infty}}}}\ind_{mn=k}.
$$

We bound the innermost sum using the divisor counting lemma (lemma \ref{divc}). Also note that $\Nm(k)\ll \Nm(\mu)$.

$$
\sum_{|k-\mu|_{\infty}\ll \sqrt{|\mu|_{\infty}}}  \ind_{k\in \mu +I} \sum_{\substack{m\in I_a\\ n\in I\\|m|_{\infty}\ll_{f_{\infty}} \sqrt{|\mu|_{\infty}}\\|n|_{\infty}\ll \sqrt{|\mu|_{\infty}}}}\ind_{mn=k} \ll (\Nm(\mu))^{\text{o}(1)} \sum_{|k-\mu|_{\infty}\ll \sqrt{|\mu|_{\infty}}}  \ind_{k\in \mu +I}.
$$
  By the lattice point counting lemma \ref{lpc}  we have
\begin{multline*}
     \sum_{m,n\in F}\ind_{m\in I_a}\ind_{n\in I} \ind_{mn\in \mu+I} \left |\wh{\phi}(m,n,\mu)\right|  \ll (\Nm(\mu))^{\text{o}(1)} \frac{(\Nm(\mu))^{1/2}}{\Nm(I)}=(\Nm(\mu))^{\text{o}(1)}\sqrt{\Nm(\mf l_1)\Nm(\mf l_2)^2}.
\end{multline*}
 
The last equality follows by substituting the expression for $(\mu)$ and $I$ in terms of $J$ that we noted earlier.

Now rewriting the $\mu$ sum under the new parametrization using $J'$ we have

\begin{multline*}
     \sqrt{\frac{\Nm(\mf p)}{\Nm(\mf l_1 \mf l_2^2)}} \sum\limits_{ \substack{\mf l_1\neq \mf  l_2\\\mf l_1,\mf l_2^2\in \Lambda\\ \mf l_1,\mf l_2 \text{ prime }}}|x_{\mf l_1}\ov{x_{\mf l_2^2}}|\sum_{a=1}^2\sum\limits_{ \substack{J'\lhd \mcO_F\\ \mf N\mid J'\\ \Nm(J')<L^a\\ [(J')^2.\mf l_1 /\mf  l_2^{2+2a}]=1}}\sum_{\substack{(J')^2.\mf l_1 /\mf  l_2^{2+2a}=(\mu)\\|\mu|_{\infty}<L^{2a-1}}}  \sum_{m,n\in F}\ind_{m\in I_a}\ind_{n\in I} \ind_{mn\in \mu+I} \left |\wh{\phi}(m,n,\mu)\right|. 
\end{multline*}

 Therefore putting everything together and using the unit counting lemma (lemma \ref{unic}) we get 
 $$\ll_{F, \mf N, f_{\infty}} \cond(K)^4 \sqrt{\Nm(\mf p)}\left(\sum_{a=1}^2 L^{a-\frac{1}{2}+\text{o}(1)}\right)\sum\limits_{ \substack{\mf l_1\neq \mf  l_2\\\mf l_1,\mf l_2^2\in \Lambda\\ \mf l_1,\mf l_2 \text{ prime }}}|x_{\mf l_1}\ov{x_{\mf l_2^2}}|. $$

\textbf{Case 6:} $\mf l_1$  prime and $\mf l_2^2$ square of prime ($\mf l_1= \mf l_2$)

\begin{multline*}
\sum_{m,n\in F}|B_2(m,n)|\leq \Bigg|C_0 . \Nm(\mf p)\frac{1}{\sqrt{\Nm(\mf l_1^3)}}\sum\limits_{ \substack{\mf l_1,\mf l_1^2\in \Lambda \\ \mf l_1,\mf l_2 \text{ prime }}}x_{\mf l_1}\ov{x_{\mf l_1}^2}\sum\limits_{\mu \in \mathcal{C}_{\mf l_1}}\ind_{ c\in \mf N\whO^{(\mf l_1)}} \ind_{ cm \in \whO^{(\mf l_1)}}\ind_{cn \in \whO^{(\mf l_1)}}\\
\Bigg(\sum_{a=0}^{2} \ind_{\un {l_1}^a c\in \mf N\mcO_{\mf l_1}} \ind_{\un {l_1}^a cm \in \mcO_{\mf l_1}}\ind_{\un{l_1}^2 cn \in \mcO_{\mf l_2}} \ind_{ mn\in \mu+\frac{1}{c\un{l_1}^2}\mcO_{\mf l_1}}\Bigg) \Cor(K,\gamma_{m,n}(\mu)) \wh{\phi}(m,n,\mu).\Bigg|
\end{multline*}

Recall that $$\mathcal{C}_{\mf l_1}=F^{\times}\cap (\Aft)^2\frac{1}{\un {l_2}}\whOt$$ and for $\mu\in \mathcal{C}_{\mf l_1}$, $c\in \Aft$ is any idele satisfying $$c^2\mu\in \un {l_2}\whOt.$$

Since $\gamma_{m,n}(\mu)$ is clearly not contained in the standard Borel subgroup and hence in the automorphism group attached to $\wh{K}$, using corollary \ref{corrsumsize}, we have square root cancellation, $$|\Cor(K,\gamma_{m,n}(\mu))|\leq (\cond(K))^4\sqrt{\Nm(\mf p)} .$$

\begin{multline*}
   \sum_{m,n\in F}|B_2(m,n)| \ll(\cond(K))^4 \sqrt{\frac{\Nm(\mf p)}{\Nm(\mf l_1^3)}} \sum\limits_{ \substack{\mf l_1,\mf l_1^2\in \Lambda\\ \mf l_1 \text{ prime }}}|x_{\mf l_1}\ov{x_{\mf l_1^2}}|\\ \sum_{a=0}^2\sum_{\mu \in \mathcal{C}_{\mf l_1}}\ind_{\un {l_1}^a c\in \mf N\whO}  \sum_{m,n\in F}\ind_{m\in I_{a}}\ind_{n\in I} \ind_{mn\in \mu+I} \left |\wh{\phi}(m,n,\mu)\right| 
\end{multline*}

where $I:=F\cap \frac{1}{\un {l_1}^2c}\whO$ and $I_a:=F\cap \frac{1}{\un {l_1}^a c}\whO$ are fractional ideals satisfying $$\Nm(I)= |\un {l_1}^2c|_{\Af}\gg_{F,f_{\infty}} L^{-\frac{3}{2}}$$ and $$\Nm(I_a)= |\un {l_1}^a c|_{\Af}\gg_{F,f_{\infty}} L^{-a+\frac{1}{2}}$$

We may think of the $\mu$-sum as a sum over ideals $J$ s.t. $\frac{1}{\mf l_1}J^2$ is a principal ideal. This is by defining the ideal $J$ to satisfy
$$(\mu)= \frac{1}{\mf l_1}.J^2.$$

Moreover we saw in observation \ref{obsdfi} that in this case , $\wh{\phi}(m,n,\mu)\neq 0$ only if $|\Nm(J)|\ll_{f,F}L^{-\frac{1}{2}}$ The non-archimedean condition implies $J\mf l_1^a$ is integral. Let us write $J'=\mf l_1^a J$ then $J'$ is integral. Since $\mf l_1^a$ is of size $L^a$, we have $\wh{\phi}(m,n,\mu)\neq 0$ only if $|\Nm(J')|\ll_{f,F} L^{a-\frac{1}{2}}$. So the $\mu$-sum is finite and of length $L^{a-\frac{1}{2}}$. (The case $a=0$ yields the empty sum.)

Note that we have
$$I=\frac{J}{\mf l_1^2}$$
under the new parametrization.

Also recall that 
$$\wh{\phi}(m,n,\mu)=f_{\infty}\begin{pmatrix}
m & \mu-mn \\
1 & -n 
\end{pmatrix}.$$

Hence (since $f_{\infty}$ is compactly supported modulo the center)

$$ |\wh{\phi}(m,n,\mu)| \ll_{f_{\infty}} \ind_{|m|_{\infty}\ll \sqrt{|\mu|_{\infty}}}\ind_{|n|_{\infty}\ll \sqrt{|\mu|_{\infty}}}\ind_{|mn-\mu|_{\infty}\ll \sqrt{|\mu|_{\infty}}}.$$

The inequalities in the indicators also depend on $f_{\infty}$ or more precisely its support.

 Now we bound each term of the $\mu$-sum

$$
    \sum_{m,n\in F}\ind_{m\in I_a}\ind_{n\in I} \ind_{mn\in \mu+I} \left |\wh{\phi}(m,n,\mu)\right|  \ll \sum_{|k-\mu|_{\infty}\ll \sqrt{|\mu|_{\infty}}}  \ind_{k\in \mu +I} \sum_{\substack{m\in I_a\\ n\in I\\ |m|_{\infty}\ll \sqrt{|\mu|_{\infty}}\\|n|_{\infty}\ll \sqrt{|\mu|_{\infty}}}}\ind_{mn=k}.
$$

We bound the innermost sum using the divisor counting lemma (lemma \ref{divc}). Also note that $\Nm(k)\ll \Nm(\mu)$.

$$
\sum_{|k-\mu|_{\infty}\ll \sqrt{|\mu|_{\infty}}}  \ind_{k\in \mu +I} \sum_{\substack{m\in I_a\\ n\in I\\|m|_{\infty}\ll_{f_{\infty}} \sqrt{|\mu|_{\infty}}\\|n|_{\infty}\ll \sqrt{|\mu|_{\infty}}}}\ind_{mn=k} \ll (\Nm(\mu))^{\text{o}(1)} \sum_{|k-\mu|_{\infty}\ll \sqrt{|\mu|_{\infty}}}  \ind_{k\in \mu +I}.
$$
  By the lattice point counting lemma \ref{lpc}  we have
\begin{multline*}
     \sum_{m,n\in F}\ind_{m\in I_a}\ind_{n\in I} \ind_{mn\in \mu+I} \left |\wh{\phi}(m,n,\mu)\right|  \ll (\Nm(\mu))^{\text{o}(1)} \frac{(\Nm(\mu))^{1/2}}{\Nm(I)}=(\Nm(\mu))^{\text{o}(1)}\sqrt{\Nm(\mf l_1)^3}.
\end{multline*}
 
The last equality follows by substituting the expression for $(\mu)$ and $I$ in terms of $J$ that we noted earlier.

Now rewriting the $\mu$ sum under the new parametrization using $J'$ we have

\begin{multline*}
     \sqrt{\frac{\Nm(\mf p)}{\Nm(\mf l_1^3)}} \sum\limits_{ \substack{\mf l_1\neq \mf  l_2\\\mf l_1,\mf l_1^2\in \Lambda\\ \mf l_1\text{ prime }}}|x_{\mf l_1}\ov{x_{\mf l_1^2}}|\sum_{a=1}^2\sum\limits_{ \substack{J'\lhd \mcO_F\\ \mf N\mid J'\\ \Nm(J')<L^a\\ [(J')^2 1/\mf  l_2^{1+2a}]=1}}\sum_{\substack{(J')^2.1 /\mf  l_2^{1+2a}=(\mu)\\|\mu|_{\infty}<L^{2a-1}}}  \sum_{m,n\in F}\ind_{m\in I_a}\ind_{n\in I} \ind_{mn\in \mu+I} \left |\wh{\phi}(m,n,\mu)\right|. 
\end{multline*}

 Therefore putting everything together and using the unit counting lemma (lemma \ref{unic}) we get 
 $$\ll_{F, \mf N, f_{\infty}} \cond(K)^4 \sqrt{\Nm(\mf p)}\left(\sum_{a=1}^2 L^{a-\frac{1}{2}+\text{o}(1)}\right)\sum\limits_{ \substack{\mf l_1,\mf l_1^2\in \Lambda\\ \mf l_1 \text{ prime }}}|x_{\mf l_1}\ov{x_{\mf l_1^2}}|. $$
\qed

\printbibliography[title={References}]
\end{document}